\documentclass[12pt]{article}
\usepackage{amsmath}
\usepackage[english]{babel}
\usepackage[psamsfonts]{amsfonts}
\usepackage{theorem}
\usepackage{version}
\usepackage{graphicx}
\usepackage{pstricks,pst-node,pst-text,pst-3d}
\usepackage{subfigure}
\usepackage{pdflscape}
\usepackage{lscape}
\usepackage{multirow}
\usepackage{setspace}
\usepackage{url} 
\usepackage[toc,page]{appendix}

\newtheorem{Lem}{Lemma}[section]

\newtheorem{defin}{Definition}[section]

\newtheorem{theor}{Theorem}[section]

\newtheorem{rem}{Remark}[section]

\newcommand{\R}{\mathbb R}

\newcommand{\Z}{\mathbb{Z}}

\newcommand{\Rb}{\mathbf{R}}

\newcommand{\Hb}{\ensuremath{\mathbf{H}}}

\newcommand{\Nb}{\ensuremath{\mathbf{N}}}

\newcommand{\eps}{\epsilon}

\DeclareMathOperator \vecc {vec}

\DeclareMathOperator \tr {tr}

\numberwithin{equation}{section}

\newcommand{\twenty}{\textwidth 16.2cm}

\onehalfspace

\begin{document}

\title {{\sc Asymptotic properties of QML estimators for VARMA models with time-dependent coefficients: Part I}}

%\thedate
\author{
%Abdelkamel {\sc Alj} \thanks{Universit\'{e} Moulay Ismail, FSJES, B.P. 3102 Toulal, Mekn\`{e}s, Morocco. This work was done while at Dept of Mathematics, Universit\'{e} libre de Bruxelles CP112, Blv du Triomphe, B-1050 Bruxelles, Belgium (e-mail: abdelkamel.alj@gmail.com).}, %ALM_V13
Abdelkamel {\sc Alj} \thanks{Universit\'{e} Moulay Ismail, Facult\'{e} des Sciences juridiques, \'{e}conomiques et sociales, B.P. 3102 Toulal, Mekn\`{e}s, Morocco. This work was started while at Universit\'{e} libre de Bruxelles, Faculty of Sciences, Dept of Mathematics, Blv du Triomphe CP 210, B-1050 Brussels, Belgium (e-mail: abdelkamel.alj@gmail.com).}, %ALM_V13
Christophe {\sc Ley}\thanks{Universit\'{e} libre de Bruxelles, Faculty of Sciences, Dept of Mathematics, Blv du Triomphe CP 210, B-1050 Brussels, Belgium (e-mail: chrisley@ulb.ac.be).} \,
{\rm and}
Guy {\sc M\'{e}lard}\thanks{Universit\'{e} libre de Bruxelles, SBS-EM, ECARES,  avenue Franklin Roosevelt 50 CP 114/04, B-1050 Brussels, Belgium (e-mail: gmelard@ulb.ac.be).} }

\maketitle
\vspace{-10mm}
\begin{abstract}\twenty
%This paper is about vector auto-regressive-moving average (VARMA) models with time-dependent coefficients to represent non-stationary time series. %ALM_V16
This paper is about vector autoregressive-moving average (VARMA) models with time-dependent coefficients to represent non-stationary time series. %ALM_V16
%Contrarily to Dahlhaus (2000) and some parts of Azrak and M\'{e}lard (2006) in the %ALM_V13
Contrarily to other papers in the %ALM_V13
 univariate case, the coefficients depend on time but not on the length of the series $n$. 
Under appropriate assumptions, it is shown that a Gaussian quasi-maximum likelihood estimator is almost surely consistent and asymptotically normal. 
The theoretical results are illustrated by means of two examples of bivariate %processes, generalizing Kwoun and Yajima (1986).  %ALM_V13
processes.  %ALM_V13
It is shown that the assumptions underlying the theoretical results apply. 
In the second example the innovations are also marginally heteroscedastic with a correlation ranging from $-0.8$ to $0.8$. 
In the two examples, the asymptotic information matrix is obtained in the Gaussian case. 
Finally, the finite-sample behaviour is checked via a Monte Carlo simulation study for $n$ going from $25$ to $400$.
The results confirm the validity of the asymptotic properties even for short series and reveal that the asymptotic information matrix deduced from the theory is correct. 
\end{abstract}\twenty


Running title: Asymptotics of QMLEs for tdVARMA models\twenty

 %, Monte-Carlo experiments.
\twenty
\twenty
\newpage
%%%%%%%%*********************************************************************************%%%%%%%%%%%%%%%%%%%%%%%%%%%%
%%%%%%%%%%%%%%%%%%%%%%%%%%%%%%%%%%%%%%%%%%%%%%%%%%%%%%%%%%%%%%%%%%%%%%%
%%%%%%%%%%%%%%%%%%%%%%%%%%%%%%%%%%%%%%%%%%%%%%%%%%%%%%%%%%%%%%%%%%%%%%%
%%%%%%%%%%%%%%%%%%%%%%%%%%%%%%%%%%%%%%%%%%%%%%%%%%%%%%%%%%%%%%%%%%%%%%%
\section{Introduction}\label{S.1}

%{\red GM: Trop de références !}
A large part of the literature on time series models is concerned with stationary models. This is of course due to the ensuing mathematical simplifications of stationarity. %GM 24/02/15
Even in that simple context, an asymptotic analysis is not necessarily easy; %indeed, the celebrated auto-regressive-moving average (ARMA) models popularized by %ALM_V16
indeed, the celebrated autoregressive-moving average (ARMA) models popularized by %ALM_V16
 Box and Jenkins (Box {et al.}, 2008) require several pages in Brockwell \& Davis (1991, pp. 375-396) for the derivation of their asymptotic properties. 
However,  the assumption of  invariance over time (especially for long time intervals) is difficult to justify in most practical situations. Therefore,  recent years have seen an increasing interest in models with time-dependent coefficients and non-stationary time series. 
%Analyzed in the seminal work Quenouille (1957), models with time-dependent coefficients for univariate time series %GM 24/02/15 %ATM_V13
Initiated by the seminal work  Quenouille (1957), models with time-dependent or time-varying coefficients 
for univariate time series %GM 24/02/15 %ATM_V13
have been investigated over the years by, \emph{inter alia}, Whittle (1965), Subba Rao (1970), %Hallin and M\'elard (1977), 
Tj\o stheim (1984), %GM 24/02/15
Kwoun \& Yajima (1986), %GM 24/02/15
Singh \& Peiris (1987), 
Priestley (1988), %GM 24/02/15
Grillenzoni (1990), Dahlhaus (1996a, b, c, 1997), Bibi \& Francq (2003), Azrak \& M\'elard (2006) and Triantafyllopoulos \& Nason (2007). We refer to the introduction of Azrak \& M\'elard (2006) or of Van Bellegem \& Dahlhaus (2006) for further references. In several of these papers, the coefficients of the ARMA models are not constant but are deterministic functions of time. Also the innovation variance can be a deterministic function of time instead of being constant, such as in Van Bellegem \& von Sachs (2004). We can speak of marginal heteroscedasticity by opposition to conditional heteroscedasticity which is encountered in ARCH and GARCH models. All these functions of time are supposed to depend on a small number of parameters. 
Other somewhat related recent approaches include generalized autoregressive score (GAS) models of Creal et al. (2013) and testing parameter constancy against deterministically time-varying parameters, e.g.\ Ter\"asvirta et al (2010, Section 6%.3) and references therein, generalized to VAR models in Ter\"asvirta and Yang (2014).  %ALM_V16 %ALM_V17
.3) and references therein, generalized to VAR models in Ter\"asvirta \& Yang (2014).  %ALM_V16 %ALM_V17

The present paper inscribes itself in this line of research %GM 24/02/15
but for multivariate time series. 
%We however clarify from the beginning that we will not assume any periodicity on the time-dependent coefficients, as supposed in, e.g., Tiao and Grupe (1980) or Basawa and Lund (2001), as such an assumption is incompatible with slow changes in the dynamics of a process. 
The generality of the models we consider evidently entails numerous challenges, since the convenient asymptotic theory of stationary ergodic processes does no longer apply. Also, the asymptotic theory of time series models makes a large use of Fourier transforms and, consequently, of what is called spectral analysis. 
%Further difficulties inherent to general models with time-dependent coefficients are that the observations are (i) not independent, (ii) not identically distributed, and (iii) not normally distributed. 
As a consequence, deriving conditions for consistency and asymptotic normality of estimators of the coefficients, as well as obtaining the asymptotic covariance matrix, becomes highly complicated. 
%This explains the absence, in the literature, of a general theory on the asymptotic properties of estimators for models with time-dependent coefficients (Tj\o stheim~1984 comes closest to the achievement of such a theory with his work on linear and non-linear models). 

%Yet another difficulty we have to face here is the fact that w
We consider multidimensional time series models, with particular emphasis on %\emph{vector auto-regressive-moving average (VARMA)} models, the multivariate 
\emph{vector ARMA} (VARMA) models, the multivariate extension of the ARMA models. %celebrated ARMA models popularized by Box and Jenkins (1970, 1976). 
The main difference between ARMA and VARMA models lies in the fact that the coefficients change from scalars to squared matrices. The main developments in the area of statistical inference of standard stationary VARMA models are due to Kohn (1978), Hannan \& Deistler (1988), Francq \& Ra\"issi (2007) and, quite recently, Boubacar Mainassara \& Francq (2011) who study the consistency and asymptotic normality of quasi-maximum likelihood estimators for weak VARMA models. However, the field of \emph{time-dependent VARMA (tdVARMA)} models  with marginally heteroscedastic %innovation variance (meaning that both the  coefficients and that variance are %ATM_V13
innovation covariance matrix remains largely unexplored. %GM 24/02/15
An exception is Dahlhaus (2000) using an entirely different approach and assuming that the coefficients depend on time $t$ but also on the length of the series $n$ trough their ratio $t/n$. 
Here we assume dependency on $t$ only. 
Even if our theory is illustrated on pure VAR examples, it should be emphasized that it is valid for VMA and VARMA models, like in Dahlhaus (2000). Note that L\"utkepohl (2005, Chap. 14) treats tdVAR models by Gaussian maximum likelihood but does not discuss asymptotic properties in the general case.

Thus we want to fill in this gap in the literature by extending to the multivariate setting the methodology of Azrak \& M\'elard (2006) who, to the best of the authors' knowledge, were the first to obtain asymptotic properties of estimators for the general class of univariate time-dependent ARMA models by having recourse to \emph{quasi-maximum likelihood estimation (QMLE)}.

%The Azrak-M\'elard quasi-maximum likelihood approach acts as if the process under investigation were Gaussian.  Of course, one may argue that for clearly non-Gaussian data, this  methodology is wrong; but this criticism would then hold as well for the classical maximum likelihood estimation, where  the error then occurs at the very beginning with assuming the model to be Gaussian. Thus, in some sense, quasi-maximum likelihood is ``safer'' than the classical method. 
Like other QMLE approaches, the estimation method in Azrak \& M\'elard (2006) does not use the true, unknown, density of the observations but rather acts as if that density were Gaussian, thus using the Gaussian log-likelihood, which is an extension of the generalized least-squares method since it takes care of possible heteroscedasticity. There is no assumption of stationarity but, although it is not illustrated in our examples, there is an adjustment in the asymptotic theory for allowing non-normal observations. 
One major advantage of QMLE  is that the Gaussian likelihood function can be computed %exactly, with an efficient algorithm,  Alj \emph{\emph{et al.}} (2015c), %ALM_V13
exactly, with an efficient algorithm,  Alj  {et al.} (2015c), %ALM_V13
and this is very important for short time series. 
The main task in the Azrak-M\'elard approach, hence also in our extension, consists in checking conditions from two crucial theorems in Klimko \& Nelson (1978), which respectively ensure existence of an almost surely (a.s.) consistent estimator and prove asymptotic normality of that estimator, whilst providing the asymptotic covariance matrix. This is precisely what we are aiming at but, as we shall see in the rest of this paper, it is all but an easy task. 

Let us briefly comment on two other univariate approaches, and explain why we have not opted for extending those. 
\begin{itemize}
\item[-] The Dahlhaus approach (Dahlhaus 1996a, b, c, 1997) has allowed to obtain asymptotic results for a class of locally stationary processes including heteroscedastic ARMA processes with time-dependent coefficients. Dahlhaus uses either a spectral-based or a maximum likelihood estimation method. His asymptotics are based on rescaling time, i.e. $t/n$; therefore this is not related to our approach where there is no such requirement. Moreover local stationarity implies that the coefficients are continuous functions of time (and even two time differentiable functions), which is not necessarily the case here. Azrak \& M\'elard (2011) show a univariate example where Azrak \& M\'elard (2006)'s theory holds but the assumptions of local stationarity are not valid. Note also that Dahlhaus (2000)'s theory for multivariate processes assumes a Gaussian process whereas we assume only existence of 8th-order moments. %: instead of increasing the length $n$ of the series, he fixes the observation period to, say, $[0,1]$, with the interval between the $n$ observations decreasing and tending to zero. \red GM: and the VARMA case? \black
\item[-] The Bibi-Francq approach (Bibi \& Francq~2003) applies the quasi-least squares estimation method and gives asymptotic results for cyclical ARMA models with non constant periods. 
Although only 4th-order moments are needed, the theory was not developed for multivariate processes until now. 
\end{itemize}

%We however clarify from the beginning that 
%We do not assume any periodicity on the time-dependent coefficients, as supposed in, %ATM_V13
Another related approach, %ATM_V13
see for instance Tiao \& Grupe (1980) or Basawa \& Lund (2001), %as such an assumption is incompatible with slow changes in the dynamics of a process. 
consists in ARMA models with coefficients that vary as periodic functions of time, see  also Hindrayanto  {et al.} (2010). If the period $s$ is an integer, $s$ consecutive variables can be stacked as a vector which satisfies a stationary VARMA model. Here we do not assume periodic coefficients  %ATM_V13
%although our two examples will have coefficients that are periodic functions of time, to simplify the derivations, but with large or irrational periods. Therefore  that the stacking  %ATM_V13 
although, to simplify the derivations, our two examples will have periodic coefficients or innovation covariance matrix, but with large or irrational periods. Therefore stacking the variables will be practically inoperative and standard asymptotic theory for stationary VARMA models will not apply. %ATM_V13 

%{\red I WOULD REMOVE THIS
%The reasons for extending rather the Azrak-M\'elard approach instead of the other two are the following. The Dahlhaus approach requires local stationarity, a too strict restriction which would not allow to obtain as general results as those we shall obtain in this paper. The same comment applies to the Bibi-Francq approach, but their conditions moreover strongly differ from the Azrak-M\'elard ones. %(for example, they only need fourth-order moments but have several other requirements) and happen to be difficult to check in this setting.
%  %ATM_V13 
%  }
Most of the technical lemmas used in the present paper and their proofs are given in a technical appendix denoted here 'TA', see Alj  {et al.} (2015a). 
The second part of this paper (Alj  {et al.}, 2015b) will deal with the more general case where the coefficients of the model depend on time $t$ but also possibly on the number of observations $n$ of the series. However, even in that case, the theory differs from Dahlhaus' approach in the sense that the coefficients do not need to be continuous functions of time. In the univariate case, an example is provided by Azrak \& M\'elard (2011). 
The technical appendix, Alj  {et al.} (2015a), will be shared by Alj  {et al.} (2015b). The reason to separate the material in two parts is that while theorems related to martingale sequences are enough in the context of the present paper, the second part requires martingale arrays, although the technicalities are more or less the same.  

Another aspect of the present paper is that it provides an alternative theory for the asymptotics of standard VARMA models that does not rely on stationarity or ergodicity arguments, although in  the standard case, the assumptions will imply the usual conditions on the roots of the autoregressive and moving average polynomials in the lag operator. Our alternative theory also avoids spectral analysis. %ATM_V13 
	
The paper is organized as follows. In Section~\ref{S2}, we first develop asymptotics for quasi-maximum likelihood estimators in a general multivariate time series model which is not necessarily stationary. Then, in Section~\ref{S.3}, we focus our attention on  tdVARMA models: after setting the notations, we analyze pure VAR and pure VMA representations, with an illustration, before finally stating the main theorem for the tdVARMA case. We %illustrate our theoretical findings by means of a few examples in Section~\ref{S.4}, and examine the finite-sample behavior of our estimators via a Monte Carlo simulation study in Section~\ref{S.5}. Finally, an appendix collects the technical details. %ATM_V13
illustrate our theoretical findings by means of two examples in Section~\ref{S.5}, and examine the finite-sample behavior of our estimators via a Monte Carlo simulation study in Section~\ref{S.5.4}. Finally, Appendix A collects the main proofs and Appendix B contains a verification of the main assumptions for the two examples studied with a few nice mathematical derivations. %ATM_V13

%{\red Remarque sur les r\'ef\'erences: il manque Golub and Van Loan (1996) dans la liste. 
%}
%%%%%%%%%%%%%%%%%%%%%%%%%%%%%%%%%%%%%%%%%%%%%%%%%%%%%%%%%%%%%%%%%%%%%%%
%%%%%%%%%%%%%%%%%%%%%%%%%%%%%%%%%%%%%%%%%%%%%%%%%%%%%%%%%%%%%%%%%%%%%%%
%%%%%%%%%%%%%%%%%%%%%%%%%%%%%%%%%%%%%%%%%%%%%%%%%%%%%%%%%%%%%%%%%%%%%%%
\section{QMLE for a general multivariate time series model}\label{S2}

%%%%%%%%%%%%%%%%%%%%%%%%%%%%%%%%%%%%%%%%%%%%%%%%%%%%%%%%%%%%%%%%%%%%%%%
%%%%%%%%%%%%%%%%%%%%%%%%%%%%%%%%%%%%%%%%%%%%%%%%%%%%%%%%%%%%%%%%%%%%%%%
\subsection{Some {preliminaries}}\label{S2.1}
%%%%%%%%*********************************************************************************%%%%%%%%%%%%%%%%%%%%%%%%%%%%
Let $\left\{x_{t} : t\in \Nb \right\}$ be a stochastic process defined on a probability space $(\Omega, F, P_\theta )$, taking values in $\Rb^r$, and whose distribution depends on a vector $\theta = (\theta_1 , ... , \theta_m )^{T} $ of unknown parameters to be estimated, with $\theta$ lying in some open set $\Theta$ of a Euclidean space $\Rb^m$. Let $E_{\theta}(.)$ and $E_{\theta}(./.)$ denote expectation and conditional expectation under $P_\theta$, respectively. The true value of $\theta$ is denoted by $\theta^{0} = (\theta^{0}_{1}, ... , \theta^{0}_{m} )^{T}$, assumed to be an interior point of $\Theta$. 
%The latter statement is taken to mean that all probabilities, all  a.s. statements and all non subscripted expectations and conditional expectations are taken relative to the measure determined by $\theta^{0}$.
Let $\left\{F_{t} : t \in \Nb\right\}$ be an increasing sequence of sub-sigma algebras of $F$ with $F_{t}$ generated by $\{x_u : u = 1,2,...,t\}$ with $F_{0}=\{\emptyset,\Omega\}$ such that, for each $t$, $x_t$ is measurable with respect to $F_{t}$. Given a set of observations $\left\{x_{t} : t= 1,2,...,n \right\}$, we want to estimate $\theta$ by trying to minimize the general real-valued objective function $Q_{n}(\theta ) = Q_{n}(\theta;x_{1}, ..., x_{n})$ which depends on $\theta$ and the observations $\{x_{t} : t=1, 2, ..., n\}$. Therefore we solve the system of equations 
\begin{eqnarray*}
\frac{\partial Q_{n}(\theta )}{\partial \theta _{i}} =0 \quad \text{for} \quad i=1, ...,m.
\end{eqnarray*}  

We suppose that the objective function $Q_{n}(\theta )$ is twice continuously differentiable in $\Theta$. Let $\widehat{\theta}_{n} = (\widehat{\theta}_1 , ... , \widehat{\theta}_m )^{T} $ be a sequence of estimators indexed by $n$. Klimko \& Nelson (1978) showed conditions for strong consistency and asymptotic normality of $\widehat{\theta}_{n}$, see also Hall \& Heyde (1980, pp. 174-176) and Taniguchi \& Kakizawa (2000, pp. 97-98). 
%%%%%%%%%%%%%%%%%%%%%%%%%%%%%%%%%%%%%%%%%%%%%%%%%%%%%%%%%%%%%%%%%%%%%%%%%%%%%%%%%%%%%%%%%%%%%%%%%%%%%%%%%%%%%%%%%%%%%%%%%%%%%%%%%%%%%%%%%%%%%%%%%%%%%%%%%%%%%%%
\subsection{General theory of quasi-maximum likelihood estimation}\label{S2.2}
%%%%%%%%%%%%%%%%%%%%%%%%%%%%%%%%%%%%%%%%%%%%%%%%%%%%%%%%%%%%%%%%%%%%%%%%%%%%%%%%%%%%%%%%%%%%%%%%%%%%%%%%%%%%%%%%%%%%%%%%%%%%%%%%%%%%%%%%%%%%%%%%%%%%%%%%%%%%%%%
%Let $\left\{x_{t} : t\in \Nb \right\}$ be a $r$-vector stochastic process defined in Section \ref{S2.1} and consider the conditional expectations and :% respectively by :
Denote
\begin {equation}
 e_{t}(\theta )= x_{t}-\hat{x}_{t/t-1}(\theta) \quad \text{with} \quad  \hat{x}_{t/t-1}(\theta) = E_{\theta}(x_{t}/F_{t-1}),\label{21.1}
\end {equation}
for which obviously $E_{\theta}(e_{t}(\theta )) = 0$. We denote by 
$$
%\Sigma_{t/t-1}(\theta )= E_{\theta}\left[e_{t}(\theta ) e^{T}_{t}(\theta )/F_
\Sigma_{t}(\theta )= E_{\theta}\left[e_{t}(\theta ) e^{T}_{t}(\theta )/F_{t-1}\right]\label{21.1.bis}
$$
the conditional covariance matrix given $F_{t-1}$. The quasi-likelihood function $L_{n}(\theta ;x_{1},...,x_{n})$ computed as if the process were Gaussian is given by
\begin{eqnarray}
L_{n}(\theta;x_{1},...,x_{n})&=& (2\pi)^{-nr/2}{\prod_{t=1}^{n}} \text{det}%\left(\Sigma_{t/t-1}(\theta)\right)^{-1/2}\exp\left\{-\frac{1}{2}{\sum_{t=1}^{n}}e_{t}^{T}(\theta )\Sigma^{-1}_{t/t-1}(\theta )e_{t}(\theta )\right\}.
\left(\Sigma_{t}(\theta)\right)^{-1/2}\exp\left\{-\frac{1}{2}{\sum_{t=1}^{n}}e_{t}^{T}(\theta )\Sigma^{-1}_{t}(\theta )e_{t}(\theta )\right\}.\nonumber
\end{eqnarray}	
We take the objective function 
\begin{eqnarray}	
Q_{n}(\theta )&=& -\log \left(L_{n}(\theta;x_{1},...,x_{n})\right) = \frac{1}{2}{\sum_{t=1}^{n}}\alpha _{t}(\theta )+ \frac{rn}{2}\log (2\pi), \label{22.1} 
\end{eqnarray}	
with
$$
%\alpha_{t}(\theta )=\log\left(\text{det}\left(\Sigma_{t/t-1}(\theta)\right) \right) + e_{t}^{T}(\theta )\Sigma^{-1}_{t/t-1}(\theta )e_{t}(\theta ). 
\alpha_{t}(\theta )=\log\left(\text{det}\left(\Sigma_{t}(\theta)\right) \right) + e_{t}^{T}(\theta )\Sigma^{-1}_{t}(\theta )e_{t}(\theta ). 
$$
Then the QMLE of $\theta $ is defined as any measurable solution $\widehat{\theta}_{n}$ of 
\begin{eqnarray}
%\widehat{\theta}_{n} = \underset{\theta \in \Theta}{\text{arg min}}\quad Q_{n}(\theta ) \label{2.3}
{\text{arg min}}_{\theta \in \Theta}\quad Q_{n}(\theta ). \label{22.3}
\end{eqnarray}

In order to check the assumptions of the Klimko \& Nelson (1978) theorems, we proceed like Azrak \& M\'elard (2006) and we make some additional assumptions as follows. Let the $r$-vector stochastic process $\left\{x_{t} : t\in \Nb\right\}$ be such that $E_{\theta} \left(\left\|x_{t}\right\|^{2}%\right)< \infty$ for all $\theta$ and that $e_{t}(\theta)$ and $\Sigma_{t/t-1
\right)< \infty$ for all $\theta$ and   $e_{t}(\theta)$ and $\Sigma_{t}(\theta)$ are almost surely twice continuously differentiable in $\Theta$. Henceforth, for simplicity, we denote $\left[E_{\theta}\{.(\theta)\}\right]_{\theta = \theta^{0}}$ by $E_{\theta^{0}}\{.(\theta)\}$.
%\newpage
We suppose that there exist two positive constants $C_1$ and $C_2$ such that for all $t \geq 1$:
\begin{description}
	\item[$\Hb_{2.1}$] $\quad \quad E_{\theta^{0}} \left\{\left|\frac{\partial \alpha_{t}(\theta )}{\partial \theta_{i}}\right|^{4}\right\} \leq  C_{1}$ for $i=1, ...,m;$
	\item[$\Hb_{2.2}$] $\quad \quad E_{\theta^{0}}\left\{\left|\frac{\partial^{2}\alpha_{t}(\theta)}{\partial\theta_{i}\partial\theta_{j}}- E_{\theta}\left(\frac{\partial^{2}\alpha_{t}(\theta)}{\partial\theta_{i}\partial\theta_{j}}/F_{t-1}\right)\right|^{2}\right\} \leq  C_{2}$ for $ i,j=1,...,m. $

\noindent Suppose further that 
	\item[$\Hb_{2.3}$] 
$$ \underset{n\rightarrow \infty }{\lim} \frac{1}{2n} {\sum_{t=1}^{n}} E_{\theta^{0}}\left\{\frac{\partial^{2}\alpha_{t}(\theta)}{\partial\theta_{i}\partial\theta_{j}}/F_{t%-1}\right\}=V_{ij}(\theta^{0}) \quad \text{a.s.}\quad \text{for} \quad  i,j =1, ..., m,$$ %ATM_V13
%\noindent where $V(\theta)=(V_{ij}(\theta^{0}))_{1\leq i,j\leq m}$ is a strictly positive definite matrix of constants;  %ATM_V13
-1}\right\}=V_{ij} \quad \text{a.s.}\quad \text{for} \quad  i,j =1, ..., m,$$ %ATM_V13
\noindent where $V=(V_{ij})_{1\leq i,j\leq m}$ is a strictly positive definite matrix of constants;  %ATM_V13

	\item[$\Hb_{2.4}$] $$ \underset{n\rightarrow \infty }{\lim} \underset{ \Delta\downarrow 0}{\sup}(n\Delta)^{-1}\left|\sum_{t=1}^{n}\left(\left\{\frac{\partial^{2} \alpha_{t}(\theta)}{\partial \theta_{i}\partial\theta_{j}}\right\}_{\theta = \theta^{*}}- \left\{\frac{\partial^{2} \alpha_{t}(\theta)}{\partial \theta_{i}\partial\theta_{j}}\right\}_{\theta = \theta^{0}} \right)\right| < \infty \quad \text{a.s.}\quad $$
\noindent for $i, j=1, ...,m$, where $\theta^{*}$ is a point of the straight line %joining $\theta^{0}$ to every $\theta$, such that $\|\theta - \theta^{0} \|< \Delta$, $0 <\Delta$. %ALM_V13
joining $\theta^{0}$ to every $\theta$, such that $\|\theta - \theta^{0} \|< \Delta$, $0 <\Delta$, where $\| . \|$ is the Euclidean norm. %ALM_V13
\end{description}

The main notations and assumptions being settled, we are ready to state the two main theorems of the present paper.

\begin{theor}\label{T1.2.3}
Suppose that Assumptions $\Hb_{2.1}-\Hb_{2.4}$ hold. Then there exists a sequence of estimators  $\widehat{\theta}_{n} = (\widehat{\theta}_1 , ... , \widehat{\theta}_m )^{T}$ such that $\widehat{\theta}_{n}\rightarrow \theta^{0}$ a.s. and, for any $\eps>0$, there exists an event $ E $ with $P_{\theta^{0}}(E)>1-\eps$ and an $n_{0}$ such that on $E$, for any $n >n_0$, $\{\partial Q_{n}(\theta)/\partial\theta_{i}\}_{\theta = \widehat{\theta}_{n}} = 0$, for $i = 1,2,...,m, $ and $Q_{n}(\theta )$ attains a relative minimum at~$\widehat{\theta}_{n}$.
\end{theor} 
The proof of Theorem \ref{T1.2.3} is given in  Appendix~\ref{ASS0}.
\begin{theor}\label{T1.2.5}
If the assumptions $\Hb_{2.1}-\Hb_{2.4} $ are satisfied, as well as
%\begin{description} %ALM_V13
% \item[$\Hb_{2.5}$]  $$\quad \quad \frac{1}{n}{\sum_{t=1}^{n}}E_{\theta^{0}}\left(\frac %ALM_V13
 %ATM_V13
$\Hb_{2.5}$  $$\quad \quad \frac{1}{n}{\sum_{t=1}^{n}}E_{\theta^{0}}\left(\frac %ALM_V13
{\partial \alpha_{t}(\theta)}{\partial \theta}\frac{\partial \alpha_{t}(\theta)}{\partial \theta^{T}}/F_{t-1}\right)-\frac{1}{n}{\sum_{t=1}^{n}}E_{\theta^{0}}\left(\frac{\partial \alpha_{t}(\theta)}{\partial \theta}\frac{\partial \alpha_{t}(\theta)}{\partial \theta^{T}}\right)\stackrel{\text{a.s.}}{\rightarrow} 0 \quad \textrm{as} \quad n\rightarrow \infty,$$ 
%\end{description} %ALM_V13
then
%  $$ n^{1/2}(\widehat{\theta}_{n} - \theta^{0}) \stackrel{L}{\rightarrow} \mathcal{N}(0,V(\theta^{0})^{-1}W(\theta^{0})V(\theta^{0})^{-1}),$$ %ATM_V13
  $$ n^{1/2}(\widehat{\theta}_{n} - \theta^{0}) \stackrel{L}{\rightarrow} \mathcal{N}(0,V^{-1} W V^{-1}),$$ %ATM_V13
%where $\stackrel{L}{\longrightarrow}$ indicates convergence in law and $W(\theta^{0})}$ is a strictly positive definite matrix defined by %ALM_V13
where $\stackrel{L}{\longrightarrow}$ indicates convergence in law and $W=(W_{ij})_{1\leq i,j\leq m}$ is a strictly positive definite matrix defined by %ALM_V13
$$
%W(\theta^{0}) = \underset{n\rightarrow \infty }{\lim}\frac{1}{4n}{\sum_{t=1}^{n}}E_{\theta^{0}}\left(\frac{\partial \alpha_{t}(\theta)}{\partial \theta}\frac{\partial \alpha_{t}(\theta)}{\partial \theta^{T}}\right). \label{Wmatrix}  %ATM_V13
W = \underset{n\rightarrow \infty }{\lim}\frac{1}{4n}{\sum_{t=1}^{n}}E_{\theta^{0}}\left(\frac{\partial \alpha_{t}(\theta)}{\partial \theta}\frac{\partial \alpha_{t}(\theta)}{\partial \theta^{T}}\right).   %ATM_V13
$$
\end{theor} 
For the proof of Theorem \ref{T1.2.5} we use the Central Limit Theorem for martingale %differences of Basawa and Prakasa Rao (1980).  It is also given in Appendix~\ref{ASS0}. %ALM_V13  
differences of Basawa \& Prakasa Rao (1980). %ALM_V13  

%%%%%%%%*********************************************************************************%%%%%%%%%%%%%%%%%%%%%%%%%%%%
%%%%%%%%%%%%%%%%%%%%%%%%%%%%%%%%%%%%%%%%%%%%%%%%%%%%%%%%%%%%%%%%%%%%%%%
%%%%%%%%%%%%%%%%%%%%%%%%%%%%%%%%%%%%%%%%%%%%%%%%%%%%%%%%%%%%%%%%%%%%%%%
%%%%%%%%%%%%%%%%%%%%%%%%%%%%%%%%%%%%%%%%%%%%%%%%%%%%%%%%%%%%%%%%%%%%%%%
\section{VARMA models with time-dependent coefficients}\label{S.3}
%%%%%%%%*********************************************************************************%%%%%%%%%%%%%%%%%%%%%%%%%%%%
%%%%%%%%%%%%%%%%%%%%%%%%%%%%%%%%%%%%%%%%%%%%%%%%%%%%%%%%%%%%%%%%%%%%%%%
%%%%%%%%%%%%%%%%%%%%%%%%%%%%%%%%%%%%%%%%%%%%%%%%%%%%%%%%%%%%%%%%%%%%%%%
\subsection{tdVARMA models: definition and notations} \label{S3.1}
The process $\{x_{t} : t\in \Nb \}$ is called a zero mean $r$-vector mixed %autoregressive moving average process of order ($p$, $q$) with time-dependent %ALM_V16
autoregressive-moving average process of order ($p$, $q$) with time-dependent %ALM_V16
 coefficients, and is denoted by tdVARMA ($p$, $q$), if and only if it satisfies the  equation
\begin{eqnarray}
x_{t} = \sum_{i=1}^{p}A_{ti}x_{t-i} + g_{t}\epsilon_{t}+ \sum_{j=1}^{q}B_{tj}g_{t-j}\epsilon_{t-j},   \label{3.1} 
\end{eqnarray}
where $p$ and $q$ are integer constants, $\{\epsilon_{t} : t\in \Nb \}$ is an independent white noise process, consisting of independent random variables, not %necessarily identically distributed, with zero mean and covariance matrix $\Sigma$ %which is  invertible, and where the coefficients $A_{t1},$ $...,A_{tp}$ and $B_{t1} %ALM_V13
necessarily identically distributed, with zero mean, covariance matrix $\Sigma$ which is  invertible, and finite fourth-order moments, and where the coefficients $A_{t1},$ $...,A_{tp}$ and $B_{t1} %ALM_V13
,..., B_{tq}$, as well as the $r\times r$ matrix $g_{t}$, are deterministic functions of time $t$. The initial values $x_{t}, t < 1$, and $\epsilon_{t}$, $t < 1$, are supposed to be equal to the zero vector. In the sequel, we will also use $A_{t0} = B_{t0} = I_{r}$ with $I_{r}$ the $r$-dimensional identity matrix and set to zero the %coefficients $A_{tk}$ with $k > p$ and $B_{tk}$ with $k > q$, for all $t$. For $k, l  %ALM_V13
%\in \Nb$ and $k \neq l $, we denote by $ \kappa(\Sigma)$ the matrix %ALM_V13
coefficients $A_{tk}$ with $k > p$ and $B_{tk}$ with $k > q$, for all $t$.  %ALM_V13
%Let us define $F_{t}$ as being generated by $\{\epsilon_u : u = 1,2,...,t\}$. 
 %ALM_V13
Writing $\otimes$ the Kronecker product, we let  %ALM_V13
\begin{eqnarray} %ALM_V13
 \kappa_{t} = E\left(\vecc(\epsilon_{t}\epsilon^{T}_{t}) \vecc(\epsilon_{t}\epsilon^{T}_{t})^{T}\right) = E\left((\epsilon_{t} \epsilon^{T}_{t})\otimes(\epsilon_{t} \epsilon^{T}_{t})\right),\nonumber \label{kappa2}  %ALM_V13
\end{eqnarray} %ALM_V13
which depends on $t$, in general.  %ALM_V13
For $k, l  %ALM_V13
\in \Nb$ and $k \neq l $, we consider the matrix %ALM_V13
$$
  E\left(\vecc(\epsilon_{t-k}\epsilon^{T}_{t-l})\vecc(\epsilon_{t-l}\epsilon^{T}_{t-k})^{T} \right) = K_{r,r}(\Sigma \otimes \Sigma),
$$
%which does not depend on $t, k$ or $l$ and where the $r^{2}\times r^{2}$  matrix $K_{r,r}$ is the \emph{commutation matrix}; see Lemma 4.2 of the technical appendix for a proof of equality~\eqref{comm}. The latter matrix can be described in the following way: in the $(i, j)$-th block of $K_{r,r}$ (each block is of size $r\times r$) the $(j, i)$-th element equals one, while all other elements in that block are zeros, see Kollo and von Rosen (2005, p. 79). Moreover, we have for $k = l =0$  %ALM_V13
which does not depend on $t, k$ or $l$ and where the $r^{2}\times r^{2}$  matrix $K_{r,r}$ is the \emph{commutation matrix}. See Kollo \& von Rosen (2005, p. 79) or TA Lemma 4.2. %ALM_V13
%\begin{eqnarray} %ALM_V13
% \kappa_{t}(\Sigma) = E\left(\vecc(\epsilon_{t}\epsilon^{T}_{t}) %\vecc(\epsilon_{t}\epsilon^{T}_{t})^{T}\right) = E\left((\epsilon_{t}.\epsilon^{T}_{t})\otimes(\epsilon_{t}.\epsilon^{T}_{t})\right),\nonumber \label{kappa2}  %ALM_V13
%\end{eqnarray} %ALM_V13
%which depends on $t$, in general. If the process $\{\epsilon_{t}, t\in \Nb\} $ is Gaussian, then $\kappa_{t}$ does not depend on $t$ and  %ALM_V13
%\begin{eqnarray*} %ALM_V13
% \kappa_{t}&=& \vecc(\Sigma) \vecc(\Sigma)^{T} + (I_{r^{2}} +  K_{r,r})(\Sigma \otimes \Sigma), %ALM_V13
%\end{eqnarray*}   %ALM_V13 
% see Kollo and von Rosen (2005, p. 207) for more details. %ALM_V13

Let us now consider the parametric model corresponding to (\ref{3.1}), namely
\begin{eqnarray}
x_{t} = \sum_{i=1}^{p}A_{ti}(\theta)x_{t-i} + e_{t}(\theta)+ \sum_{j=1}^{q}B_{tj}(\theta)e_{t-j}(\theta),  \label{3.2}
\end{eqnarray}
where the $e_{t}(\theta)$ can be considered as the residuals of the model and are defined as in (\ref{21.1}) and where $A_{ti}(\theta)$, for $i=1,...,p,$ and $B_{tj}(\theta)$, for $j=1,...,q,$ are the parametric coefficients. Furthermore the covariance matrix of $e_{t}(\theta)$ is parametrized as 
%$$E_{\theta}\left(e_{t}(\theta)e^{T}_{t}(\theta)\right) = g_{t}(\theta)\Sigma g^{T}_{t}(\theta) = \Sigma_{t}(\theta ).$$   %ALM_V17
$$\Sigma_{t}(\theta ) =^{\text{def}} E_{\theta}\left(e_{t}(\theta)e^{T}_{t}(\theta)\right) = g_{t}(\theta)\Sigma g^{T}_{t}(\theta).$$   %ALM_V17
%ALM_V13
For $\theta = \theta^{0}$, we have $A_{ti}(\theta^{0})= A_{ti}$, 
%for $i=1,...,p$ and 
$B_{tj}(\theta^{0})= B_{tj}$,  %for $j=1,...,q$, 
$g_{t}(\theta^{0})= g_{t}$, $e_{t}(\theta^{0})= g_{t}\epsilon_{t}$ and $\Sigma_{t}=^{\text{def}} \Sigma_{t}(\theta^{0})$ $= E\left(e_{t}(\theta^{0})e^{T}_{t}(\theta^{0})\right)= g_{t}\Sigma g^{T}_{t}$. We assume that the $m$-dimensional vector $\theta$ contains all the parameters of interest to be estimated, those in the coefficients $A_{t1}(\theta), ..., $  $A_{tp}(\theta)$, $B_{t1}(\theta), ..., B_{tq}(\theta)$ and $g_{t}(\theta)$ but not the nuisance parameters in the scale factor matrix $\Sigma$ which are estimated separately. In usual VARMA($p$, $q$) models, the coefficients $A_{1}(\theta), ..., A_{p}(\theta)$, $B_{1}(\theta), ..., B_{q}(\theta)$ and $g_{t}(\theta)$ do not depend on $t$, and the parameters are the coefficients themselves. 
%ALM_V13
Note that for a given $\theta$ we have% by using the same notations in Section \ref{S1.2.2}
\begin{eqnarray*}
\hat{x}_{t/t-1}(\theta) &=& E_{\theta}(x_{t}/F_{t-1}) = \sum_{i=1}^{p}A_{ti}(\theta)x_{t-i} + \sum_{j=1}^{q}B_{tj}(\theta) e_{t-j}(\theta).
\end{eqnarray*}
%\begin{eqnarray*}
%e_{t}(\theta ) &=& x_{t}-\hat{x}_{t/t-1}(\theta)\\ \nonumber
%               &=& g_{t}(\theta)\epsilon_{t}
%\end{eqnarray*}.
%The conditional covariance matrix given $F_{t-1}$ is %ALM_V13
%\begin{eqnarray*} %ALM_V13
%%\Sigma_{t/t-1}(\theta ) &=& E_{\theta}\left[e_{t}(\theta ) e^{T}_{t}(\theta )
%\Sigma_{t}(\theta ) &=& E_{\theta}\left[e_{t}(\theta ) e^{T}_{t}(\theta )/F_{t-1}\right] = g_{t}(\theta)\Sigma g^{T}_{t}(\theta).\nonumber %ALM_V13
%&=& \Sigma_{t}(\theta ),
%\end{eqnarray*} %ALM_V13
%As we can see,  %ALM_V13
%%the conditional covariance matrix $\Sigma_{t/t-1}(\theta)$ is equal to $\Sigma_{t}(\theta)$ in this case and $\Sigma_{t/t-1}(\theta^{0}) = 
%$\Sigma_{t}=^{\text{def}} \Sigma_{t}(\theta^{0}) = g_{t}\Sigma g^{T}_{t}.$  %ALM_V13

%Let $e_{t}(\theta)$ be the residual like the section before. Of course  $e_{t}(\theta^{0})= g_{t}(\theta^{0})e_{t}$ are the innovation of the process $\{x_{t}\}$.

According to the assumptions made about initial values, it is possible to write out properly the pure autoregressive and the pure moving average representation of the model (\ref{3.2}), as we shall see in the next section.

%%%%%%%%%%%%%%%%%%%%%%%%%%%%%%%%%%%%%%%%%%%%%%%%%%%%%%%%%%%%%%%%%%%%%%%
%%%%%%%%%%%%%%%%%%%%%%%%%%%%%%%%%%%%%%%%%%%%%%%%%%%%%%%%%%%%%%%%%%%%%%%
\subsection{The pure autoregressive and the pure moving average representations} \label{S1.3.1.2}
%By using the assumption about initial values and using (\ref{3.2}) recurrently, M\'elard (1985) and Azrak \& M\'elard (2012) have established expressions for the pure %ALM_V13
By using the assumption about initial values and using (\ref{3.2}) recurrently, M\'elard (1985) and Azrak \& M\'elard (2015) have established expressions for the pure %ALM_V13
 autoregressive representation and the pure moving average representation of tdVARMA processes. In the univariate case these representations can be found in Azrak \& M\'elard (2006).  In our setting, for any $\theta$ the pure autoregressive representation corresponds to 
\begin{eqnarray}
x_{t} =  \sum_{k=1}^{t-1}\pi_{tk}(\theta)x_{t-k} + e_{t}(\theta),\label{3.3}
\end{eqnarray}
where the coefficients $\pi_{tk}(\theta)$ can be obtained from the autoregressive and moving average coefficients by using the following recurrences (see M\'elard, 1985, pp. 43-45):
\begin{eqnarray*}
\pi^{(0)}_{t0}(\theta) &=& I_{r}, \quad \pi^{(0)}_{tj}(\theta) = A_{tj}(\theta), \quad \widetilde{\pi}^{(0)}_{tj}(\theta) = B_{tj}(\theta), \quad \text{for} \quad j = 1,...,t-1,\\
\pi^{(k)}_{tj}(\theta) &=& \pi^{(k-1)}_{tj}(\theta) - \widetilde{\pi}^{(k-1)}_{tk}(\theta)A_{t-k, j-k}(\theta), \quad \text{for} \quad j = k,...,t-1,\\
\widetilde{\pi}^{(k)}_{tj}(\theta) &=& \widetilde{\pi}^{(k-1)}_{tj}(\theta) - \widetilde{\pi}^{(k-1)}_{tk}(\theta)B_{t-k, j-k}(\theta), \quad \text{for} \quad j = k+1,...,t-1,
\end{eqnarray*}
and $\pi_{tk}(\theta) = \pi^{(k)}_{tk}(\theta)$ for $k = 1,...,t-1$. By (\ref{3.3}) we of course have $
e_{t}(\theta) = x_{t} - \sum_{k=1}^{t-1}\pi_{tk}(\theta)x_{t-k}$, and consequently its first three derivatives with respect to $\theta$ are given by  
\begin{eqnarray}
\frac{\partial e_{t}(\theta)}{\partial\theta_{i}}&=&-\sum_{k=1}^{t-1}\frac{\partial\pi_{tk}(\theta)}{\partial\theta_{i}}x_{t-k}, \nonumber\\
\frac{\partial^{2} e_{t}(\theta)}{\partial\theta_{i}\partial\theta_{j}}&=& -\sum_{k=1}^{t-1}\frac{\partial^{2}\pi_{tk}(\theta)}{\partial\theta_{i}\partial\theta_{j}}x_{t-k},\label{32.11} \\
\frac{\partial^{3} e_{t}(\theta)}{\partial\theta_{i}\partial\theta_{j}\partial\theta_{l}} &=& -\sum_{k=1}^{t-1}\frac{\partial^{3}\pi_{tk}(\theta)}{\partial\theta_{i}\partial\theta_{j}\partial\theta_{l}}x_{t-k},\nonumber
\end{eqnarray}
for $i, j, l =1,..., m$. 

On the other hand, for the pure moving average representation we have 
\begin{eqnarray}
x_{t} = e_{t}(\theta) + \sum_{k=1}^{t-1}\psi_{tk}(\theta)e_{t-k}(\theta), \label{31.12}
\end{eqnarray}
where the coefficients $\psi_{tk}(\theta) = \psi^{(k)}_{tk}(\theta),$ $ k =0, 1,..., t-1,$ can be obtained from the autoregressive and moving average coefficients by using the following recurrences (see M\'elard, 1985, pp. 36-38): 
\begin{eqnarray*}
\psi^{(0)}_{t0}(\theta) &=& I_{r}, \quad \psi^{(0)}_{tj}(\theta) = B_{tj}(\theta), \quad \widetilde{\psi}^{(0)}_{tj}(\theta) = A_{tj}(\theta), \quad \text{for} \quad j = 1,...,t-1,\\
\psi^{(k)}_{tj}(\theta) &=& \psi^{(k-1)}_{tj}(\theta) + \widetilde{\psi}^{(k-1)}_{tk}(\theta)B_{t-k, j-k}(\theta), \quad \text{for} \quad j = k,...,t-1,\\
\widetilde{\psi}^{(k)}_{tj}(\theta) &=& \widetilde{\psi}^{(k-1)}_{tj}(\theta) + \widetilde{\psi}^{(k-1)}_{tk}(\theta)A_{t-k, j-k}(\theta), \quad \text{for} \quad j = k+1,...,t-1,
\end{eqnarray*}
for each $ k = 1,..., t-1$. Hence for $\theta = \theta^{0}$ we have 
$$
%x_{t} = g_{t}(\theta^{0})\epsilon_{t}+ \sum_{k=1}^{t-1}\psi_{tk}(\theta^{0})g_{t-k}(\theta^{0})\epsilon_{t-k}. \label{32.12} %ALM_V13
x_{t} = g_{t} \epsilon_{t}+ \sum_{k=1}^{t-1}\psi_{tk} g_{t-k} \epsilon_{t-k},  %ALM_V13
$$
where we denote $\psi_{tk}=\psi_{tk}(\theta^{0})$. %ALM_V13
Then, by using (\ref{31.12}), $e_{t}(\theta)$ and its first three derivatives in (\ref{32.11}) can be written as a pure moving average in terms of the innovations process: 
\begin{equation}
% e_{t}(\theta) = g_{t}(\theta^{0})\epsilon_{t} + \sum_{k=1}^{t-1}\psi_{t0k}(\theta,\theta^{0})g_{t-k}(\theta^{0})\epsilon_{t-k}, \label{32.13new} %ALM_V13
 e_{t}(\theta) = g_{t} \epsilon_{t} + \sum_{k=1}^{t-1}\psi_{t0k}(\theta,\theta^{0})g_{t-k} \epsilon_{t-k}, \label{32.13new} %ALM_V13
\end{equation}
\begin{equation}
%\frac{\partial e_{t}(\theta)}{\partial\theta_{i}}= \sum_{k=1}^{t-1}\psi_{tik}(\theta,\theta^{0})g_{t-k}(\theta^{0})\epsilon_{t-k}, \label{32.13} %ALM_V13
\frac{\partial e_{t}(\theta)}{\partial\theta_{i}}= \sum_{k=1}^{t-1}\psi_{tik}(\theta,\theta^{0})g_{t-k} \epsilon_{t-k}, \label{32.13} %ALM_V13
\end{equation}
\begin{equation}
%\frac{\partial^{2} e_{t}(\theta)}{\partial\theta_{i}\partial\theta_{j}}=\sum_{k=1}^{t-1}\psi_{tijk}(\theta,\theta^{0}) g_{t-k}(\theta^{0})\epsilon_{t-k}, \label{32.14}  
\frac{\partial^{2} e_{t}(\theta)}{\partial\theta_{i}\partial\theta_{j}}=\sum_{k=1}^{t-1}\psi_{tijk}(\theta,\theta^{0}) g_{t-k} \epsilon_{t-k}, \label{32.14}  %ALM_V13
\end{equation}
\begin{equation}
%\frac{\partial^{3} e_{t}(\theta)}{\partial\theta_{i}\partial\theta_{j}\partial\theta_{l}}=\sum_{k=1}^{t-1}\psi_{tijlk}(\theta,\theta^{0}) g_{t-k}(\theta^{0})\epsilon_{t-k},\label{32.15}  %ALM_V13
\frac{\partial^{3} e_{t}(\theta)}{\partial\theta_{i}\partial\theta_{j}\partial\theta_{l}}=\sum_{k=1}^{t-1}\psi_{tijlk}(\theta,\theta^{0}) g_{t-k} \epsilon_{t-k},\label{32.15}  %ALM_V13
\end{equation}
for $i, j, l = 1,..., m$, where the coefficients $\psi_{t0k}(\theta,\theta^{0})$, $\psi_{tik}(\theta,\theta^{0}),\psi_{tijk}(\theta,\theta^{0})$ and $\psi_{tijlk}(\theta,\theta^{0})$ are obtained from the autoregressive and moving average coefficients by the following relations: 
\begin{eqnarray}
\psi_{t0k}(\theta,\theta^{0})&=&\psi_{tk}(\theta^{0}) - \sum_{u=1}^{k} \pi_{tu}(\theta)% \psi_{t-u,k-u}(\theta^{0}),\label{p.0} \\ %ALM_V13
 \psi_{t-u,k-u},\label{p.0} \nonumber\\ %ALM_V13
\psi_{tik}(\theta,\theta^{0})&=&-\sum_{u=1}^{k}\frac{\partial\pi_{tu}(\theta)}{%\partial\theta_{i}}\psi_{t-u,k-u}(\theta^{0}),\label{p.1} \\ %ALM_V13
\partial\theta_{i}}\psi_{t-u,k-u},\label{p.1} \\ %ALM_V13
\psi_{tijk}(\theta,\theta^{0})&=&-\sum_{u=1}^{k}\frac{\partial^{2}\pi_{tu}(\theta)}{%\partial\theta_{i}\partial\theta_{j}}\psi_{t-u,k-u}(\theta^{0}),  \label{p.2} \\ 
\partial\theta_{i}\partial\theta_{j}}\psi_{t-u,k-u},  \label{p.2} \\ %ALM_V13
\psi_{tijlk}(\theta,\theta^{0})&=&-\sum_{u=1}^{k}\frac{\partial^{3}\pi_{tu}(\theta)}{%\partial\theta_{i}\partial\theta_{j}\partial\theta_{l}}\psi_{t-u,k-u}(\theta^{0}).
\partial\theta_{i}\partial\theta_{j}\partial\theta_{l}}\psi_{t-u,k-u}.\label{p.3} %ALM_V13
\end{eqnarray}
We denote 
$$ \psi_{t0k} = \psi_{t0k}(\theta^{0},\theta^{0}) ,\quad  \psi_{tik} = \psi_{tik}(\theta^{0},\theta^{0}) ,\quad  \psi_{tijk} = \psi_{tijk}(\theta^{0}, \theta^{0})\quad  \text {and}\quad   \psi_{tijlk} = \psi_{tijlk}(\theta^{0},\theta^{0}).$$
                  
\begin{rem}                      
If the process were not started at time $t=1$, it should be necessary to impose a causality and an invertibility condition, see for example Hallin \& Ingenbleek (1983) %and Hallin (1986). Note that $\psi_{t0k}(\theta^{0},\theta^{0})=0$, for $k \geq 1$.  %ALM_V13     
and Hallin (1986). Note that $\psi_{t0k}(\theta^{0},\theta^{0})=0$, for $k \geq 1$, and $1$, for $k=0$.  %ALM_V13     
\end{rem}  

%%%%%%%%%%%%%%%%%%%%%%%%%%%%%%%%%%%%%%%%%%%%%%%%%%%%%%%%%%%%%%%%%%%%%%%
%%%%%%%%%%%%%%%%%%%%%%%%%%%%%%%%%%%%%%%%%%%%%%%%%%%%%%%%%%%%%%%%%%%%%%%
\subsection{An illustration: tdVARMA($1, 1$)} \label{S2.4}
Let $\{ x_{t} : t\in \Nb \}$ be an $r$-vector time series satisfying 
$$
x_{t}  = A_{t}(\theta)x_{t-1} + e_{t}(\theta)+ B_{t}(\theta)e_{t-1}(\theta).  %
$$
It is easy to see that the tdVARMA(1,1) model is considered as a special case of the model defined in (\ref{31.12}) with $p=1$ and $q=1$. Now following M\'elard (1985) and %Azrak \& M\'elard (2012), in this special case, the coefficients of the pure moving average representation are given by: %ALM_V13
Azrak \& M\'elard (2015), in this special case, the coefficients of the pure moving average representation are given by: %ALM_V13
$$\psi_{tk}(\theta)=\Big\{\prod_{l=0}^{k-2}A_{t-l}(\theta)\Big\}
\Big\{B_{t-k+1}(\theta) + A_{t-k+1}(\theta)\Big\}, \quad \text{for} \quad k = 1, 2,..., t - 1,$$
%\end{equation}
where a product for $l = 0$ to $-1$ is set to $I_r$. The coefficients of the pure autoregressive form are 
%\begin{equation}%\nonumberight
$$\pi_{tk}(\theta)=\Big\{(-1)^{k-1}\prod_{l=0}^{k-2}B_{t-l}(\theta)\Big\}
\Big\{A_{t-k+1}(\theta) + B_{t-k+1}(\theta)\Big\},$$
%\end{equation}
so, for $i=1,...,m$, their derivatives are given by 
$$\frac{\partial\pi_{t1}(\theta)}{\partial\theta_{i}}=\frac{\partial A_{t}(\theta)}{\partial\theta_{i}} + \frac{\partial B_{t}
(\theta)}{\partial \theta_{i}}$$
$$\frac{\partial\pi_{t2}(\theta)}{\partial\theta_{i}} =-\frac{\partial B_{t}(\theta)}{\partial\theta_{i}}\{A_{t-1}(\theta) + B_{t-1}(\theta)\}
- B_{t}(\theta)\left\{\frac{\partial A_{t-1}(\theta)}{\partial\theta_{i}} + \frac{\partial B_{t-1}(\theta)}{\partial\theta_{i}}\right\}$$
\begin{eqnarray*}
\frac{\partial\pi_{t3}(\theta)}{\partial\theta_{i}}&=& \frac{\partial B_{t}(\theta)}{\partial\theta_{i}} B_{t-1}(\theta)\big\{A_{t-2}(\theta) + B_{t-2}(\theta)\big\} \nonumber \\
\nonumber &+& B_{t}(\theta)\frac{\partial B_{t-1}(\theta)}{\partial\theta_{i}}\big\{A_{t-2}(\theta) + B_{t-2}(\theta)\big\} \nonumber \\
\nonumber &+& B_{t}(\theta)B_{t-1}(\theta)\left\{\frac{\partial A_{t-2}(\theta)}{\partial\theta_{i}} + \frac{\partial B_{t-2}(\theta)}{\partial\theta_{i}}\right\},
... 
\end{eqnarray*}
Consequently
$$\frac{\partial\pi_{tk}(\theta)}{\partial\theta_{i}}=(-1)^{k-1}\sum_{l=1}^{k}\left(\prod_{h=1}^{k}
\chi_{t+1-h,k,l,h,i}(\theta)\right),$$
where
$$\chi_{t,k,l,h,i}(\theta)=\left\{
\begin{array}{lll}
  \frac{\partial\chi_{t,k,l,h}(\theta)}{\partial \theta_{i}}& {\rm if} & h=l , \\
  \chi_{t,k,l,h}(\theta) & {\rm if} & h \neq l , \\
\end{array}\right.$$
and
$$\chi_{t,k,l,h}(\theta) =\left\{\begin{array}{ll}
                            B_{t}(\theta) & $if$ \quad  h<k\\
                            A_{t}(\theta) + B_{t}(\theta) &  $if$ \quad h = k. \\
                          \end{array}\right.$$
Then  
$$\psi_{tik}(\theta, \theta^{0})=\sum_{u=1}^{k}\left\{\sum_{l=1}^{u}\left(\prod_{h=1}^{u}
\chi_{t+1-h,k,l,h,i}(\theta)\right)\right\}\left\{ \prod_{h=u+1}^{k}
\widetilde{\chi}_{t+1-h,k,h}(\theta^{0})  \right\},$$ %\text{To be checked }$$                       

$$\widetilde{\chi}_{t+1-h,k,h}(\theta) =\left\{\begin{array}{ll}
                            A_{t}(\theta)  & $if$ \quad  h < k , \\
                            A_{t}(\theta) + B_{t}(\theta) &  $if$ \quad h = k. \\
                          \end{array}\right.$$
In the univariate case these results can be found in Azrak \& M\'elard (2015, Chapter 4), correcting Azrak \& M\'elard (2006). 

%\red a voir par rapport aux exemples qui seront reellement faits; peut-etre rajouter une phrase que ces calulcs servent par apres pour verifier les hypotheses \black

%%%%%%%%%%%%%%%%%%%%%%%%%%%%%%%%%%%%%%%%%%%%%%%%%%%%%%%%%%%%%%%%%%%%%%%
%%%%%%%%%%%%%%%%%%%%%%%%%%%%%%%%%%%%%%%%%%%%%%%%%%%%%%%%%%%%%%%%%%%%%%%
\subsection{QMLE of tdVARMA($p$, $q$) models: asymptotic results}\label{S.4}
%%%%%%%%*********************************************************************************%%%%%%%%%%%%%%%%%%%%%%%%%%%%
% Let $\{x_{t} : t=1, 2, ..., n\}$ be a realization of the process $\{x_{t}\}$ of %ALM_V17
% length $n$ defined in (\ref{3.1}). In the present section, we shall apply the %ALM_V17
 Let $\{x_{t} : t=1, 2, ..., n\}$ be a partial realization of  length $n$ of the process $\{x_{t} : t \in \Nb \}$ %ALM_V17
 defined in (\ref{3.1}). In the present section, we shall apply the %ALM_V17
 general results of Section~\ref{S2.2} to the tdVARMA($p$, $q$) setting after formulating the minimal requirements for satisfying Assumptions $\Hb_{2.1}-\Hb_{2.4} $ (resp., $\Hb_{2.1}-\Hb_{2.5} $). The notations $Q_{n}(\theta), \alpha_{t}(\theta )$ as well as the QMLE solution~(\ref{22.3}) remain of course the same here.

%Theorem \ref{T1.3.1} below establishes the strong consistency of the QML estimators %ALM_V13
Theorem \ref{T1.3.1} below establishes the strong consistency of the QMLE %ALM_V13
 and further the asymptotic normality of this estimator. 
%To achieve that we will need the following assumptions.
For convenience, we suppose that the parameters in $A_{ti}(\theta)$ for $i=1,...,p$, in $B_{tj}(\theta)$ for $j=1,...,q$, and in $ g_{t}(\theta)$ are functionally independent. Without loss of generality we suppose that the vector $\theta$ is composed of three sub-vectors $A$, $B$ and $g$, more concretely $\theta = (A^{T}, B^{T}, g^{T})^{T}$, $A$ being the sub-vector of the parameters included in $A_{ti}(\theta)$ for  $i=1,...,p$, with dimension $s_{1}$, $B$  the sub-vector of the parameters included in $B_{tj}(\theta)$ for  $j=1,...,q$, with dimension $s_{2}$ and $g$  the sub-vector of the parameters included in $g_{t}(\theta)$ with dimension $m-s_{1}-s_{2}$.  Let us define the following Schur or Frobenius matrix norm.
\begin{defin} 
%The Schur or Frobenius norm, sometimes also called the Euclidean norm, is a matrix norm of an $m \times n$ matrix $A$ defined as  %ALM_V13
The Schur or Frobenius norm is a matrix norm of an $m \times n$ matrix $A$ defined as  %ALM_V13
$$ \left\|A\right\|_{F} = \sqrt{\tr\left(A^{T}A\right)}. $$
\end{defin} \vspace{2mm}
For further information about this matrix norm, see Golub \& Van Loan~(1996, p.~55). %The Frobenius norm can also be considered as a vector norm
% $\left\| \bullet \right\|$ already used in $\Hb_{2.4} $.  %ALM_V13

We now introduce a set of assumptions that will allow us, as mentioned above, to use results from Section~\ref{S2.2}. We assume for all $t \in \Nb$:
\begin{description}
	\item[$\Hb_{3.1}$]:  The matrices $A_{ti}(\theta)$, $B_{tj}(\theta)$ and $ g_{t}(\theta)$ are three times continuously differentiable with respect to $\theta$, in an open set $\Theta$ which contains the true value $\theta^{0}$ of $\theta$.  
	
	\item[$\Hb_{3.2}$]:  There exist positive constants $N_{1}$, $N_{2}$, $N_{3}$, $N_{4}$, $N_{5}$ and $0< \Phi < 1$ such that, for $\nu = 1,...,t - 1 $,  \\
	$$ \hspace{-2cm}	\sum_{k=\nu}^{t-1}\left\| \psi_{tik} \right\|^{2}_{F}<N_{1}\Phi^{\nu-1},\quad \sum_{k=\nu}^{t-1}\left\|\psi_{tik} \right\|^{4}_{F} <N_{2}\Phi^{\nu-1}, $$	
	$$ \hspace{-2cm}	\sum_{k=\nu}^{t-1}\left\|\psi_{tijk} \right\|^{2}_{F}<N_{3}\Phi^{\nu-1}, \quad \sum_{k=\nu}^{t-1}\left\|\psi_{tijk} \right\|^{4}_{F} <N_{4}\Phi^{\nu-1}, $$	
$$\hspace{-2cm}\sum_{k=1}^{t-1} \left\|\psi_{tijlk} \right\|^{2}_{F}<N_{5}, \quad i, j , l = 1,..., m,$$	

	\item[$\Hb_{3.3}$]:  There exist positive constants $K_{1}$, $K_{2}$, $K_{3}$, $K_{4}$, $K_{5}$ such that 
%	$$\hspace{-2.5cm}\left\|\left\{\frac{\partial \Sigma_{t}(\theta)}{\partial \theta_{i}}\right\}_{\theta=\theta^{0}}\right\|^{2}_{F}\leq K_{1}, %ALM_V13
	$$\hspace{-2.0cm}\left\|\left\{\frac{\partial \Sigma_{t}(\theta)}{\partial \theta_{i}}\right\}_{\theta=\theta^{0}}\right\|^{2}_{F}\leq K_{1}, %ALM_V13
	\quad \left\|\left\{\frac{\partial^{2} \Sigma_{t}(\theta)}{\partial \theta_{i}\partial\theta_{j}}\right\}_{\theta=\theta^{0}}\right\|^{2}_{F}\leq K_{2},\quad
%\left\|\left\{\frac{\partial^{3} \Sigma_{t}(\theta)}{\partial \theta_{i}\partial\theta_{j}\theta_{l}}\right\}_{\theta=\theta^{0}}\right\|^{2}_{F}\leq K_{3}, $$ %ALM_V13
\left\|\left\{\frac{\partial^{3} \Sigma_{t}^{-1}(\theta)}{\partial \theta_{i}\partial\theta_{j}\theta_{l}}\right\}_{\theta=\theta^{0}}\right\|^{2}_{F}\leq K_{3}, $$ %ALM_V13
$$\hspace{-2cm}\left\|\left\{\frac{\partial \Sigma^{-1}_{t}(\theta)}{\partial \theta_{i}}\right\}_{\theta=\theta^{0}}\right\|^{2}_{F}\leq K_{4},\quad \left\|\left\{\frac{\partial^{2}\Sigma^{-1}_{t}(\theta)}{\partial \theta_{i}\partial\theta_{j}}\right\}_{\theta=\theta^{0}}\right\|^{2}_{F}\leq K_{5}, \quad  i, j, l = 1,...,m.$$
	 
	\item[$\Hb_{3.4}$]:  There exist positive constants $M_{1}$,  $M_{2}$, and $ M_{3}$ such that 	
  $$ E\left[({\epsilon_{t}}^{T}\epsilon_{t})^{4}\right] \leq M_{1}, \quad \left\|E\left(\epsilon_{t}\epsilon^{T}_{t} \otimes \epsilon^{T}_{t}\right)\right\|_{F}\leq M_{2}, %\quad \left\|\Xi_{t}(\Sigma) \right\|_{F}\leq M_{3}, \quad \text{and } \quad \left\|\kappa_{t}(\Sigma)\right\|^{2}_{F}\quad \leq M_4, %ALM_V13
%	\quad \left\|\Xi_{t}(\Sigma) \right\|_{F}\leq M_{3},% \quad \text{and } \quad \left\|\kappa_{t}(\Sigma)\right\|^{2}_{F}\quad \leq M_4, %ALM_V13
$$
$$ %ALM_V13
\left\| \kappa_{t} \right\|_{F} %ALM_V13
+ \left\| \vecc(\Sigma).\vecc(\Sigma)^{T} \right\|_{F} %ALM_V13
+ \left\| \Sigma \otimes \Sigma \right\|_{F} %ALM_V13
+ \left\| K_{r,r}(\Sigma \otimes \Sigma) \right\|_{F} %ALM_V13
\leq M_{3}. %ALM_V13
$$ %ALM_V13
%where  %ALM_V13
%$$\Xi_{t}(\Sigma) = \kappa_{t}(\Sigma) - \vecc(\Sigma).\vecc(\Sigma)^{T} - (\Sigma \otimes \Sigma) - K_{rr}(\Sigma \otimes \Sigma). $$ %ALM_V13
% 	\item[$\Hb_{3.5}$]: There exist positive constants $m_1$, $m_2$, and $m_3$  such that  %ALM_V13
 	\item[$\Hb_{3.5}$]: There exist positive constants $m_1$ and $m_2$ such that  %ALM_V13
%$$ \left\|g_{t}(\theta^{0})\right\|^{4}_{F} \leq m_1 , \quad  \left\|\Sigma^{-1}_
%$$ \left\| g_{t} \right\|^{2}_{F} \leq m_1 , \quad  \left\|\Sigma^{-1}_{t} \right\|^{2}_{F}\leq m_2, \quad  \left\|\Sigma_{t} \right\|^{2}_{F}\leq m_{3}%,  %ALM_V13
$$ \left\| g_{t} \right\|^{2}_{F} \leq m_1 , \quad  \left\|\Sigma^{-1}_{t} \right\|^{2}_{F}\leq m_2.  %ALM_V13
%\quad  \left\|g^{T}_{t}(\theta^{0})\otimes g^{T}_{t}(\theta^{0})\right\|_{F}\leq m_4
%.  $$ %ALM_V13
$$ %ALM_V13

\noindent Furthermore, we suppose that 

	\item[$\Hb_{3.6}$]: 
\begin{eqnarray*}	
&&\hspace{-0.5 cm}
%\lim_{n\rightarrow\infty}\frac{1}{n}\sum_{t=1}^{n} \left(2 E_{\theta^{0}}\left(\frac{\partial e^{T}_{t}(\theta)}{\partial\theta_{i}} %ALM_V13
%\Sigma^{-1}_{t}(\theta)\frac{\partial e_{t}(\theta)}{\partial \theta_{j}}\right)\right.\\ %ALM_V13
%&&\left.\quad+ \tr \left[\Sigma^{-1}_{t} \frac{\partial \Sigma_{t}(\theta)}{\partial \theta_{i}} \Sigma^{-1}_{t} \frac{\partial \Sigma_{t}(\theta)}{\partial\theta_{j}} \right]_{\theta=\theta^{0}} \right)= V_{ij}(\theta^{0}), \quad   %ALM_V13
\lim_{n\rightarrow\infty}\frac{1}{n}\sum_{t=1}^{n} \left\{ E_{\theta^{0}}\left(\frac{\partial e^{T}_{t}(\theta)}{\partial\theta_{i}} %ALM_V13
\Sigma^{-1}_{t}(\theta)\frac{\partial e_{t}(\theta)}{\partial \theta_{j}}\right) \right. \\ %ALM_V13
&&\left. \quad+ \tr \frac{1}{2} \left[\Sigma^{-1}_{t} \frac{\partial \Sigma_{t}(\theta)}{\partial \theta_{i}} \Sigma^{-1}_{t} \frac{\partial \Sigma_{t}(\theta)}{\partial\theta_{j}} \right]_{\theta=\theta^{0}} \right\}= V_{ij}, \quad   %ALM_V13
\end{eqnarray*}   
%\mbox{a.s.} for $i, j =1,..., m$, where the matrix $V(\theta^{0})=(V_{1 \leq i,j \leq m}(\theta^{0}))_{i,j}$ is a strictly positive definite matrix. %ALM_V13
for $i, j =1,..., m$, where the matrix $V =(V_{i,j})_{1 \leq i,j \leq m}$ is a %strictly positive definite matrix. %ALM_V13  %ALM_V17
strictly positive definite matrix; %ALM_V13  %ALM_V17
	\item[$\Hb_{3.7}$]:
	$$
\frac{1}{n^{2}}\sum_{d= 1}^{n-1}\sum_{t=1}^{n-d} \sum^{t-1}_{k=1}\left\|g_{t-k} \right\|^{2}_{F}\left\| \psi^{}_{tik} \right\|_{F}\left\| \psi^{}_{t+d,i,k+d} \right\|_{F}=O\left(\frac{1}{n}\right),
$$
\begin{eqnarray*}
%	&& \frac{1}{n^{2}}\sum_{d=1}^{n-1}\sum_{t=1}^{n-d}\left[\sum^{t-1}_{k=1} M^{{jiT}}_{t0k,k})^{T} \Xi_{t} M^{ij}_{tdkk} \right.\\ %ALM_V13
	&& \frac{1}{n^{2}}\sum_{d=1}^{n-1}\sum_{t=1}^{n-d}\left[\sum^{t-1}_{k=1} M^{{jiT}}_{t0kk} \Xi_{t-k} M^{ij}_{tdkk} \right.\\ %ALM_V13
%&&	+ \underset{}{{\sum^{t-1}\limits_{k_{1} = 1}} {\sum^{t-1}\limits_{k_{2}=1}}} M^{jiT}_{t0k_{2}k_{1}} K_{rr} (\Sigma \otimes \Sigma) M^{ij}_{tdk_{1}k_{2}} \\ %ALM_V13
&&	+ \underset{}{{\sum^{t-1}\limits_{k_{1} = 1}} {\sum^{t-1}\limits_{k_{2}=1}}} M^{jiT}_{t0k_{2}k_{1}} K_{r,r} (\Sigma \otimes \Sigma) M^{ij}_{tdk_{1}k_{2}} \\ %ALM_V13
&&	+ \left. \sum^{t-1}_{k_{1} = 1} \sum^{t-1}_{k_{2}=1} M^{jiT}_{t0k_{2}k_{1}} (\Sigma \otimes \Sigma) M^{ij}_{tdk_{2}k_{1}} \right] =O\left(\frac{1}{n}\right),%ALM_V13
\end{eqnarray*}
with 	
%where  %ALM_V13
$$\Xi_{t}(\Sigma) = \kappa_{t}(\Sigma) - \vecc(\Sigma).\vecc(\Sigma)^{T} - (\Sigma \otimes \Sigma) - K_{r,r}(\Sigma \otimes \Sigma), $$ %ALM_V13
and, for $k', k'' = k, k_{1}, k_{2}$, %ALM_V13
%$$M^{ij}_{tfk_{1}k_{2}} = \vecc( g^{T}_{t-k_{1}} \psi^{T}_{t+f,i,k_{1}+f} \Sigma_{t+f}^{-1} \psi_{t+f,j,k_{2}+f} g_{t-k_{2}} ), \quad f=0, d.$$	 %ALM_V13
$$M^{ij}_{tfk'k''} = \vecc( g^{T}_{t-k'} \psi^{T}_{t+f,i,k'+f} \Sigma_{t+f}^{-1} \psi_{t+f,j,k''+f} g_{t-k''} ), \quad f=0, d.$$	 %ALM_V13
\end{description}

\begin{rem}
These assumptions are a generalization of those in Azrak \& M\'elard (2006). Note %however that {their} assumption  about upper bounds of the 4th power of the process is  %not needed  in their proof.  It will also not be used here, hence is left out. 
 however that their assumption  about upper bounds of the 4th order moment of the process is  %not needed  in their proof.  It will also not be used here, hence is left out. %ALM_V13  %ALM_V17
not needed in their proof.  It will also not be used here, hence it is left out. %ALM_V13 
\end{rem}

With these assumptions in hand, we are able to show that the conditions for Theorems~%\ref{T1.2.3}-~\ref{T1.2.5} hold (see Appendix A.1 for a detailed proof) and hence obtain the following result. %ALM_V13
\ref{T1.2.3}-~\ref{T1.2.5} hold (see Appendix A.1 for a sketch of the proof) and hence obtain the following result. %ALM_V13

\begin{theor} \label{T1.3.1}
Suppose that Assumptions $\Hb_{3.1}$-$\Hb_{3.7}$ hold. Then there exists a sequence of estimators  $\widehat{\theta}_{n} = (\widehat{\theta}_1 , ... , \widehat{\theta}_m )^{T}$ such that 
\begin{itemize}
	\item  $\widehat{\theta}_{n} \rightarrow \theta^{0}$ a.s., and  for every $\eps>0$ there exists an event $E$ with $P_{\theta^{0}}(E)> 1 -\eps $ and an $n_0$ such that, for $n >n_0$ on $E$, $Q_{n}(\theta )$ reaches a relative minimum at the point $\widehat{\theta}_{n}$;
	\item $ n^{1/2}(\widehat{\theta}_{n} - \theta^{0}) \stackrel{L}{\rightarrow}%\mathcal{N}(0,V(\theta^{0})^{-1}W(\theta^{0})V(\theta^{0})^{-1})$, with %ALM_V17
\mathcal{N}(0,V^{-1} W V^{-1})$, with %ALM_V17
%	  $$W(\theta^{0}) = \underset{n\rightarrow \infty }{\lim}\frac{1}{ 4  n}{\sum_{t=1}^{n}}E_{\theta^{0}}\left(\frac{\partial \alpha_{t}(\theta)}{\partial \theta}\frac{\partial \alpha_{t}(\theta)}{\partial \theta^{T}}\right). $$ %ALM_V17
	  $$W = \underset{n\rightarrow \infty }{\lim}\frac{1}{ 4  n}{\sum_{t=1}^{n}}E_{\theta^{0}}\left(\frac{\partial \alpha_{t}(\theta)}{\partial \theta}\frac{\partial \alpha_{t}(\theta)}{\partial \theta^{T}}\right). $$ %ALM_V17
\end{itemize}
\end{theor} 

\begin{rem}
For the sake of simplicity in the proof of Theorem \ref{T1.3.1} (see Appendix A.3),  we have made assumptions on $\Sigma_{t}$ in addition to those on $g_{t}$. The proof is %very similar to that of Azrak \& M\'{e}lard (2006) but is extended to multivariate %ALM_V13
somewhat similar to that of Azrak \& M\'{e}lard (2006) but is extended to multivariate %ALM_V13
% processes. Note however several corrections (e.g. the treatment of the third term of \eqref{dr3bis}) and improvements ($\theta^{0}$ and $\theta$ were sometimes not %distinguished). %ALM_V13
 processes. Note however several corrections with respect to that paper (e.g. the treatment of the third term of \eqref{dr3bis}; $\theta^{0}$ and $\theta$ were sometimes not distinguished where they should, especially in Section \ref{S1.3.1.2}) and improvements (the treatment of the last three terms of \eqref{dr3bis} is more detailed; also the role of the assumptions is better enlightened). %ALM_V13

\end{rem}

%%%%%%%%%%%%%%%%%%%%%%%%%%%%%%%%%%%%%%%%%%%%%%%%%%%%%%%%%%%%%%%%%%%%%%%
%%%%%%%%%%%%%%%%%%%%%%%%%%%%%%%%%%%%%%%%%%%%%%%%%%%%%%%%%%%%%%%%%%%%%%%
%%%%%%%%%%%%%%%%%%%%%%%%%%%%%%%%%%%%%%%%%%%%%%%%%%%%%%%%%%%%%%%%%%%%%%%
\section{Some examples}\label{S.5}

%\red on devrait dire un mot sur le but et la structure de cette section-ci avant de commencer \black
The two examples will show that the theory can be applied and that the assumptions can be verified. 

%%%%%%%%%%%%%%%%%%%%%%%%%%%%%%%%%%%%%%%%%%%%%%%%%%%%%%%%%%%%%%%%%%%%%%%
%%%%%%%%%%%%%%%%%%%%%%%%%%%%%%%%%%%%%%%%%%%%%%%%%%%%%%%%%%%%%%%%%%%%%%%
\subsection{Example 1: tdVAR(1) a generalization of Kwoun \& Yajima (1986) }\label{S.5.1}
%Kwoun and Yajima (1986) have investigated a non stationary AR model whose coefficients are deterministic functions of time, and they have given an tdAR(1) example defined as follows 
%\begin{eqnarray}
%x_{t} = A_{t}(\theta) x_{t-1}  + e_{t}(\theta), \label{SC.3}
%\end{eqnarray}
%where $A_{t}(\theta) = A' \sin(\alpha t + A'')$, $\alpha$ is a known constant, $A'$ and $A''$ are unknown parameters such that $A' \in [\epsilon, 1-\epsilon]$ and $A''\in [0, 2\pi-\epsilon]$ for $\epsilon >0$. So $\theta = (A', A'')$. As an illustration they have considered 
%\begin{eqnarray}
%x_{t} = 0.5 \sin(\frac{\pi}{2} t + 2)x_{t-1}  + \epsilon_{t}\label{KY1}
%\end{eqnarray}
%hence the true values of the parameters are $ A'^{0} = 0.5$ and $A''^{0} = 2$ with $\alpha = 1/2$. This example satisfies a number of conditions to obtaining a good asymptotic properties of the estimator. 
In the following example we discuss a generalization of Kwoun \& Yajima (1986). To achieve this we consider the bivariate tdVAR(1) model
%a constant innovation covariance matrix, in the second and third cases the innovations have a bounded time-dependent covariance matrix. 
\begin{eqnarray}
\left(\begin{array}{c} x_{t1}\\ x_{t2} \end{array}\right)= \left(\begin{array}{cc} A^{11}_{t} & A^{12}_{t} \\ A^{21}_{t} & A^{22}_{t}\end{array}\right) \left(\begin{array}{c} x_{t-1,1}\\ x_{t-1,2} \end{array}\right) + \left(\begin{array}{c} \epsilon_{t1}\\ \epsilon_{t2} \end{array}\right), \label{SC.6}
\end{eqnarray}  
where the coefficients $A^{ij}_{t} (\theta)$ for $i,j = 1,2$ are defined as 
\begin{eqnarray}
A^{ij}_{t} (\theta) = A'_{ij} \sin(\alpha_{ij} t + A''_{ij}).\label{SC.7}
\end{eqnarray} 
The unknown parameters $A_{ij}'$ and $A_{ij}''$  are such that $A_{ij}' \in [\delta, 1-\delta]$ and $A_{ij}''\in [0, 2\pi-\delta]$ for  some fixed $1/2>\delta >0$  and $\alpha_{ij}$ are known constants. Then 
\begin{eqnarray}
\theta = (A'_{11}, A'_{21}, A'_{12}, A'_{22}, A''_{11}, A''_{21}, A''_{12}, A''_{22})^{T}. \nonumber
\end{eqnarray}
The numerical example proposed by Kwoun \& Yajima (1986) for $r=1$ contains a process with periodic coefficients of period $4$ because $\alpha_{ij} = \pi/2$. However, it is well known, see \mbox{e.g.} Tiao \& Grupe (1980) and Azrak \& M\'elard (2006), that an $r$-dimensional autoregressive process with periodic coefficients of period $s \in \Nb $ can be embedded into an $s$-dimensional stationary autoregressive process. 
To avoid this { simplification} we consider  coefficients $A^{ij}_{t}(\theta)$ either with %distinct irrational periods or at least with large integral periods. We check 
distinct irrational periods or at least with large relatively prime periods. We check the assumptions of Theorem \ref{T1.3.1} in Appendix \ref{E.Ch.1} in the simplified case where we have 
\begin{eqnarray}
A_{t}(\theta) &=& \left(\begin{array}{cc}	A_{11}' \sin(at ) & \frac{1}{2} \\ 0 & A_{22}'\sin(bt) \end{array}\right),\label{btheor}
\end{eqnarray}
with $\theta = (A'_{11}, A'_{22})^{T}$. For the simulation study, we shall rather have recourse to the model
\begin{eqnarray}
A_{t}(\theta) &=& \left(\begin{array}{cc}	A_{11}' \sin(at ) & A_{12}' \\ 0 & A_{22}'\sin(bt) \end{array}\right),\label{b1.6}
\end{eqnarray}
with
\begin{eqnarray}
\theta = (A'_{11}, A_{12}',A'_{22})^{T}, \quad a = \frac{2\pi}{\sqrt{2499}}\quad\text{and}\quad b = \frac{2\pi}{\sqrt{2399}} .\label{b1.7}
\end{eqnarray} 
%We impose that $A'_{11}, A'_{22} \in [0, 1]$, partly similarly to Kwoun \& Yajima (1986), and write $\theta^{0} = (A^{'0}_{11}, 0.5, A^{'0}_{22})^{T}$ for the true value of $\theta$. For simplicity, we take here $\Sigma = I_{2}$. %ALM_V16
We impose that $A'_{11}, A'_{22} \in [0, 1[$, partly similarly to Kwoun \& Yajima (1986), and write $\theta^{0} = (A^{'0}_{11}, 0.5, A^{'0}_{22})^{T}$ for the true value of $\theta$. For simplicity, we take here $\Sigma = I_{2}$. %ALM_V16

\subsection{Example 2: tdVAR(1) with heteroscedasticity}\label{S.5.2}

Let us  re-consider the model defined in (\ref{SC.6})-(\ref{btheor}), except that $A_{12}'$ is no longer a parameter, with the added difficulty that the innovations are now $g_{t}\epsilon_{t}$ instead of $\epsilon_{t}$. Therefore we introduce a matrix $g_{t}(\theta)$ and we have a bounded time-dependent covariance matrix $\Sigma_{t}(\theta) = g_{t}(\theta)\Sigma g^{T}_{t}(\theta)$.  We have extended the Kwoun \& Yajima~(1986) parametrization by taking 
\begin{eqnarray}
g_{t}(\theta) &:=& \left(\begin{array}{cc}	\exp \left(-\eta_{11} \sin(ct)\right) & 1 \\ -1 & \exp \left(-\eta_{22} \sin(ct)\right) \\ \end{array}\right)\label{SC.821}
\end{eqnarray} 
with $c\in\R$. Also we use a matrix $\Sigma$ which is no longer the identity matrix:
$$
\Sigma := \left(\begin{array}{cc}	s_{11} & s_{12} \\ s_{12} & s_{22} \\ \end{array}\right).
$$
A close examination of $\Sigma_{t}(\theta)$ shows that if $g_{t}(\theta)$ were diagonal, then the correlation between the residuals would be constant, which is not very realistic for a time-dependent process. This is why we have put off-diagonal elements different from 0 in \eqref{SC.821}.   
%Note that similar examples are of course obtained for $\exp \left(-\eta_{ii} \cos(t)\right)$, $i=1,2$. 
Here, the vector $\theta$ reduces to %takes on the guise
\begin{eqnarray}
\theta = (A'_{11},  A'_{22}, \eta_{11},  \eta_{22})^{T}. \nonumber
\end{eqnarray}
The assumptions of Theorem \ref{T1.3.1} are checked in Appendix \ref{E.Ch.2}, and
simulation results  shown in Section \ref{S.6.2}.

%%%%%%%%%%%%%%%%%%%%%%%%%%%%%%%%%%%%%%%%%%%%%%%%%%%%%%%%%%%%%%%%%%%%%%%
%%%%%%%%%%%%%%%%%%%%%%%%%%%%%%%%%%%%%%%%%%%%%%%%%%%%%%%%%%%%%%%%%%%%%%%
%%%%%%%%%%%%%%%%%%%%%%%%%%%%%%%%%%%%%%%%%%%%%%%%%%%%%%%%%%%%%%%%%%%%%%%
\section{Simulation results}\label{S.5.4}

%%%%%%%%%%%%%%%%%%%%%%%%%%%%%%%%%%%%%%%%%%%%%%%%%%%%%%%%%%%%%%%%%%%%%%%
%%%%%%%%%%%%%%%%%%%%%%%%%%%%%%%%%%%%%%%%%%%%%%%%%%%%%%%%%%%%%%%%%%%%%%%
\subsection{Example 1: tdVAR(1) a generalization of Kwoun \& Yajima (1986) }\label{S.6.1}

The simulation experiment is performed in Matlab by using  the  program, which we call AJM, %described in Alj, Jonasson and M\'elard~(2015) and based on a special case of tdVAR(1)%ALM_V13
described in Alj  {et al.} (2015c) and based on a special case of tdVAR(1)
 process defined in (\ref{b1.6})-\eqref{b1.7},  with $A_{11}' =0.8$ and $A_{22}' = -0.9$ and $(\epsilon_{t1},\epsilon_{t2})^{T}$ %\red $\left(\begin{array}{c c} \epsilon_{t1} & \epsilon_{t2} \end{array}\right)^{T}$ \black
 has a bivariate normal distribution with covariance matrix $\Sigma = I_{2}$. A simulated series using these specifications is shown in Fig. 1. 
The true value of $\theta$ is  
\begin{eqnarray}
%\theta^{0} = (A'^{0}_{11}, A^{0}_{12}, A'^{0}_{22})^{T} = (0.8, 0.5, -0.9)^{T}. \nonumber %ALM_V17
\theta^{0} = (A'^{0}_{11}, A^{0}_{12}, A'^{0}_{22})^{T} = (0.8, 0.5, -0.9)^{T}, \nonumber %ALM_V17
\end{eqnarray}
We take 
\begin{eqnarray}
\theta^{i} = ( 0.1 \quad 0.1 \quad 0.1)^{T}. \nonumber
\end{eqnarray}
as initial value of  $\theta$.

%, in which lines (a) give the average of the parameter estimates, lines (b) give the averages across simulations of estimated standard errors of the corresponding estimates for the 1000 replications, lines (c) the sample standard deviations of the corresponding estimates for the 1000 replications and lines (d) give percentage of simulations where we reject the hypothesis $H_{0}(\theta_{i} = \theta^{0}_{i})$ with significance level 5\%.

\begin{table}[!t]
\caption {\footnotesize{Estimation results for the model \eqref{SC.6} under \eqref{b1.6}-\eqref{b1.7} via the program AMJ, where lines (a) give the averages of the parameter estimates, lines (b) give the averages across simulations of estimated standard errors of the corresponding estimates for the 1000 replications, lines (c) the sample standard deviations of the corresponding estimates for the 1000 replications and lines (d) give percentages of simulations where we reject the hypothesis $H_{0}(\theta_{i} = \theta^{0}_{i})$ at significance level 5\%. }}
\vspace{1cm}
\hspace{-0.5cm}
\begin{tabular}{ c c c c c c c c }  
\hline
\hline\\
%\vspace{0.5cm}
Sample size &  & \multicolumn{1}{c}{$\widehat{A'}_{11}$} & \multicolumn{1}{c}{$\widehat{A'}_{12}$} & \multicolumn{1}{c}{$\widehat{A'}_{22}$} & \multicolumn{1}{c}{$\widehat{\Sigma}_{11}$}& \multicolumn{1}{c}{$\widehat{\Sigma}_{12}$} & \multicolumn{1}{c}{$\widehat{\Sigma}_{22}$}  \tabularnewline \\
         &  & \multicolumn{1}{c}{$A'^{0}_{11}= 0.8$} & \multicolumn{1}{c}{$A'^{0}_{12}=0.5$} & \multicolumn{1}{c}{$A'^{0}_{22}=-0.9$}&\multicolumn{1}{c}{$\Sigma_{11}=1$}& \multicolumn{1}{c}{$\Sigma_{12}=0$} & \multicolumn{1}{c}{$\Sigma_{22}=1$}  \tabularnewline \\
%\hline
\hline  
%    &   &    &   &   &   &     & \tabularnewline     

%25 & (a) &    0.7481  &   0.5014    &     -0.8036     &     0.9804     &   -0.0026  &   1.0760 \tabularnewline %ALM_V13
25 & (a) &    0.7481  &   0.5014    &     $-0.8036$     &     0.9804     &   $-0.0026$  &   1.0760 \tabularnewline %ALM_V13
       
 & (b) &  0.2023    &    0.1543     &   0.2118     &  -  & - & - \tabularnewline
       
 & (c) &   0.2161   &  0.1654   &  0.2226   &  -  & - & - \tabularnewline
   
 & (d) &    7.1  &   7.7  &   4.8 &  -  & - & -\tabularnewline 
 
%   &   &    &   &   &   &   &    \tabularnewline         

%50 & (a) &  0.7714   &  0.5035  &   -0.8410     &    0.9845   &  0.0022  &    1.0392  \tabularnewline
50 & (a) &  0.7714   &  0.5035  &   $-0.8410$     &    0.9845   &  0.0022  &    1.0392  \tabularnewline
    
 & (b) &   0.1397  &  0.1049      &   0.1355    &   -  & - & - \tabularnewline

 & (c) &   0.1344  &   0.1112  &  0.1525    &  -  & - & - \tabularnewline
	       
 & (d) &  4.5  &  6  & 6.6 & -  & - & -\tabularnewline 
 
%   &   &    &   &   &   &   &      \tabularnewline  

%100 & (a) &    0.7855   &   0.4975    & - -0.8650    &     0.9927   &  0.0065 &   1.0224 \tabularnewline %ALM_V13
100 & (a) &    0.7855   &   0.4975    & $-0.8650$    &     0.9927   &  0.0065 &   1.0224 \tabularnewline %ALM_V13
           
 & (b) &   0.0963   &    0.0735    &   0.0926   &  -  & - & - \tabularnewline

 & (c) &   0.0995  &    0.0733    &   0.0964 &  -  & - & - \tabularnewline
    
 & (d) &  5.4  &     4.8  &   5  & -  & - & -\tabularnewline 
 
%   &   &    &   &   &   &   &      \tabularnewline    

%200 & (a) &   0.7905   &  0.4984     &   -0.8905   &     0.9959 &   0.0037     &  1.0087 \tabularnewline %ALM_V13
200 & (a) &   0.7905   &  0.4984     &   $-0.8905$   &     0.9959 &   0.0037     &  1.0087 \tabularnewline %ALM_V13

 & (b) &      0.0677  &     0.0510       &    0.0628    &  -  & - & - \tabularnewline

 & (c) &    0.0713 &    0.0513   &   0.0635 &    -  & - & - \tabularnewline
     
 & (d) &   5.4   &   5.7  &   4.2 & -  & - & -\tabularnewline 
 
%   &   &    &   &   &   &   &    \tabularnewline   
        
%400 & (a) &    0.7976    &   0.5000  &  -0.8932  &   0.9946 &  -0.0001  &     1.0073 \tabularnewline %ALM_V13
400 & (a) &    0.7976    &   0.5000  &  $-0.8932$  &   0.9946 &  $-0.0001$  &     1.0073 \tabularnewline %ALM_V13

 & (b) &   0.0474&     0.0358 &      0.0440  &  -  & - & - \tabularnewline

 & (c) &  0.0471  &    0.0363    &   0.0448    &  -  & - & - \tabularnewline

 & (d) &     5.2 &   4.7  &  5.1 & -  & - & -\tabularnewline 
 
%   &   &    &   &   &   &   &      \tabularnewline     
	
\hline
\end{tabular}
\label{Tab:SC.10} % labeling to refer it inside the text
\end{table}	
%\footnotemark[2] 
%\footnotetext[2]{AMJ: Alj, M\'elard and Jonasson Matlab program.}

%\newpage
The experiment was replicated 1000 times. The results are summarized in Table \ref{Tab:SC.10}.
 As the sample size becomes larger, we can see  that
\begin{itemize}
	\item the averages of the estimates become closer to their true value in accordance with the theory,  
	\item the sample standard deviation on line (b) becomes also closer to the averages across simulations of estimated standard errors in line (c) showing that the standard errors are well estimated, and
	\item the percentage of rejecting the hypothesis $H_{0}$ is close to 5\%.			
\end{itemize}
We compare, for the sample size $n=400$, a histogram of the 1000 replications of $\widehat{\theta}_{i}$  to the normal probability curve with mean equal to $\theta^{0}$ and standard deviation given in line (b) of Table~\ref{Tab:SC.10}. As we can see from the corresponding Figure \ref{hist_par_200}, this histogram shows empirically consistency and normality of the estimates.

%\newpage
%\begin{figure}[h!]
\begin{figure}%[h]
%\vspace{0.5cm}
%\hspace{-2.5cm} %ALM_V13
\hspace{-1cm} %ALM_V13
\scalebox{0.4}{% rescale the figure by a factor of 0.8}
\begin{minipage}{12cm}  
\includegraphics{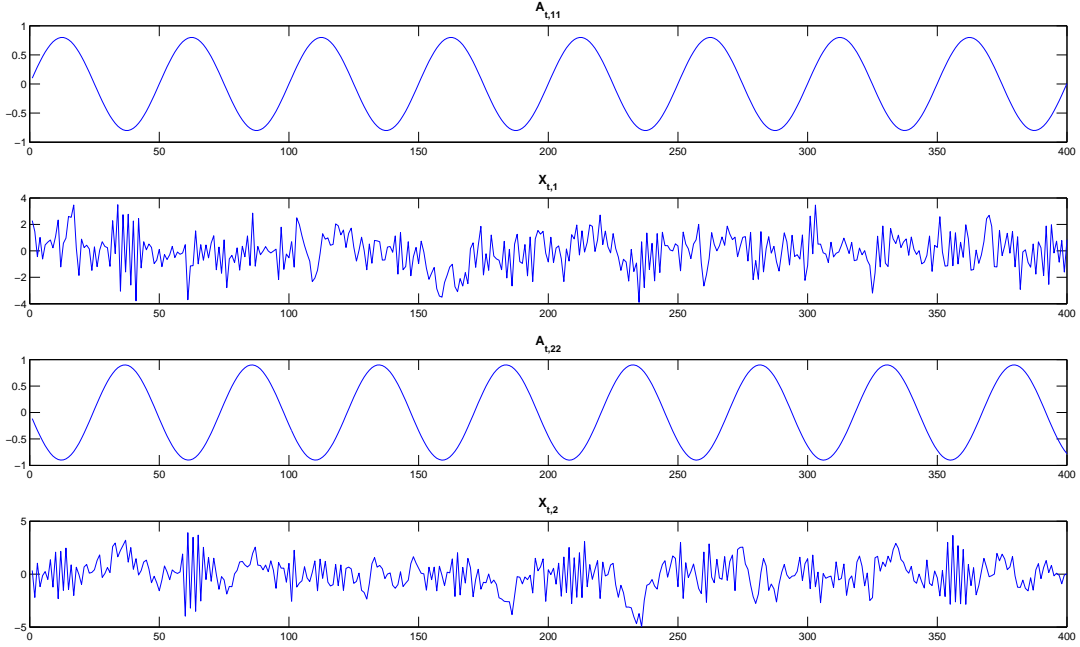}  
\end{minipage}
}
%\vspace{-0.5cm}
\caption{	\footnotesize{Time plots of the coefficients and the simulated tdVAR(1) generated by the process defined in (\ref{b1.6}) of length $n = 400$.}}
\label{Fig:ch.2.1} % labeling to refer it inside the text
\end{figure}
%\begin{landscape}
%\newpage
\begin{center}
%\begin{figure}[h!]
\begin{figure}%[h]
%\hspace{-2cm} %ALM_V13
\hspace{-1cm} %ALM_V13
\scalebox{0.4}{% rescale the figure by a factor of 0.8
% \begin{minipage}{10cm}  
\includegraphics{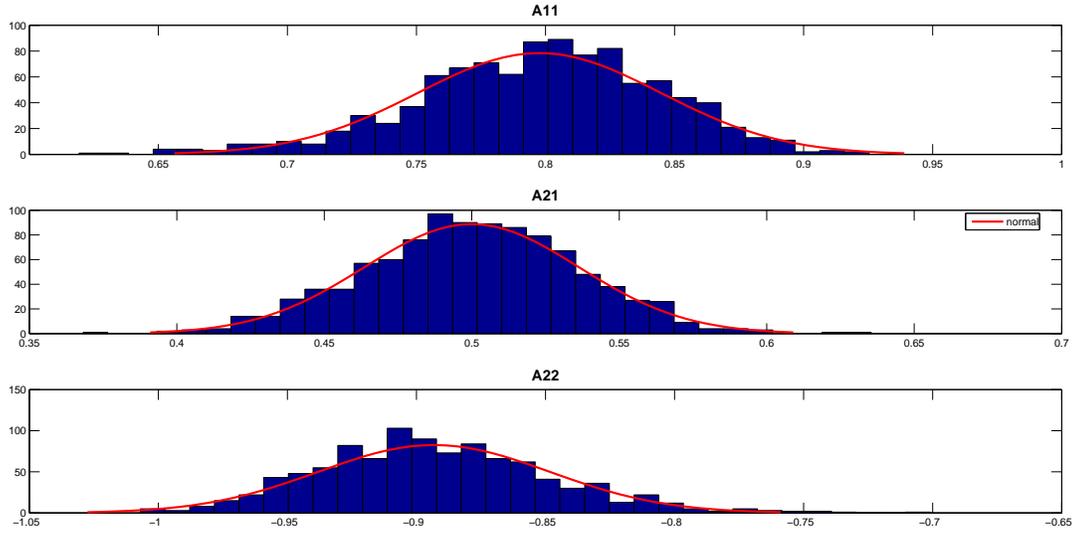} % importing figure
%\end{minipage}
} 
% \vspace{-1cm}
\caption{	\footnotesize{Histograms with the normal density function of 1000 replications of the parameters with $n = 400$.}}
\label{hist_par_200}  % labeling to refer it inside the text
\end{figure}
\end{center}
%\end{landscape}
%%%%%%%%%%%%%%%%%%%%%%%%%%%%%%%%%%%%%%%%%%%%%%%%%%%%%%%%%%%%%%%%%%%%%%%
%%%%%%%%%%%%%%%%%%%%%%%%%%%%%%%%%%%%%%%%%%%%%%%%%%%%%%%%%%%%%%%%%%%%%%%
\subsection{Example 2: tdVAR(1) with heteroscedasticity}\label{S.6.2}

We keep the bivariate model defined by (\ref{SC.6})-(\ref{btheor}), with the same numerical values for $a$ and $b$ but without $A'_{12}(\theta)$, with $g_t\epsilon_t$ instead of $\epsilon_t$ and a covariance matrix $\Sigma_{t}(\theta) = g_{t}(\theta)\Sigma g^{T}_{t}(\theta)$ bounded but time-dependent, where 
\begin{eqnarray}
g_{t}(\theta) := \left(\begin{array}{cc}	\exp \left(-\eta_{11} \sin(ct)\right) & 1 \\ -1 & \exp \left(-\eta_{22} \sin(ct)\right) \\ \end{array}\right), \label{SC.82}
\end{eqnarray} 
with $\Sigma_{11}=\Sigma_{22}=1$, $\Sigma_{12}=0.5$, $c=2\pi/25$ and $\eta_{11}=1$, $\eta_{22}=-1$ so that the innovation correlation coefficient varies between $-0.8$ and $0.8$.   

Again, the number of replications is $1000$. The results are presented in Table \ref{Tab:SC.11}. 
Moreover a program for computing the asymptotic information matrix on the basis of the formulae given in Appendix \ref{E.Ch.2} gave, for $n=50$ for example, the  standard errors $0.0905$, $0.908$, $0.1995$, $0.1995$, whereas the averages of the %standard errors estimated by the QMLE program were $0.963$, $0.1217$, $0.2027$, $0.1516$, respectively, and the standard deviations of the $1000$ estimates were $0.0917$, $0.1227$, $0.1879$, $0.1587$.  %ALM_V13
standard errors estimated by the QMLE program were $0.0963$, $0.1217$, $0.2027$, $0.1516$, respectively, and the standard deviations of the $1000$ estimates were $0.0917$, $0.1227$, $0.1879$, $0.1587$.  %ALM_V13

\begin{table}[!t]
\caption {\footnotesize{ Estimation results for the model \eqref{SC.6} under \eqref{b1.6}-\eqref{b1.7} and \eqref{SC.82}, via the program AMJ, where columns (a) give the averages of the parameter estimates and columns (d) give percentages of 1000 simulations where we reject the hypothesis $H_{0}(\theta_{i} = \theta^{0}_{i})$ at significance level 5\%. }}
\vspace{1cm}
\hspace{-0.5cm}
%\begin{tabular}{ c c c c c c c c }  
\begin{tabular}{ c c c c c c c c c c }  
\hline
\hline\\
%\vspace{0.5cm}
Sample size &  
& \multicolumn{2}{c}{$\widehat{A'}_{11}$ ($A_{11}^{0}=0.8$)} 
%& \multicolumn{1}{c}{$\widehat{A'}_{12}$} 
& \multicolumn{2}{c}{$\widehat{A'}_{22}$ ($A_{22}^{0}=-0.9$)} 
%& \multicolumn{1}{c}{$\widehat{\Sigma}_{11}$}
%& \multicolumn{1}{c}{$\widehat{\Sigma}_{12}$} 
%& \multicolumn{1}{c}{$\widehat{\Sigma}_{22}$}  
& \multicolumn{2}{c}{$\widehat{\eta}_{11}$ ($\eta_{11}^{0}=1.0$)} 
& \multicolumn{2}{c}{$\widehat{\eta}_{22}$ ($\eta_{22}^{0}=-1.0$)} 
\tabularnewline \\
& & (a) & (d) & (a) & (d) & (a) & (d) & (a) & (d) 
%         &  & \multicolumn{1}{c}{$A'^{0}_{11}= 0.8$} & \multicolumn{1}{c}{$A'^{0}_{12}=0.5$} & \multicolumn{1}{c}{$A'^{0}_{22}=-0.9$}&\multicolumn{1}{c}{$\Sigma_{11}=1$}& \multicolumn{1}{c}{$\Sigma_{12}=0$} & \multicolumn{1}{c}{$\Sigma_{22}=1$}  
\tabularnewline \\
%\hline
\hline  
 25 & & $0.7671$ & 3.8 & $-0.8567$ & 4.5 & $0.9897$ & 3.2 & $-0.9848$ & 5.1 \\
 50 & & $0.7808$ & 4.2 & $-0.8766$ & 5.0 & $0.9913$ & 2.9 & $-0.9964$ & 5.8 \\
100 & & $0.7910$ & 4.4 & $-0.8864$ & 5.8 & $0.9977$ & 4.1 & $-0.9975$ & 6.7 \\
200 & & $0.7963$ & 5.3 & $-0.8931$ & 5.1 & $0.9997$ & 5.3 & $-1.0000$ & 6.4 \\
400 & & $0.7972$ & 6.3 & $-0.8970$ & 4.2 & $0.9980$ & 4.7 & $-0.9984$ & 5.6 \\
\hline
\end{tabular}
\label{Tab:SC.11} % labeling to refer it inside the text
\end{table}	

\

\noindent ACKNOWLEDGMENTS\vspace{2mm}

\noindent Christophe Ley, who is also a member of ECARES, thanks the Fonds National de la Recherche Scientifique, Communaut\'e Fran\c caise de Belgique, for financial support via a Mandat de Charg\'e de Recherche FNRS. Guy M\'elard has benefited from a Belgian research grant F.R.S.-FNRS 1.5.261.09.

%\newpage
\makeatletter
\renewcommand{\@biblabel}[1]{}
\makeatother

%%%%%%%%%%%%%%%%%%%%%%%%%%%%%%%%%%%%%%%%%%%%%%%%%%%%%%%%%%%%%%%%%%%%%%%
%%%%%%%%%%%%%%%%%%%%%%%%%%%%%%%%%%%%%%%%%%%%%%%%%%%%%%%%%%%%%%%%%%%%%%%
%%%%%%%%%%%%%%%%%%%%%%%%%%%%%%%%%%%%%%%%%%%%%%%%%%%%%%%%%%%%%%%%%%%%%%%

\newpage
%\appendix
%\appendixpage
\begin{appendices}
%\numberwithin{equation}{section}
%\chapter{.\hspace{0.2cm}Proofs of Lemmas}
%\numberwithin{equation}{section}
%\chapter{.\hspace{0.2cm}Proofs of Lemmas}
%\numberwithin{equation}{section}
%%%%%%%%%%%%%%%%%%%%%%%%%%%%%%%%%%%%%%%%%%%%%%%%%%%%%%%%%%%%%%%%%%%%%%%
%%%%%%%%%%%%%%%%%%%%%%%%%%%%%%%%%%%%%%%%%%%%%%%%%%%%%%%%%%%%%%%%%%%%%%%
%%%%%%%%%%%%%%%%%%%%%%%%%%%%%%%%%%%%%%%%%%%%%%%%%%%%%%%%%%%%%%%%%%%%%%%
\section{Proof of theorems}\label{ASS0}

Since  we proceed as in Azrak \& M\'elard (2006),  we just sketch the %proof. Lemmas preceded by TA refer to the Technical Appendix, Alj \emph{\emph{et al.}} (2015a). %ALM_V13
proof, at least in Section \ref{Theor2}. Lemmas preceeded by TA refer to the Technical Appendix, Alj { {et al.}}~(2015a). %ALM_V13

%%%%%%%%%%%%%%%%%%%%%%%%%%%%%%%%%%%%%%%%%%%%%%%%%%%%%%%%%%%%%%%%%%%%%%%
%%%%%%%%%%%%%%%%%%%%%%%%%%%%%%%%%%%%%%%%%%%%%%%%%%%%%%%%%%%%%%%%%%%%%%%
\subsection{Proof of Theorem 2.1 and of Theorem 2.2}\label{Theor2}

As mentioned in Section~\ref{S2.2}, we check the following four assumptions of Klimko \& Nelson (1978) for $Q_{n}(\theta)$ defined by \eqref{22.1}.

%\paragraph{Assumption $H_{1.1}$}\label{H1_1} %ALM_V13
\paragraph{Assumption $\Hb_{1.1}$}\label{H1_1} %ALM_V13
$\quad n^{-1}\left\{\frac{\partial Q_{n}(\theta )}{\partial \theta _{i}}\right\}_{\theta = \theta^{0}}  \stackrel{\text{a.s.}}{\longrightarrow} 0$, for $i=1, ...,m$; \\

%\paragraph{Assumption $H_{1.2}$}\label{H1_2} %ALM_V13
\paragraph{Assumption $\Hb_{1.2}$}\label{H1_2} %ALM_V13
$\quad n^{-1}\left\{\frac{\partial^{2} Q_{n}(\theta )}{\partial \theta _{i}\theta _{j}}%\right\}_{\theta = \theta^{0}}  \stackrel{\text{a.s.}}{\longrightarrow} V_{ij}(\theta^{0})$, for $i,j=1, ..., m, $ where $V(\theta^{0})=(V_{ij}(\theta^{0}))_{1\leq i,j\leq m}$ is a strictly positive definite matrix of constants;\\ %ALM_V13
\right\}_{\theta = \theta^{0}}  \stackrel{\text{a.s.}}{\longrightarrow} V_{ij}$, for $i,j=1, ..., m, $ where $V =(V_{ij} )_{1\leq i,j\leq m}$ is a strictly positive definite matrix of constants;\\ %ALM_V13

%\paragraph{Assumption $H_{1.3}$}\label{H1_3} %ALM_V13
\paragraph{Assumption $\Hb_{1.3}$}\label{H1_3} %ALM_V13
$\quad \underset{n\rightarrow \infty }{\lim} \underset{\Delta \downarrow 0}{\sup}(n\Delta)^{-1}\left|\left\{\frac{\partial^{2} Q_{n}(\theta )}{\partial \theta_{i}\theta_{j}}\right\}_{\theta = \theta^{*}}- \left\{\frac{\partial^{2} Q_{n}(\theta )}{\partial \theta _{i}\theta _{j}}\right\}_{\theta = \theta^{0}} \right| < \infty \quad$ a.s., for $i, j =1, ...,m,$ 
where $\theta^{*}$ is a point of the straight line joining $\theta^{0}$ to $\theta$, %such that $\|\theta -\theta^{0} \|<\Delta$, $0 <\Delta$ and the vector norm $\left\|.\right\|$ is the $L_{2}$-norm. %ALM_V13
such that $\|\theta -\theta^{0} \|<\Delta$, $0 <\Delta$. %ALM_V13

%\paragraph{Assumption $H_{1.4}$}\label{H1_4} %ALM_V13
\paragraph{Assumption $\Hb_{1.4}$}\label{H1_4} %ALM_V13
$\quad\quad\quad n^{-1/2}\left\{\frac{\partial Q_{n}(\theta )}{\partial \theta}\right\}_{\theta = \theta^{0}} \stackrel{L}{\rightarrow} \mathcal{N}(0,W)$,	
%where $W$ is a ${m\times m}$ strictly positive definite matrix. %ALM_V13
where $W=(W_{ij} )_{1\leq i,j\leq m}$ is a strictly positive definite matrix. %ALM_V13

\begin{rem}
%Assumption $H_{1.3}$ coincides with $H_{2.4}$.  %ALM_V13
Assumption $\Hb_{1.3}$ coincides with $\Hb_{2.4}$.  %ALM_V13
\end{rem}

%%%%%%%%%%%%%%%%%%%%%%%%%%%%%%%%%%%%%%%%%%%%%%%%%%%%%%%%%%%%%%%%%%%%%%%
%\subsubsection{Proof of $H_{1.1}$}\label{PH1_1} %ALM_V13
\subsubsection{Proof of $\Hb_{1.1}$}\label{PH1_1} %ALM_V13
By TA Lemma 4.4, we have that $\{\partial \alpha_{t}(\theta)/\partial\theta_{i}, F_t \}$ is a martingale difference sequence.  Then we can use a strong law of large numbers for %martingale sequences (Stout, 1974) since $H_{2.1}$ implies the condition for it, more precisely %ALM_V13
martingale sequences (Stout, 1974, p. 154) since $\Hb_{2.1}$ implies the condition for it, more precisely %ALM_V13
\begin{eqnarray} 
 \sum_{t=1}^{\infty} \frac{E_{\theta^{0}}\left| \frac{\partial \alpha _{t}(\theta )}{\partial \theta_{i}} \right|^{p}}{t^{1+p/2}} < \infty, \label{Stout}
\end{eqnarray}
%for $p \geq 2$. This is true for $p = 4$ thanks to TA Lemma 4.11, proving that $H_{1.1}$ is fulfilled. %ALM_V13
for $p = 4$. %ALM_V13
$\hspace{1cm} \hfill \square $ %ALM_V13

%%%%%%%%%%%%%%%%%%%%%%%%%%%%%%%%%%%%%%%%%%%%%%%%%%%%%%%%%%%%%%%%%%%%%%%
%\subsubsection{Proof of $H_{1.2}$}\label{PH1_2}
\subsubsection{Proof of $\Hb_{1.2}$}\label{PH1_2}
By TA Lemma 4.5, we have $\{\partial^{2}\alpha_{t}(\theta)/\partial\theta_{i} \partial\theta_{j} - E_{ \theta} [ \partial^{2}\alpha_{t}(\theta)/\partial\theta_{i} \partial\theta_{j} ], F_{t} \}$ is a martingale difference sequence. Then we can again %use Stout (1974) strong law of large numbers, by adapting \eqref{Stout}, since $H_{2.2}$ implies the condition for it for $p = 2$, thanks to Lemma TA 4.12. Hence %ALM_V13
use the Stout (1974) strong law of large numbers, by adapting \eqref{Stout}, since $\Hb_{2.2}$ implies the condition for it for $p = 2$. Hence %ALM_V13
$$
\frac{1}{n}\sum_{t=1}^{n} \left\{ \frac{\partial^{2}\alpha _{t}(\theta)}{\partial\theta_{i}\partial\theta_{j}}\right\}_{\theta=\theta^{0}} - \frac{1}{n}\sum_{t=1}^{n}E_{\theta^{0}}\left(\frac{\partial^{2}\alpha_{t}(\theta)}{\partial\theta_{i}\partial\theta_{j}}/F_{t-1}\right)\stackrel{\text{a.s.}}{\longrightarrow}0, \quad i,j=1,...,m. 
$$
Then one half of the \mbox{a.s.} limit of the first term 
%is $n^{-1} \partial^{2} Q_{n}(\theta)/\partial \theta_{i} \partial \theta_{j}$ and  %ALM_V13
is $n^{-1} \partial^{2} Q_{n}(\theta)/\partial \theta_{i} \partial \theta_{j}$ for $\theta=\theta^{0}$ and  %ALM_V13
%defines $V_{ij}(\theta^{0})$, by $H_{2.3}$ so that $H_{1.2}$ is verified.   %ALM_V13
defines $V_{ij}$, by $\Hb_{2.3}$.   %ALM_V13
%\vspace{0.5cm} %ALM_V13
\vspace{0.5cm}$\hspace{1cm} \hfill \square $ %ALM_V13

%Up to now, all the assumptions of Theorem 2.1 is verified. To prove Theorem 2.2, there remains to check $H_{1.4}$.  %ALM_V13
Up to now, all the assumptions of Theorem 2.1 are verified. To prove Theorem 2.2, there remains to check $\Hb_{1.4}$.  %ALM_V13

%%%%%%%%%%%%%%%%%%%%%%%%%%%%%%%%%%%%%%%%%%%%%%%%%%%%%%%%%%%%%%%%%%%%%%%
%\subsubsection{Proof of $H_{1.4}$}\label{PH1_4} %ALM_V13
\subsubsection{Proof of $\Hb_{1.4}$}\label{PH1_4} %ALM_V13
We proceed by using the Central Limit Theorem for martingale difference sequences of Basawa \& Prakasa Rao (1980). Using the Cram\'{e}r-Rao device with any vector $\lambda$, this requires to prove that, for $\xi_{t} = \lambda^{T} [ \partial \alpha_{t}(\theta)/\partial\theta ]_{\theta = \theta^{0}}$, $E|\xi_{t}|^{4}$ is bounded. 

%\noindent This is done again using $H_{2.1}$.  %ALM_V13
%Then, with $W$ defined by $H_{2.5}$, asymptotic normality follows with the stated covariance matrix. $\hspace{1cm} \hfill \square $ %ALM_V13
\noindent This is done again using $\Hb_{2.1}$.  %ALM_V13
Then, with $W$ defined by $\Hb_{2.5}$, asymptotic normality follows with the stated covariance matrix. $\hspace{1cm} \hfill \square $ %ALM_V13

\subsection{Further preliminaries}\label{Further_preliminaries}
In order to prove Theorem 3.1, we need the following Lemma, due to Hamdoune (1995) (see also Azrak \& M\'elard 2006), and a strong law of large numbers for mixingale sequences. 

\begin{Lem} \label{LemmaA1}
Let $\{w_{t}^{}, t = 1,..., n\}$ be, for each $n \in N$, a scalar process with finite second-order moments, i.e.

\noindent \textrm{Lemma \ref{LemmaA1}\text{-i}}: $E( w_{t}^{2}) < \infty$

\noindent \textrm{Lemma \ref{LemmaA1}\text{-ii}}: $E\left(n^{-1}\sum_{t=1}^{n}w_{t}^{2}\right) = O(n^{-\delta}) \quad \text{with}\quad \delta >0.$

\noindent Then, $n^{-1}\sum_{t=1}^{n} w_{t}^{}$ converges \mbox{a.s.} to zero when $n$ tends to infinity.
\end{Lem} 
We also need a strong law of large numbers for mixingale sequences, \mbox{e.g.}  Hall \& Heyde (1980, Theorem 2.21) in the special case where their sequence $b_{n} = n$. Let us recall the definition from Hall \& Heyde (1980, Section 2.3). 

\begin{defin}\label{Mixingale}
Let $\{w^{}_{t}, t \geq 1\}$ be square-integrable random variables on a probability space $(\Omega, F, P)$ and $\{F_{t}, -\infty < t< \infty\}$ be an increasing sequence of $\sigma$-fields of $F$. Then $\{w^{}_{t}, F_{t}\}$ is a $L_{2}$-mixingale sequence if for sequences of nonnegative constants 
$\psi_{\nu}$ and $c_{t}$ where $\psi_{\nu} \rightarrow 0$ as $\nu \rightarrow \infty$, we have

\noindent i. $E \{E(w_{t} | F_{t-\nu} )^{2} \} \leq \psi_{\nu} c_{t}$, and

\noindent ii. $E( w_{t} - E(w_{t} | F_{t+\nu} ))^{2} \leq \psi_{\nu+1} c_{t}$
. 
\end{defin}

\begin{Lem} \label{LemmaA2}
If $\{w^{}_{t}, F_{t}\}$ is a $L_{2}$-mixingale sequence, and if $\sum_{t=1}^{n} c_{t}^{2} < \infty$ and $\psi_{n} = O(n^{-1/2}(\log n)^{-2})$ as $n \rightarrow \infty$, then $n^{-1} \sum_{t=1}^{n} w_{t} \stackrel{\text{a.s.}}{\rightarrow}0$ as $n \rightarrow \infty$. 
\end{Lem} 

%%%%%%%%%%%%%%%%%%%%%%%%%%%%%%%%%%%%%%%%%%%%%%%%%%%%%%%%%%%%%%%%%%%%%%%
%%%%%%%%%%%%%%%%%%%%%%%%%%%%%%%%%%%%%%%%%%%%%%%%%%%%%%%%%%%%%%%%%%%%%%%
\subsection{Proof of Theorem 3.1}\label{S1.51}
%\paragraph{Consistancy}
%\begin{rem}
%\end{rem}
First of all, as shown in Section \ref{S1.3.1.2}, assumption $\Hb_{3.1}$ is used to define the $\psi$'s used in $\Hb_{3.2}$ and in the derivatives in the other assumptions. %ALM_V13
The idea is to check the five assumptions of Theorems 2.1 and 2.2. 
%First $H_{2.1}$ and $H_{2.2}$ are direct consequences of TA Lemma 2.1 and TA Lemma 2.2, respectively.  %ALM_V13
Then $\Hb_{2.1}$ and $\Hb_{2.2}$ are direct consequences of TA Lemma 2.1 (using TA Lemma 4.11) and TA Lemma 2.2 (using TA Lemma 4.12), respectively. This makes use of assumptions $\Hb_{3.2}$, $\Hb_{3.3}$, $\Hb_{3.4}$ and $\Hb_{3.5}$. %ALM_V13 

%%%%%%%%%%%%%%%%%%%%%%%%%%%%%%%%%%%%%%%%%%%%%%%%%%%%%%%%%%%%%%%%%%%%%%%
%\subsubsection{Proof of $H_{2.3}$}\label{PH2_3} %ALM_V13
\subsubsection{Proof of $\Hb_{2.3}$}\label{PH2_3} %ALM_V13

%Let us consider the process $Z^{}_{t,ij}$ defined by %ALM_V13
Let us consider the process $Z^{}_{tij}$ defined by %ALM_V13
\begin{eqnarray}
%Z^{}_{t,ij} &=& \left\{ \frac{\partial e^{T}_{t}(\theta)}{\partial  \theta_{i}}\Sigma^{-1}_{t}( \theta) %ALM_V13
Z^{}_{tij} &=& \left\{ \frac{\partial e^{T}_{t}(\theta)}{\partial  \theta_{i}}\Sigma^{-1}_{t}( \theta) %ALM_V13
\frac{\partial e^{}_{t}( \theta)}{\partial \theta_{j}} \right\}_{\theta=\theta^{0}} \nonumber - E_{\theta^{0}}\left(\frac{\partial e^{T}_{t}( \theta)}{\partial  \theta_{i}}\Sigma^{-1}_{t}( \theta)
\frac{\partial e^{}_{t}( \theta)}{\partial \theta_{j}}\right).\label{6.30bis} 
\end{eqnarray}
%By using TA Lemma 4.13, the two assumptions i and ii of Lemma A.1 are fulfilled for $Z^{}_{t,ij}$, hence  %ALM_V13
By using TA Lemma 4.13 and assumptions $\Hb_{3.2}$, $\Hb_{3.4}$, $\Hb_{3.5}$ and $\Hb_{3.7}$, the two assumptions i and ii of Lemma A.1 are fulfilled for $Z^{}_{tij}$, hence  %ALM_V13
\begin{eqnarray}
\frac{1}{n}\sum_{t=1}^{n}\left[\left\{
\frac{ \partial e^{T}_{t}(\theta)}{\partial \theta_{j}} \Sigma^{-1}_{t}(\theta)
 \frac{\partial e^{}_{t}( \theta)}{\partial \theta_{j}}\right\}_{ \theta = \theta^{0}} 
- E_{\theta^{0}}\left(\frac{ \partial e^{T}_{t}(\theta)}{\partial \theta_{j}} \Sigma^{-1}_{t}( \theta) \frac{\partial e^{}_{t}( \theta)}{\partial \theta_{j}}\right)\right]\stackrel{\text{a.s.}}{\longrightarrow} 0 \nonumber \\
\label{6.30ter} 
\end{eqnarray}
when $n\rightarrow\infty$. However, from TA Lemma 4.5, we have for $i,j = 1,...,m$:
\begin{eqnarray}
\frac{1}{2n}\sum_{t=1}^{n}E_{\theta^{0}}\left(\frac{\partial^{2} \alpha^{}_{t}( \theta)}{\partial \theta_{i}\partial \theta_{j}}/F_{t-1}\right) &=&\frac{1}{n}\sum_{t=1}^{n} \left\{ \frac{\partial e^{T}_{t}( \theta)}{\partial  \theta_{i}}\Sigma^{-1}_{t}( \theta) 
\frac{\partial e^{}_{t}( \theta)}{\partial \theta_{j}}\right\}_{ \theta =  \theta^{0}}\nonumber 
\end{eqnarray}
\vspace{-0.5cm}
\begin{eqnarray}
\hspace{-.5cm}
&& + \frac{1}{2n}\sum_{t=1}^{n} \tr \left\{ \Sigma^{-1}_{t}(\theta) \frac{\partial \Sigma^{}_{t}(\theta)}{\partial  \theta_{i}} \Sigma^{-1}_{t}(\theta) \frac{\partial \Sigma^{}_{t}(\theta)}{\partial \theta_{j}}\right\}_{ \theta =  \theta^{0}}. \label{6.30} 
\end{eqnarray}
Then \eqref{6.30ter} implies that the \mbox{a.s.} limit of (\ref{6.30}) for $n\rightarrow\infty$ will be equal to
\begin{eqnarray}
&&\lim_{n\rightarrow\infty}\frac{1}{n}\sum_{t=1}^{n}E_{\theta^{0}}\left(\frac{\partial e^{T}_{t}(\theta)}{\partial \theta_{i}}
\Sigma^{-1}_{t}( \theta) \frac{\partial e^{}_{t}( \theta)}{\partial \theta_{j}}\right) \nonumber \\
&&+ \lim_{n\rightarrow\infty} \frac{1}{2n}\sum_{t=1}^{n}\tr \left\{ \Sigma^{-1}_{t}(\theta) \frac{\partial \Sigma^{}_{t}(\theta)}{\partial  \theta_{i}} \Sigma^{-1}_{t}( \theta) \frac{\partial \Sigma^{}_{t}(\theta)}{\partial \theta_{j}}\right\}_{\theta =\theta^{0}}, \nonumber
\end{eqnarray}
%and this is $V_{ij}(\theta^{0})$ as defined in $H_{3.6}$.  %ALM_V13
and this is $V_{ij}$ as defined in $\Hb_{3.6}$.  %ALM_V13

%%%%%%%%%%%%%%%%%%%%%%%%%%%%%%%%%%%%%%%%%%%%%%%%%%%%%%%%%%%%%%%%%%%%%%%
%\subsubsection{Proof of $H_{2.4}$}\label{PH2_4} %ALM_V13
\subsubsection{Proof of $\Hb_{2.4}$}\label{PH2_4} %ALM_V13
By using TA Lemma 2.3, it suffices to show that
\begin{eqnarray*}
\lim_{n\rightarrow\infty}\frac{1}{n}\left|\sum_{t=1}^{n}\left\{\frac{\partial^{3}\alpha_{t}(\theta)}{\partial \theta_{i}\partial \theta_{j}\partial
\theta_{l}}\right\}_{\theta = \theta^{0}}\right| &<& \infty \quad \quad {\rm a.s.} \quad {\rm for} \quad i, j, l = 1,..., m. \label{H6.2}
\end{eqnarray*}
%As shown in TA Lemma 2.4, the expression in the left hand side can be bounded by %ALM_V13
As a consequence of TA Lemma 2.4, the expression on the left-hand side can be bounded by %ALM_V13
\begin{eqnarray}
%\hspace{-2cm}
%\tilde{\Phi}_{1} + \frac{1}{n}\tilde{\Phi}_{2}\sum_{t=1}^{n}\left| \epsilon^{T}_{t}\epsilon_{t}\right| + \frac{1}{n}\sum_{t=1}^{n} \tilde{\Psi}_{1t} + \frac{1}{n}\sum_{t=1}^{n}\tilde{\Psi}_{2t}+ \frac{1}{n}\sum_{t=1}^{n}\tilde{\Psi}_{3t}, \label{dr3bis}
%\end{eqnarray} %ALM_V13
\tilde{\Phi}_{1} + \frac{1}{n} \sum_{t=1}^{n} \tilde{\Phi}_{2t} + \frac{1}{n}\sum_{t=1}^{n} \tilde{\Psi}_{1t} + \frac{1}{n}\sum_{t=1}^{n}\tilde{\Psi}_{2t}+ \frac{1}{n}\sum_{t=1}^{n}\tilde{\Psi}_{3t}, \label{dr3bis}
\end{eqnarray} %ALM_V13
%where $\tilde{\Phi}_{1} $ and $\tilde{\Phi}_{2}$ are defined in TA Lemmas 4.14 and 4.15 and the last three terms $\tilde{\Psi}_{1t}$, $\tilde{\Psi}_{2t}$, and $\tilde{\Psi}_{3t}$ are briefly described now.   %ALM_V13
where $\tilde{\Phi}_{1} $ is shown in TA Lemma 4.14 (using assumptions $\Hb_{3.3}$ and $\Hb_{3.5}$) to be bounded, and the last four terms $\tilde{\Phi}_{2t}$, $\tilde{\Psi}_{1t}$, $\tilde{\Psi}_{2t}$, and $\tilde{\Psi}_{3t}$ are briefly described now:  %ALM_V13
\begin{itemize}
\item $\tilde{\Phi}_{2t}$ contains terms like  %ALM_V13
$e^{T}_{t}(\theta)$ $(\partial^{3} \Sigma^{-1}_{t}(\theta) /\partial\theta_{i}\partial\theta_{j}\partial \theta_{l})$ $e^{}_{t}(\theta)$;   %ALM_V13
\item $\tilde{\Psi}_{1t}$ contains terms like 
$e^{T}_{t}(\theta)$ $(\partial^{2} \Sigma^{-1}_{t}(\theta) /\partial\theta_{i}\partial\theta_{j})$ $(\partial e^{}_{t}(\theta)/\partial \theta_{l})$, 
$e^{T}_{t}(\theta)$ $(\partial \Sigma^{-1}_{t}(\theta) /\partial\theta_{i})$ $(\partial e^{2}_{t}(\theta)/\partial \theta_{j} \partial \theta_{l})$, 
and 
$e^{T}_{t}(\theta)$ $\Sigma^{-1}_{t}(\theta)$ $(\partial e^{3}_{t}(\theta)/\partial \theta_{i} \partial \theta_{j} \partial \theta_{l})$, where $i, j, l = 1,...,m$;
\item the last two terms  $\tilde{\Psi}_{2t}$ and $\tilde{\Psi}_{3t}$ 
contain respectively terms like 
$(\partial e^{T}_{t}(\theta)/\partial\theta_{i})$ $(\partial \Sigma^{-1}_{t}(\theta) /\partial\theta_{l})$ $(\partial e^{}_{t}(\theta)/\partial \theta_{j})$ and 
$(\partial^{2}e^{T}_{t} (\theta)/\partial \theta_{i}\partial \theta_{j})$ $\Sigma^{-1}_{t}(\theta)$ $(\partial e^{}_{t}(\theta)/\partial \theta_{l})$, where $i, j, l = 1,...,m$.  %ALM_V13
\end{itemize}
For $i, j, l = 1,...,m$, let us define the following six sequences of random variables 
\begin{eqnarray*}
\hspace{-3cm}
X^{ilj}_{t}&=&  \left(\frac{\partial e_{t}^{T}(\theta)}{\partial\theta_{i}} \frac{\partial \Sigma^{-1}_{t}(\theta)}{\partial\theta_{l}} \frac{\partial e_{t}^{}(\theta)}{\partial\theta_{j}}\right)_{\theta = \theta^{0}} 
- E_{\theta^{0}}\left[ \frac{\partial e_{t}^{T}(\theta)}{\partial\theta_{i}} \frac{\partial \Sigma^{-1}_{t}(\theta)}{\partial\theta_{l}} \frac{\partial e_{t}^{}(\theta)}{\partial\theta_{j}}\right] ,\label{h1.41} \\
%\end{eqnarray}
%\begin{eqnarray}
Y^{ijl}_{t}&=&  \left(\frac{\partial e_{t}^{T}(\theta)}{\partial\theta_{i}\partial\theta_{j}} \Sigma^{-1}_{t}(\theta) \frac{\partial e_{t}^{}(\theta)}{\partial\theta_{l}}\right)_{\theta = \theta^{0}} 
- E_{\theta^{0}}\left[ \frac{\partial e_{t}^{^{T}}(\theta)}{\partial\theta_{i}\partial\theta_{j}} \Sigma^{-1}_{t}(\theta) \frac{\partial e_{t}^{}(\theta)}{\partial\theta_{l}}\right], \label{h1.41bis} \\
%\end{eqnarray}
%Z^{}_{t}&=&  e_{t}^{^{T}}(\theta) e_{t}^{}(\theta), \label{h1.41quater} \\ %ALM_V13
Z^{}_{t}&=& \left( e_{t}^{T}(\theta) \frac{\partial^{3} \Sigma^{-1}_{t}(\theta)}{\partial\theta_{i}\partial\theta_{j}\partial\theta_{l}} e_{t}^{}(\theta) \right)_{\theta = \theta^{0}}  %ALM_V13
- E_{\theta^{0}}\left[ e_{t}^{T}(\theta) \frac{\partial^{3} \Sigma^{-1}_{t}(\theta)}{\partial\theta_{i}\partial\theta_{j}\partial\theta_{l}} e_{t}^{}(\theta) \right], 
\label{h1.41quater} \\ %ALM_V13
%W^{ijl}_{1t}&=&  \left( e_{t}^{T}(\theta) \frac{\partial^{2} \Sigma^{-1}_{t}(\theta)}{\partial\theta_{i}\partial\theta_{j}} \frac{\partial e_{t}^{}(\theta)}{\partial\theta_{l}} \right)_{\theta = \theta^{0}}  %ALM_V13
W^{(1)ijl}_{t}&=&  \left( e_{t}^{T}(\theta) \frac{\partial^{2} \Sigma^{-1}_{t}(\theta)}{\partial\theta_{i}\partial\theta_{j}} \frac{\partial e_{t}^{}(\theta)}{\partial\theta_{l}} \right)_{\theta = \theta^{0}}  %ALM_V13
- E_{\theta^{0}}\left[ e_{t}^{T}(\theta) \frac{\partial^{2} \Sigma^{-1}_{t}(\theta)}{\partial\theta_{i}\partial\theta_{j}} \frac{\partial e_{t}^{}(\theta)}{\partial\theta_{l}} \right], \label{h1.41ter}
%\begin{eqnarray}
\end{eqnarray*}
%and $W^{ijl}_{2t}$ and $W^{ijl}_{3t}$ defined similarly by replacing $(\partial^{2} \Sigma^{-1}_{t}(\theta)/\partial\theta_{i} \partial\theta_{j}) (\partial e_{t}^{}(\theta)/\partial\theta_{l})$ by, respectively,  %ALM_V13
and $W^{(2)ijl}_{t}$ and $W^{(3)ijl}_{t}$ defined similarly by replacing $(\partial^{2} \Sigma^{-1}_{t}(\theta)/\partial\theta_{i} \partial\theta_{j}) (\partial e_{t}^{}(\theta)/\partial\theta_{l})$ with, respectively,  %ALM_V13
$(\partial^{} \Sigma^{-1}_{t}(\theta)/\partial\theta_{i}) (\partial e_{t}^{2}(\theta)/ \partial\theta_{j} \partial\theta_{l})$ 
and $\Sigma^{-1}_{t}(\theta) (\partial e_{t}^{3}(\theta)/ \partial\theta_{i} \partial\theta_{j}\partial\theta_{l})$. 

%TO BE CHANGED Let $Y_{t}$ be defined by
%$$Y_{t}= {\epsilon}^{T}_{t}\epsilon_{t} - E({\epsilon}^{T}_{t}\epsilon_{t}).$$ 
%First, we have that $E(Z_{t}/ F_{t-1}) = 0$ and $E(Z_{t}^{2})$ is uniformly bounded by a constant $M_{1}{1/2}$, using $H_{3.4}$ and Cauchy-Schwarz inequality. Then, %under Lemma \ref{L3.2.1} %ALM_V13
First, we have that $E(Z_{t}/ F_{t-1}) = 0$ and, according to TA Lemma 4.15 (using assumptions $\Hb_{3.3}$, $\Hb_{3.4}$ and $\Hb_{3.5}$) , that $E(Z_{t}^{2})$ is uniformly bounded by a constant. Then, %under Lemma \ref{L3.2.1} %ALM_V13
the strong law of large numbers (Stout, 1974, p. 154) implies that 
$$
%\lim_{n\rightarrow\infty}\frac{1}{n}\sum_{t=1}^{n}Y_{t}&=&\lim_{n\rightarrow\infty}
\frac{1}{n}\sum_{t=1}^{n}Z_{t} \stackrel{{\rm a.s.}}{\rightarrow}0. %ALM_V13
$$
%The arguments for the other sequences is more involved. Using expressions in TA Lemmas 4.14 and 4.15, the coefficients  %ALM_V13
%$\tilde{\Phi}_{1} $ and $\tilde{\Phi}_{2}$ in the first two terms of \eqref{dr3bis} can be bounded.  %ALM_V13
The arguments for the other sequences are more involved. %ALM_V13
%From TA Lemmas 4.17, 4.19 and 4.21 we have that $\{W^{ijl}_{qt}, F_{t}\}$, $q=1, 2, 3$, $\{X^{ilj}_{t}, F_{t}\}$ and $\{Y^{ijl}_{t}, F_{t}\}$ are $L_{2}$-mixingale sequences.  %ALM_V13
From TA Lemmas 4.16, 4.18 and 4.20 we have that $\{W^{(q)ijl}_{t}, F_{t}\}$, $q=1, 2, 3$, $\{X^{ilj}_{t}, F_{t}\}$ and $\{Y^{ijl}_{t}, F_{t}\}$ are $L_{2}$-mixingale sequences.  %ALM_V13
%From TA Lemmas 4.18, 4.20 and 4.22 we have that these $L_{2}$-mixingale sequences $\{W^{ijl}_{qt}, F_{t}\}$, $q=1, 2, 3$, $\{X^{ilj}_{t}, F_{t}\}$ and $\{Y^{ijl}_{t}, F_{t}\}$ fulfil the conditions in Lemma A.2, the strong law of large numbers for a mixingale sequence (Hall and Heyde, 1980, p. 41, Theorem 2.21). Hence  %ALM_V13
From TA Lemmas 4.17, 4.19 and 4.21 we have that these $L_{2}$-mixingale sequences $\{W^{(q)ijl}_{t}, F_{t}\}$, $q=1, 2, 3$, $\{X^{ilj}_{t}, F_{t}\}$ and $\{Y^{ijl}_{t}, F_{t}\}$ fulfil the conditions in Lemma A.2, the strong law of large numbers for a mixingale sequence (Hall \& Heyde, 1980, p. 41, Theorem 2.21). This makes use of assumptions $\Hb_{3.2}$, $\Hb_{3.3}$, $\Hb_{3.4}$ and $\Hb_{3.5}$. Hence  %ALM_V13
\begin{eqnarray*}
%n^{-1}\sum_{t=1}^{n}W^{ijl}_{1t} \stackrel{\text{a.s.}}{\rightarrow}0 ,\quad %ALM_V13
%n^{-1}\sum_{t=1}^{n}W^{ijl}_{2t} \stackrel{\text{a.s.}}{\rightarrow}0 ,\quad %ALM_V13
%n^{-1}\sum_{t=1}^{n}W^{ijl}_{3t} \stackrel{\text{a.s.}}{\rightarrow}0 ,\\ %ALM_V13
n^{-1}\sum_{t=1}^{n}W^{(1)ijl}_{t} \stackrel{\text{a.s.}}{\rightarrow}0 ,\quad %ALM_V13
n^{-1}\sum_{t=1}^{n}W^{(2)ijl}_{t} \stackrel{\text{a.s.}}{\rightarrow}0 ,\quad %ALM_V13
n^{-1}\sum_{t=1}^{n}W^{(3)ijl}_{t} \stackrel{\text{a.s.}}{\rightarrow}0 ,\\ %ALM_V13
n^{-1}\sum_{t=1}^{n}X^{ilj}_{t} \stackrel{\text{a.s.}}{\rightarrow}0, \quad \text{and} \quad n^{-1}\sum_{t=1}^{n}Y^{ijl}_{t} \stackrel{\text{a.s.}}{\rightarrow}0,\end{eqnarray*}
and this for $i, j, l = 1,...,m$.
%which allow the determination of $\tilde{\Psi}_{1}$ and $\tilde{\Psi}_{2}$. 
%Furthermore TA Lemma 4.17 implies that an upper bound for the third term is related to the upper bound on the second term of \eqref{dr3bis} and that the two assumptions i and ii of TA Lemma 4.8 are fulfilled. Hence 
%\begin{eqnarray}
%\frac{1}{n}\sum_{t=1}^{n}\left|\tr\left(\epsilon^{T}_{t}\epsilon_{t}\right)\right|<\infty  \quad \text{a.s.}, \label{H6.10}
%\end{eqnarray}
%and 
%\begin{eqnarray}
%\frac{1}{n}\sum_{t=1}^{n}\left|\tr\left(\epsilon^{T}_{t-k}\epsilon_{t}\right)\right|<\infty  \quad \text{a.s.}, \label{H6.10bis}
%\end{eqnarray}
%for all finite $k$.

%Consequently, the proof is completed so that $H_{2.4}$ is checked. Then for every $\eps>0$, there exists an event $E$ with $P_{\theta^{0}}(E)> 1 -\eps $ and an $n_0$ such that, for $n >n_0$ on $E$, $Q_{n}(\theta )$ reaches a relative minimum at the point $\widehat{\theta}_{n}$. Consequently, there exists an estimator $\hat{\theta}_{n}$ such that $\hat{\theta}_{n}\stackrel{\text{a.s.}}{\rightarrow} \theta^{0}$. %ALM_V13
Consequently, the proof is completed so that $\Hb_{2.4}$ is checked. Then for every $\eps>0$, there exists an event $E$ with $P_{\theta^{0}}(E)> 1 -\eps $ and an $n_0$ such that, for $n >n_0$ on $E$, $Q_{n}(\theta )$ reaches a relative minimum at the point $\widehat{\theta}_{n}$. Consequently, there exists an estimator $\hat{\theta}_{n}$ such that $\hat{\theta}_{n}\stackrel{\text{a.s.}}{\rightarrow} \theta^{0}$ as $n\rightarrow\infty$. %ALM_V13

%%%%%%%%%%%%%%%%%%%%%%%%%%%%%%%%%%%%%%%%%%%%%%%%%%%%%%%%%%%%%%%%%%%%%%%
%\subsubsection{Proof of $H_{2.5}$}\label{PH2_5} %ALM_V13
\subsubsection{Proof of $\Hb_{2.5}$}\label{PH2_5} %ALM_V13
%From TA Lemma 4.10 we can determine the explicit form of the left hand side of $H_{2.5}$ for all $1 \leq i, j \leq m$: %ALM_V13
From TA Lemma 4.10 we can determine the explicit form of the left-hand side of $\Hb_{2.5}$ for all $1 \leq i, j \leq m$: %ALM_V13
\begin{eqnarray}
%\hspace{-4.5cm}
&&\frac{4}{n}{\sum_{t=1}^{n}}\left[ \left\{ \frac{\partial e_{t}^{T}(\theta) }{\partial \theta_{j}} 
\Sigma^{-1}_{t}(\theta) 
\frac{\partial e_{t}(\theta) }{\partial \theta_{i}}\right\}_{\theta = \theta^{0}} - E_{\theta^{0}}\left(\frac{\partial e_{t}^{T}(\theta) }{\partial \theta_{j}}\Sigma^{-1}_{t}(\theta) \frac{\partial e_{t}^{}(\theta)}{\partial \theta_{i}}\right)\right]\nonumber\\
&& + \frac{2}{n}{\sum_{t=1}^{n}}\left[ \left\{ \frac{\partial e_{t}^{T}(\theta)}{\partial \theta_{i}}\right\}_{\theta = \theta^{0}} - E_{\theta^{0}}\left(\frac{\partial e_{t}^{T}(\theta) }{\partial \theta_{i}}\right)\right]K^{}_{t,j}  \nonumber\\
&& +  \frac{2}{n}{\sum_{t=1}^{n}}\left[ \left\{ \frac{\partial e_{t}^{T}(\theta)}{\partial \theta_{j}} \right\}_{\theta = \theta^{0}} - E_{\theta^{0}}\left(\frac{\partial e_{t}^{T}(\theta) }{\partial \theta_{j}}\right)\right]K^{}_{t,i},  \label{H5.01}    
\end{eqnarray}
where $K_{t,i}$ is defined in TA Lemma 4.10. 
%While checking $H_{2.3}$, we have shown that   %ALM_V13
While checking $\Hb_{2.3}$, we have shown that   %ALM_V13
$$ \frac{4}{n}{\sum_{t=1}^{n}}\left[ 
\left\{ \frac{\partial e_{t}^{T}(\theta) }{\partial \theta_{j}} 
\Sigma^{-1}_{t}(\theta) 
\frac{\partial e_{t}^{}(\theta)}{\partial \theta_{i}}\right\}_{\theta =\theta^{0}} 
- E_{\theta^{0}}\left(\frac{\partial e_{t}^{T}(\theta) }{\partial \theta_{j}}\Sigma^{-1}_{t}(\theta) \frac{\partial e_{t}^{}(\theta)}{\partial \theta_{i}}\right)\right] \stackrel{\text{\mbox{a.s.}}}{\rightarrow} 0.$$
There remains to prove that the second and third terms of (\ref{H5.01}) also tend \mbox{a.s.} to zero. To achieve that, let us consider  
$$
\tilde{Z}^{}_{t,ij}(\theta)=\left( \left\{ \frac{\partial e_{t}^{T}(\theta) }{\partial \theta_{j}}\right\}_{\theta =\theta^{0}} -E_{\theta}\left(\frac{\partial e_{t}^{T}(\theta) }{\partial \theta_{j}}\right)\right)K^{}_{t,i}, 
$$
%for $i,j=1,...,m$. Then, by TA Lemma 4.22, the two assumptions i and ii of Lemma A.1 are verified so the last two terms of (\ref{H5.01}) also tend to zero almost surely.   %ALM_V13
for $i,j=1,...,m$. Then, by TA Lemma 4.22 (involving assumptions $\Hb_{3.3}$, $\Hb_{3.4}$, $\Hb_{3.5}$ and $\Hb_{3.7}$), the two assumptions i and ii of Lemma A.1 are verified, entailing that the last two terms of (\ref{H5.01}) also tend to zero almost surely.   %ALM_V13

As a conclusion, the asymptotic convergence of the estimator $\hat{\theta}_{n}$ towards the normal distribution is ensured and the proof of Theorem 3.1 is achieved.
$\hspace{1cm} \hfill \square $

%%%%%%%%%%%%%%%%%%%%%%%%%%%%%%%%%%%%%%%%%%%%%%%%%%%%%%%%%%%%%%%%%%%%%%%
%%%%%%%%%%%%%%%%%%%%%%%%%%%%%%%%%%%%%%%%%%%%%%%%%%%%%%%%%%%%%%%%%%%%%%%
%%%%%%%%%%%%%%%%%%%%%%%%%%%%%%%%%%%%%%%%%%%%%%%%%%%%%%%%%%%%%%%%%%%%%%%
\section{Assumptions checking}\label{ASS1}
%In this appendix we consider a tdVAR(1) process (\ref{SC.1}).
We check the assumptions of Theorem \ref{T1.3.1} for the two examples of Section~\ref{S.5}, hereby providing a theoretical foundation for the simulation results of Section~\ref{S.5.4}. Since several of the assumptions are somewhat similar, we have only covered once each argument. Also, since Example 2 is a generalization of Example 1, we have avoided to repeat some of the verifications when they are too similar.  
%%%%%%%%%%%%%%%%%%%%%%%%%%%%%%%%%%%%%%%%%%%%%%%%%%%%%%%%%%%%%%%%%%%%%%%
%%%%%%%%%%%%%%%%%%%%%%%%%%%%%%%%%%%%%%%%%%%%%%%%%%%%%%%%%%%%%%%%%%%%%%%
\subsection{Example \ref{S.5.1} }\label{E.Ch.1}
%%%%%%%%%%%%%%%%%%%%%%%%%%%%%%%%%%%%%%%%%%%%%%%%%%%%%%%%%%%%%%%%%%%%%%%
\subsubsection{Assumption $\Hb_{3.1}$}
Trivial
%%%%%%%%%%%%%%%%%%%%%%%%%%%%%%%%%%%%%%%%%%%%%%%%%%%%%%%%%%%%%%%%%%%%%%%
\subsubsection{Assumption $\Hb_{3.2}$} \label{h32}
In order check this hypothesis, we shall have recourse to the results of Section \ref{S1.3.1.2}. In this  example the coefficients of the pure moving average representation of (\ref{SC.6}) are given by
$$
\psi_{tk}(\theta)= \prod_{l=0}^{k-1}A_{t-l}(\theta), 
$$
for $k = 1, 2,..., t - 1$. The coefficients of the pure autoregressive form are 
$$
\pi_{t1}(\theta) = A_{t}(\theta), 
$$
and $\pi_{tk}(\theta)=0 $  if $k = 2, ..., t-1.$ 

Then, by using (\ref{p.1})-(\ref{p.3}), we can calculate $\psi_{tik}(\theta,\theta^{0})$, $\psi_{tijk}(\theta,\theta^{0})$ and \\ $\psi_{tijlk}(\theta,\theta^{0})$. For example, denoting $A_{t}^{(k-1)} = \prod_{l = 1}^{k-1} A_{t-l}(\theta^{0})$, and its $(i,j)$ element $A_{t,i,j}^{(k-1)}$, $i, j= 1, 2$, we have  
\begin{eqnarray}
\psi_{tik}(\theta,\theta^{0})= \frac{\partial \pi_{t1}(\theta)}{\partial\theta_{i}} \prod_{l=1}^{k-1} A_{t-l}(\theta^{0})&=& \frac{\partial A_{t}(\theta)}{\partial\theta_{i}}  A_{t}^{(k-1)}, \label{b1.1}
\end{eqnarray}                          
\begin{eqnarray}
\psi_{tijk}(\theta,\theta^{0})= \frac{\partial^{2} \pi_{t1}(\theta)}{\partial\theta_{i}\partial\theta_{j}} \prod_{l=1}^{k-1} A_{t-l}(\theta^{0})&=& \frac{\partial^{2} A_{t}(\theta)}{\partial\theta_{i}\partial\theta_{j}} A_{t}^{(k-1)}, \label{b1.2}
\end{eqnarray}
\begin{eqnarray}
\psi_{tijlk}(\theta,\theta^{0})= \frac{\partial^{3} \pi_{t1}(\theta)}{\partial\theta_{i}\partial\theta_{j}\partial\theta_{l}} \prod_{l=1}^{k-1} A_{t-l}(\theta^{0}) &=& \frac{\partial^{3} A_{t}(\theta)}{\partial\theta_{i}\partial\theta_{j}\partial\theta_{l}} A_{t}^{(k-1)}, \label{b1.3}
\end{eqnarray}
%Then we have %By specializing the results of Section \ref{S1.3.1.2}, we have 
%%\begin{eqnarray*}
%%\psi^{(n)}_{tk}(\theta)&=& \prod_{l=1}^{k-1}A_{t-l}(\theta), \label{b1.1}
%%\end{eqnarray*}
%%and 
%\begin{eqnarray*}
%\psi_{tik}(\theta,\theta^{0})&=& \frac{\partial A_{t}(\theta)}{\partial\theta_{i}} \prod_{l=1}^{k-1}A_{t-l}(\theta^{0}). \label{b1.3}
%\end{eqnarray*}  
for $k = 1, 2,..., t - 1$, where $A_{t}^{(0)} = I_2$. 
%The coefficients of the pure autoregressive form are as follows
%\begin{eqnarray}
%\pi^{(n)}_{tk}(\theta) &=& A^{(n)}_{t}(\theta), \label{b1.2}
%\end{eqnarray}
%for $k = 1$ and equal to $0$ if $k = 2, ..., t-1.$ 
%GMIn addition from (\ref{SC.6}) and \ref{SC.7})  
%GM\begin{eqnarray}
%GM\prod_{l=1}^{k-1}A_{t-l}(\theta^{0})&=& \prod_{l=1}^{k-1} \left(\begin{array}{cc}	A_{t-l,11}(\theta^{0}) & A_{t-l,12}(\theta^{0}) \\ A_{t-l,21}(\theta^{0}) & A_{t-l,22}(\theta^{0}) \end{array}\right),\label{b1.4}
%GM\end{eqnarray} 
%GMThus the product defined in (\ref{b1.4}) can be written like  %of those matrices 
%GM\begin{eqnarray}
%GM\left(\prod_{l=1}^{k-1}A_{t-l}(\theta^{0})\right)_{j_{1}j_{2}}&=& \sum_{i_{1},i_{2},...,i_{k-2} \in \{1,2\}} A_{t-1,i_{0}i_{1}}(\theta^{0})A_{t-2,i_{1}i_{2}}(\theta^{0})...A_{t-k+1,i_{k-2}i_{0}}(\theta^{0})\quad \quad \label{b1.5}
%GM\end{eqnarray}
%GMwith $i_{0} =j_{1}$, $i_{k-1}=j_{2}$ and $j_{1},j_{2}\in \{1,2\}$, as an example for $i_{0} =1$, $i_{k-1}=1$, (\ref{b1.5}) becomes 
%GM\begin{eqnarray}
%GM\left(\prod_{l=1}^{k-1}A_{t-l}(\theta^{0})\right)_{j_{1}j_{2}}&=& \sum_{i_{1},i_{2},...,i_{k-2} \in \{1,2\}} \prod_{l=1}^{k-1} A_{t-l,i_{l-1}i_{l}}(\theta^{0}) \label{b1.51}
%GM\end{eqnarray}
Obviously, checking the assumptions in the general setup happens to be a complicated and tedious task, hence, for the sake of simplification, we shall consider the model defined in (\ref{btheor}). It can be shown by induction that 
\begin{equation}
A_{t}^{(k-1)} = \left(\begin{array}{cc} A_{t,1,1}^{(k-1)} & A_{t,1,2}^{(k-1)}\\ 0 & A_{t,2,2}^{(k-1)} \end{array}\right),	
\label{At}
\end{equation}
for $k \geq 2$, where
\begin{equation}
A_{t,1,1}^{(k-1)} = (A^{'0}_{11})^{k-1}\prod_{l=1}^{k-1}\sin(a(t-l)), \quad A_{t,2,2}^{(k-1)} = (A^{'0}_{22})^{k-1}\prod_{l=1}^{k-1}\sin(b(t-l)), 
\label{prodA11_22}
\end{equation}
\begin{equation}
A_{t,1,2}^{(k-1)} = \frac{1}{2} \sum_{l = 1}^{k-1} (A^{'0}_{11})^{k-l-1}(A^{'0}_{22})^{l-1} \prod_{f=1}^{k-2} \sin(c_{lf} (t - f - \delta_{lf})) ,
\label{prodA12}
\end{equation}
and 
$c_{lf}=a$ and $\delta_{lf}=0$, for $l + f \leq k - 1$, and $c_{lf}=b$ and $\delta_{lf}=1$, for $l + f > k - 1$. 
For example, for $k=4$
\begin{eqnarray*}
A_{t,1,2}^{(3)} = \frac{1}{2} \left( (A^{'0}_{11})^{2} \sin(a(t-1))\sin(a(t-2)) + A^{'0}_{11}A^{'0}_{22}\sin(a(t-1))\sin(b(t-3)) \right.\\
+ \left. (A^{'0}_{22})^{2} \sin(b(t-2))\sin(b(t-3)) \right).
\end{eqnarray*}

For $i = 1$, \mbox{i.e.} $\theta_{1} = A'_{11}$, and using (\ref{At}) we have 
\begin{eqnarray}
%\psi_{t1k}(\theta^{0})&=& \left(\begin{array}{cc}	\sin(2t) &  0 \\ 
%\psi_{t1k}(\theta,\theta^{0})&=& \left(\begin{array}{cc}	\sin(at) &  0 \\ 0 & 0 \end{array}\right) \left(\begin{array}{cc} A_{t,1,1}^{(k-1)} & A_{t,1,2}^{(k-1)}\\ 0 & - \end{array}\right) \nonumber\\ %ALM_V17
\psi_{t1k} &=& \left(\begin{array}{cc}	\sin(at) &  0 \\ 0 & 0 \end{array}\right) \left(\begin{array}{cc} A_{t,1,1}^{(k-1)} & A_{t,1,2}^{(k-1)}\\ 0 & - \end{array}\right) \nonumber\\ %ALM_V17
%                             &=& 	\left(\begin{array}{cc}	\sin(2t)
%\left(\prod_{l=1}^{k-1}A_{t-l}(\theta^{0})\right)_{11}&  \sin(2t)
                             &=& \sin(at)	\left(\begin{array}{cc}	A_{t,1,1}^{(k-1)} &  A_{t,1,2}^{(k-1)} \\ 0 & 0 \end{array}\right) \label{b9} 
                             %&=& 	\left(\begin{array}{cc} \left(A^{'0}_{11}\right)^{k-1}\prod_{l=1}^{k-1} \sin(2(t-l))&  \sin(2t) \\ 0 & 0 \end{array}\right),                      
\end{eqnarray} 
where a dash will always indicate an element that will not be used. By the same way for $\theta_{2} = A'_{22}$
\begin{eqnarray}
%\psi_{t2k}(\theta,\theta^{0})&=& \left(\begin{array}{cc}	0 &  0 \\ 0 & %\sin(t) \end{array}\right) \left(\begin{array}{cc}-& - \\ 0 & \left %ALM_V17
\psi_{t2k} &=& \left(\begin{array}{cc}	0 &  0 \\ 0 & %\sin(t) \end{array}\right) \left(\begin{array}{cc}-& - \\ 0 & \left %ALM_V17
\sin(bt) \end{array}\right) \left(\begin{array}{cc}-& - \\ 0 & A_{t,2,2}^{(k-1)} \end{array}\right)\nonumber\\
                            &=& \sin(bt) \left(\begin{array}{cc} 0 &  0 \\% 0 & \sin(t)\left(\prod_{l=1}^{k-1}A_{t-l}(\theta^{0})\right)_{22}
0 &  A_{t,2,2}^{(k-1)} \end{array}\right). \label{b1.8}                     
\end{eqnarray} 
We define the constant
\begin{eqnarray}
%\Phi = \max\left\{ \sup_{t=1,...,n}\left[\left|A^{'0}_{11}\sin(2t)\right|\right], \sup_{t=1,...,n}\left[A^{'0}_{22}\left|\sin(t)\right|\right]\right\},\label{b1.62}
\Phi^{1/2} = \max\left\{ \left|A^{'0}_{11}\right|, \left|A^{'0}_{22}\right|\right\},\label{b1.62}
\end{eqnarray} 
which is such that  $0<\Phi< 1$. 
Now by using (\ref{prodA11_22})-(\ref{b1.62}) and bounding sines by 1, we can show that 
\begin{eqnarray}
A_{t,1,1}^{(k-1)} \leq   \Phi^{\frac{k-1}{2}},\label{b1.631}
\end{eqnarray}
\begin{eqnarray}
%A_{t,1,2}^{(k-1)} \leq \red (k-1) \Phi^{\frac{k}{2}-1}\black\label{b1.632} %ALM_V13
A_{t,1,2}^{(k-1)} \leq (k-1) \Phi^{\frac{k}{2}-1} \label{b1.632} %ALM_V13
\end{eqnarray}
and 
\begin{eqnarray}
A_{t,2,2}^{(k-1)} \leq  \Phi^{\frac{k-1}{2}}.\label{b1.633}
\end{eqnarray}
Consequently, from (\ref{b1.8}) and (\ref{b1.633}),
%\begin{eqnarray*}
%\sum_{k=\nu}^{t-1}\left\|\psi_{t1k}(\theta^{0})\right\|^{2}_{F} &= & \sum_{k=\nu}^{t-1} \left[\left(\sin(2t)\left(\prod_{l=1}^{k-1}A_{t-l}(\theta^{0})\right)_{11}\right)^{2} + \left(\sin(2t)\left(\prod_{l=1}^{k-1}A_{t-l}(\theta^{0})\right)_{12}\right)^{2} \right] \\
%&\leq & \sum_{k=\nu}^{t-1}\left[\Phi^{k-1} + (k-1) \Phi^{k-1} \right] \\
%&\leq & \sum_{k=\nu}^{t-1}\Phi^{k-1} + \sum_{k=\nu}^{t-1}(k-1) \Phi^{k-1}  \\
%&\leq & ... \\
%&\leq & N_{1}\Phi^{\nu-1},
%\end{eqnarray*}
$$
%\sum_{k=\nu}^{t-1}\left\|\psi_{t2k}(\theta^{0})\right\|^{2}_{F} =  \sum_{k=\nu}^{t-1} \left[\sin(bt) A_{t,2,2}^{(k-1)} \right]^{2}  %ALM_V17
\sum_{k=\nu}^{t-1}\left\|\psi_{t2k} \right\|^{2}_{F} =  \sum_{k=\nu}^{t-1} \left[\sin(bt) A_{t,2,2}^{(k-1)} \right]^{2}  %ALM_V17
\leq \sum_{k=\nu}^{t-1} \Phi^{k-1} \leq N_{1}\Phi^{\nu-1},
$$
where $N_{1} = 1/(1 - \Phi)$. By the same way, from (\ref{b9}) and (\ref{b1.631})-(\ref{b1.632}),  
%\begin{eqnarray*}
%\sum_{k=\nu}^{t-1}\left\|\psi_{t2k}(\theta^{0})\right\|^{2}_{F} &= & \sum_{k=\nu}^{t-1} \left[\sin(t)\left(\prod_{l=1}^{k-1}A_{t-l}(\theta^{0})\right)_{22}\right]^{2} \\
%&\leq &\sum_{k=\nu}^{t-1} \Phi^{k-1}\\
%&\leq &...\\
%&\leq & N_{1}\Phi^{\nu-1}, 
%\end{eqnarray*}
\begin{eqnarray}
%\sum_{k=\nu}^{t-1}\left\|\psi_{t1k}(\theta^{0})\right\|^{2}_{F} &= & \sum_{k=\nu}^{t-1} \sin^{2}(at) \left[ \left( A_{t,1,1}^{(k-1)} \right)^{2} + \left( A_{t,1,2}^{(k-1)} \right)^{2} \right] \label{A14ter} \\  %ALM_V17
\sum_{k=\nu}^{t-1}\left\|\psi_{t1k} \right\|^{2}_{F} &= & \sum_{k=\nu}^{t-1} \sin^{2}(at) \left[ \left( A_{t,1,1}^{(k-1)} \right)^{2} + \left( A_{t,1,2}^{(k-1)} \right)^{2} \right] \label{A14ter} \\  %ALM_V17
%&\leq & \sum_{k=\nu}^{t-1}\Phi^{k-1} + \sum_{k=\nu}^{t-1}(k-1)^{2} \red\Phi^{k-1} \black , \nonumber
&\leq & \sum_{k=\nu}^{t-1}\Phi^{k-1} + \sum_{k=\nu}^{t-1}(k-1)^{2} \Phi^{k-2}, \nonumber
\end{eqnarray}
but this cannot be bounded by $N_{1}\Phi^{\nu-1}$ for some constant $N_{1}$, independently of $\nu$, so a more subtle upper bound needs to be found. Using \eqref{prodA12}, an upper bound of the element $(1,2)$ is equal to
\begin{equation}
\frac{1}{2} \Phi^{k/2-1} \sum_{\ell = 1}^{k-1} \prod_{f=1}^{k-2} \sin(c_{lf} (t - f - \delta_{lf})) .
\label{UBprodA12}
\end{equation}
It should be possible to adapt some results on products of sines (e.g., Freiman \& Halberstam, 1988, and Janous \& King, 2000) to the general case. 
More precisely, let $Q_{k} = \max_{P \geq 2} \prod_{t=1}^k \sin(2\pi t/P)$, then $\lim_{k \rightarrow \infty} (Q_k)^{1/k} = 0.6098579...$. 
But it will be tricky given the specific form \eqref{UBprodA12}, and we will restrain the proof by assuming a sufficient (but not necessary) condition that $a = 2\pi /P_{1}$ and $b = 2\pi /P_{2}$, for some strictly positive integers $P_{1}$ and $P_{2}$. 
%Then $\sin(a (t - f)) = 0$, for $t - f = g_{1} P_{1}$, for $g_{1} = 1, 2, \dots$, and $\sin(b (t - f - 1)) = 0$, for $t - f - 1 = g_{2} P_{2}$, for $g_{2} = 1, 2, \dots$. 
Then $\sin(a (t - f)) = 0$, for $t - f = g_{1} P_{1}$, for $g_{1} \in \Z$, and $\sin(b (t - f - 1)) = 0$, for $t - f - 1 = g_{2} P_{2}$, for $g_{2} \in \Z$. 
%Let us consider $k \geq \tilde{k} =^{\text{def}} P_{1} + P_{2} +2$ (TO BE CHECKED, PERHAPS A PRODUCT INSTEAD OF A SUM). Then, for any given $t$, it is 
Let $\tilde{k} =^{\text{def}} P_{1} + P_{2} +1$. Then, for any given $t$, 
take the remainders $g_{1}$ and $g_{2}$ of the Euclidean divisions of $t$ by 
$P_{1}$ and  of $t - 1$ by $P_{2}$, respectively. 
These remainders are between $0$ and, respectively, $P_{1} - 1$ and $P_{2} - 1$, hence smaller than $\tilde{k} - 2$. Hence 
%it is possible to find two integers $g_{1}$ and $g_{2}$ such that $t - g_{1}P_{1}$ or $t - g_{2}P_{2} -1$ is an integer between 1 and $k-2$ and consequently 
for any $k \geq \tilde{k}$, there is some integer $f$, either $t - g_{1} P_{1}$, or $t - 1 - g_{2} P_{2}$ between $0$ and $k - 2$ such that 
either $\sin(a(t-f))=0$ or $\sin(b(t-f-1))=0$ %for some $f = 1, \dots, k - 2$
 and finally each term of \eqref{UBprodA12} vanishes. Hence, taking $\nu > \tilde{k}$, the second term of \eqref{A14ter} vanishes for all $t$ and the Frobenius norm can be bounded by some expression $N_{1} \Phi^{\nu -1}$ for some $N_{1}$. Note also that, under the conditions of integral periods, we have 
$$
%\left\|\psi_{t1k}(\theta^{0})\right\|^{2}_{F} = \left\|\psi_{t2k}(\theta^{0})\right\|^{2}_{F} = 0, \quad k > \tilde{k}. %ALM_V17 
\left\|\psi_{t1k} \right\|^{2}_{F} = \left\|\psi_{t2k} \right\|^{2}_{F} = 0, \quad k > \tilde{k}. %ALM_V17 
$$

By using (\ref{b1.1})-(\ref{b1.2})-(\ref{b1.3}) the other inequalities of this assumption are checked.  
%\frac{1}{(1-\Phi)}
%%%%%%%%%%%%%%%%%%%%%%%%%%%%%%%%%%%%%%%%%%%%%%%%%%%%%%%%%%%%%%%%%%%%%%%
\subsubsection{Assumption $\Hb_{3.3}$} 
Trivial.
%%%%%%%%%%%%%%%%%%%%%%%%%%%%%%%%%%%%%%%%%%%%%%%%%%%%%%%%%%%%%%%%%%%%%%%
\subsubsection{Assumption $\Hb_{3.4}$} 
Trivial.
%%%%%%%%%%%%%%%%%%%%%%%%%%%%%%%%%%%%%%%%%%%%%%%%%%%%%%%%%%%%%%%%%%%%%%%
\subsubsection{Assumption $\Hb_{3.5}$} 
Trivial.
%%%%%%%%%%%%%%%%%%%%%%%%%%%%%%%%%%%%%%%%%%%%%%%%%%%%%%%%%%%%%%%%%%%%%%%
\subsubsection{Assumption $\Hb_{3.6}$} \label{h36}
The second term in $\Hb_{3.6}$ is equal to 0, so it remains to show that 
\begin{eqnarray*}
\lim_{n\rightarrow\infty}\frac{1}{n}\sum_{t=1}^{n} \left( E_{\theta^{0}}\left(\frac{\partial e^{T}_{t}(\theta)}{\partial\theta_{i}}
\Sigma^{-1}_{t}(\theta)\frac{\partial e_{t}(\theta)}{\partial \theta_{j}}\right)\right)% =  V_{ij}(\theta^{0}) %ALM_V13
 =  V_{ij} %ALM_V13
\end{eqnarray*} 
%\mbox{a.s.} for $i, j =1,2$, where the matrix $V(\theta^{0})=(V_{1 \leq i,j \leq 2}(\theta^{0}))_{i,j}$ is a strictly positive definite matrix. From (\ref{32.13}) we have =  V_{ij}(\theta^{0}) %ALM_V13
for $i, j =1,2$, where the matrix $V=(V_{i,j})_{1 \leq i,j \leq 2}$ is a strictly positive definite matrix. From (\ref{32.13}) we have % %ALM_V13
\begin{eqnarray*}
%\frac{\partial e^{T}_{t}(\theta)}{\partial\theta_{1}} &=& \sum_{k=1}^{t-1}\psi_{t1k}(\theta,\theta^{0})\epsilon_{t-k}, %ALM_V13
\frac{\partial e_{t}(\theta)}{\partial\theta_{1}} &=& \sum_{k=1}^{t-1}\psi_{t1k}(\theta,\theta^{0})\epsilon_{t-k}, %ALM_V13
\end{eqnarray*} 
entailing, since we have taken $\Sigma = I_{2}$,  
\begin{eqnarray}
E_{\theta^{0}}\left(\frac{\partial e^{T}_{t}(\theta)}{\partial\theta_{1}}
\Sigma^{-1}\frac{\partial e_{t}(\theta)}{\partial \theta_{1}}\right)&=& E_{\theta^{0}}\left[\left(\sum_{k=1}^{t-1}\psi_{t1k}(\theta,\theta^{0})\epsilon_{t-k}\right)^{T}\Sigma^{-1}\sum_{k=1}^{t-1}\psi_{t1k}(\theta,\theta^{0})\epsilon_{t-k}\right]\nonumber\\
&=& \sum_{k=1}^{t-1} \sin^{2}(at) \tr\left[ \left(\begin{array}{cc}  A_{t,1,1}^{(k-1)}  &   A_{t,1,2}^{(k-1)}  \\ 0 & 0 \end{array}\right)^{T}
\left(\begin{array}{cc}  A_{t,1,1}^{(k-1)} &  A_{t,1,2}^{(k-1)}  \\ 0 & 0 \end{array}\right)\right]\nonumber\\
&=& \sin^{2}(at) \sum_{k=1}^{t-1} \left[ \left( A_{t,1,1}^{(k-1)} \right)^{2} + \left( A_{t,1,2}^{(k-1)} \right)^{2} \right].\label{B19}
\end{eqnarray}
The terms for $k > \tilde{k}$ vanish, as explained above. The remaining terms are bounded by $\Phi^{k-1}$, also as above, and the sum is strictly positive, %ALM_V13
 at least for some $t$'s. Therefore, taking the average for $t=1$ to $n$ gives a finite %strictly positive limit $V_{11}(\theta^{0})$ when $n \rightarrow \infty$.  %ALM_V13
strictly positive limit $V_{11}$ when $n \rightarrow \infty$.  %ALM_V13

%  A_{t,2,2}^{(k-1)} 
%Therefore
%\begin{eqnarray}
%V_{11}(\theta^{0}) &=& \lim_{n\rightarrow\infty}\frac{1}{n}\sum_{t=1}^{n}\sum_{k=1}^{t-1}\left(A^{'0}_{11}\right)^{2k-2}\left[\prod_{l=0}^{k-1} \sin(b(t-l))\right]^{2} \nonumber\\
%&\leq&  \lim_{n\rightarrow\infty}\frac{1}{n}\sum_{t=1}^{n}\sum_{k=1}^{t-1}\Phi^{k-1} .\label{b1.81}  
%\end{eqnarray}
%\red je ne vois pas l'utilie de ce $>0$; il est clair que notre premiere expression est positive, donc la borne le sera d'office aussi 
\black By the same method for $i, j = 2$
\begin{eqnarray*}
E_{\theta^{0}}\left(\frac{\partial e^{T}_{t}(\theta)}{\partial\theta_{2}}
\Sigma^{-1}\frac{\partial e_{t}(\theta)}{\partial \theta_{2}}\right)&=& E_{\theta^{0}}\left[\left(\sum_{k=1}^{t-1}\psi_{t2k}(\theta,\theta^{0})\epsilon_{t-k}\right)^{T}\Sigma^{-1}\sum_{k=1}^{t-1}\psi_{t2k}(\theta,\theta^{0})\epsilon_{t-k}\right]\nonumber\\
&=& \sum_{k=1}^{t-1} \sin^{2}(bt) \tr\left[ \left(\begin{array}{cc} 0 & 0\\  0  &   A_{t,2,2}^{(k-1)}   \end{array}\right)^{T}
\left(\begin{array}{cc} 0 & 0  \\ 0 &  A_{t,2,2}^{(k-1)}  \end{array}\right) \right]\nonumber\\
&=& \sin^{2}(bt) \sum_{k=1}^{t-1} \left( A_{t,2,2}^{(k-1)} \right)^{2}  ,
\end{eqnarray*}
%with a similar consequence to provide a finite strictly positive limit $V_{22}(\theta^{0})$ when $n \rightarrow \infty$.  %ALM_V13
with a similar consequence to provide a finite strictly positive limit $V_{22}$ when $n \rightarrow \infty$.  %ALM_V13

%hence
%\begin{eqnarray}
%V_{22}(\theta^{0}) &=& \lim_{n\rightarrow\infty}\frac{1}{n}\sum_{t=1}^{n}\sum_{k=1}^{t-1}\left(\cos(t + A^{''}_{22})\right)^{2}\left[\prod_{l=1}^{k-1} \frac{1}{2}\sin(t-l + A^{''}_{22})\right]^{2} \nonumber\\
%&\leq&  \lim_{n\rightarrow\infty}\frac{1}{n}\sum_{t=1}^{n}\sum_{k=1}^{t-1}\Phi^{k-1}.\label{b1.82}  
%\end{eqnarray}

Furthermore, it can easily be seen 
$$
%V_{12}(\theta^{0}) &=& V_{21}(\theta^{0}) = 0.\label{b1.83}   %ALM_V13
V_{12} = V_{21} = 0.   %ALM_V13
$$
%Consequently, the matrix $V =(V_{1 \leq i,j \leq 2} )_{i,j}$ is strictly positive definite. %ALM_V13 %ALM_V17
Consequently, the matrix $V =(V_{i,j})_{1 \leq i,j \leq 2} $ is strictly positive definite. %ALM_V13 %ALM_V17

\begin{rem}
If we handle the model described in Section \ref{S.5.4}, but without the parameter $A'_{12}$, for the case where the number of observations is $25$, we obtain the following standard errors from that theoretical formula for the two %remaining parameters: $0.2175$ and $0.2303$, respectively, which largely agree 
%remaining parameters: $0.2161$ and $0.2303$, respectively, which largely agree %ALM_V16 %
remaining parameters: $0.2175$ and $0.2303$, respectively, which largely agree %ALM_V16 %
%with the averages drawn from the simulation results: $0.2161$ and $0.2226$, as 
with the averages drawn from the simulation results: $0.2161$ and $0.2226$, as %ALM_V16
 shown in Table \ref{Tab:SC.10}. 
\end{rem}

%%%%%%%%%%%%%%%%%%%%%%%%%%%%%%%%%%%%%%%%%%%%%%%%%%%%%%%%%%%%%%%%%%%%%%%
\subsubsection{Assumption $\Hb_{3.7}$}  \label{h37}
%For the first part of this assumption we have to show that  %ALM_V13
For the first part of this assumption, since $g_{t}$ is bounded, it remains to show that  %ALM_V13
$$
%\frac{1}{n^{2}}\sum_{d= 1}^{n-1}\sum_{t=1}^{n-d} \sum^{t-1}_{k=1}\left\|\psi_{tik}(\theta^{0})\right\|_{F}\left\|\psi_{t+d,i,k+d}(\theta^{0})\right\|_{F}=O\left(\frac{1}{n}\right).  %ALM_V17
\frac{1}{n^{2}}\sum_{d= 1}^{n-1}\sum_{t=1}^{n-d} \sum^{t-1}_{k=1}\left\|\psi_{tik} \right\|_{F}\left\|\psi_{t+d,i,k+d} \right\|_{F}=O\left(\frac{1}{n}\right).  %ALM_V17
$$

For $i = 1$
\begin{eqnarray}
%&&\sum^{t-1}_{k=1}\left\|\psi_{t1k}(\theta^{0})\right\|_{F}\left\|\psi_{t+d,1,k+d}(\theta^{0})\right\|_{F}   %ALM_V17
&&\sum^{t-1}_{k=1}\left\|\psi_{t1k} \right\|_{F}\left\|\psi_{t+d,1,k+d} \right\|_{F}   %ALM_V17
= %&=& 
|\sin(at) \sin(a(t+d))|\nonumber \\%
&& 
\hspace{1.5cm}\times\sum^{t-1}_{k=1} \left[ \left( A_{t,1,1}^{(k-1)} \right)^{2} + \left( A_{t,1,2}^{(k-1)} \right)^{2} \right]^{1/2} \left[ \left( A_{t+d,1,1}^{(k+d-1)} \right)^{2} + \left( A_{t+d,1,2}^{(k+d-1)} \right)^{2} \right]^{1/2} . \nonumber \\
&& \label{b1.10}
\end{eqnarray}
To simplify the proof we assume again that 
$a = 2\pi /P_{1}$ and $b = 2\pi /P_{2}$, for some integers $P_{1}$ and $P_{2}$ and use $\tilde{k}$ as defined above. Then \eqref{b1.10} becomes a sum for $k = 1$ to $\tilde{k} - d$. 

The general term can be bounded by $[\Phi^{k-2}(1+(k-1)^2)\Phi^{k+d-2}(1+(k+d-1)^2)]^{1/2}$ which is of order $d \Phi^{d/2}$. 
Hence
\begin{eqnarray*}
\sum_{d=1}^{n-1}\sum_{t=1}^{n-d}(\ref{b1.10})&=& 
\sum_{t=1}^{n-1}\sum_{d=1}^{n-t}(\ref{b1.10}) \nonumber\\
&\leq& \sum_{t=1}^{n-1}\sum_{d=1}^{n-t} k_{1}(\Phi) \Phi^{d/2} d
\nonumber\\
&\leq&k_{1}(\Phi) \frac{\Phi^{1/2}}{(1-\Phi^{1/2})^{2}} \sum_{t=1}^{n-1} 1\nonumber\\
&\leq& k_{2}(\Phi) n,
\end{eqnarray*}
where $k_{1}(\Phi)$ and $k_{2}(\Phi)$ are constants and where we have used the formula $\sum_{j=1}^{\infty}j x^{j} = x/(1-x)^{2}$, provided $|x|<1$. 
Dividing by $n^2$ thus gives $O(1/n)$, as requested. 

Applying the same method for $i = 2$, with again a sum for $k=1$ to $\tilde{k}$ but of 
$A_{t,2,2}^{(k-1)} A_{t+d,2,2}^{(k+d-1)}$
instead of the product of square roots of sums of squares in \eqref{b1.10}, with a general term bounded this time by $\Phi^{(k-1)/2}\Phi^{(k+d-1)/2}$, we can show that 
\begin{eqnarray*}
%\sum_{d=1}^{n-1}\sum_{t=1}^{n-d}\sum^{t-1}_{k=1}\left\|\psi^{(n)}_{t2k}(\theta^{0})\right\|_{F}\left\|\psi^{(n)}_{t+d,2,k+d}(\theta^{0})\right\|_{F} &\leq & \sum_{d=1}^{n-1} \sum_{t=1}^{n-1} k_{3}(\Phi) \Phi^{d/2} \nonumber \\ %ALM_V17
\sum_{d=1}^{n-1}\sum_{t=1}^{n-d}\sum^{t-1}_{k=1}\left\|\psi^{(n)}_{t2k} \right\|_{F}\left\|\psi^{(n)}_{t+d,2,k+d} \right\|_{F} &\leq & \sum_{d=1}^{n-1} \sum_{t=1}^{n-1} k_{3}(\Phi) \Phi^{d/2} \nonumber \\ %ALM_V17
&\leq & k_{4}(\Phi) n , 
\end{eqnarray*}
with other constants $k_{3}(\Phi)$ and $k_{4}(\Phi)$ and reach the same final conclusion as for $i=1$. 

%For the second part of this assumption, since $\Xi_{t}(\Sigma)$,  $\Sigma \times \Sigma$ and $\kappa(\Sigma)$ are finite constants, similarly to the first part we can show that  %ALM_V13
For the second part of this assumption, since $\Xi_{t}(\Sigma)$,  $\Sigma \otimes \Sigma$ and $\vecc(\Sigma)\vecc(\Sigma)^{T}$ are finite constants, similarly to the first part we can show that the second part of assumption $\Hb_{3.7}$ is fulfilled. %ALM_V13
%\begin{eqnarray*} %ALM_V13
%	&& \frac{1}{n^{2}}\sum_{d=1}^{n-1}\sum_{t=1}^{n-d}\left[\sum^{t-1}_{k=1}\vecc(M^{{ij}^{T}}_{t,k,k}(\theta^{0}))^{T} \Xi_{t}(\Sigma)\vecc(M^{ij}_{t+d,k+d,k+d}(\theta^{0}))\right.\\ %ALM_V13
%&&	+ \underset{}{{\sum^{t-1}\limits_{k_{1} = 1}} {\sum^{t-1}\limits_{k_{2}=1}}}\vecc(M^{{ij}^{T}}_{t,k_{1},k_{2}}(\theta^{0}))^{T}(\Sigma \otimes \Sigma)\vecc(M^{ij}_{t+d,k_{1}+d,k_{2}+d}(\theta^{0})) \\ %ALM_V13
%&&	+ \left. \sum^{t-1}_{k_{1} = 1} \sum^{t-1}_{k_{2}=1}\vecc(M^{{ij}^{T}}_{t,k_{1},k_{2}}(\theta^{0}))^{T} \kappa(\Sigma)\vecc(M^{ij}_{t+d,k_{2}+d,k_{1}+d}(\theta^{0}))\right] =O\left(\frac{1}{n}\right),\quad \text{a.s.} %ALM_V13
%\end{eqnarray*} %ALM_V13
%with 	 %ALM_V13
%$$M^{ij}_{t, k_{1}, k_{2}}(\theta^{0}) = \psi^{T}_{tik_{1}}(\theta^{0})\Sigma^{-1}\psi_{tjk_{2}}(\theta^{0}).$$	 %ALM_V13

%%%%%%%%%%%%%%%%%%%%%%%%%%%%%%%%%%%%%%%%%%%%%%%%%%%%%%%%%%%%%%%%%%%%%%%
%%%%%%%%%%%%%%%%%%%%%%%%%%%%%%%%%%%%%%%%%%%%%%%%%%%%%%%%%%%%%%%%%%%%%%%
\subsection{Example \ref{S.5.2} }\label{E.Ch.2}

The assumptions $\Hb_{3.1}-\Hb_{3.5}$ are easily checked: $\Hb_{3.1}$ is trivial, $\Hb_{3.2}$ remains unchanged but requires the same conditions on integers $P_1,P_2$ as for the previous example, and $\Hb_{3.3}-\Hb_{3.5}$ are readily checked given boundedness of the sine function.  It remains to discuss about $\Hb_{3.6}$ and $\Hb_{3.7}$ which we are going to check now. Note that from (\ref{b1.6}) the vector of the parameters is $\theta = (A'_{11}, A'_{22}, \eta_{11},\eta_{22})^{T}$ and  $\theta^{0} = (A^{'0}_{11}, A^{'0}_{22}, \eta^{0}_{11},\eta^{0}_{22})^{T}$ is its true value. Note that we assume again that $\Phi$ is of the form \eqref{b1.62}.  Finally, we write $s_{11},s_{12}$ and $s_{22}$ the entries of the matrix $\Sigma$ (without loss of generality, we could as well take $\Sigma$ the identity matrix, as was done for the proof of the previous example).

%%%%%%%%%%%%%%%%%%%%%%%%%%%%%%%%%%%%%%%%%%%%%%%%%%%%%%%%%%%%%%%%%%%%%%%
\subsubsection{Assumption $\Hb_{3.6}$} \label{E.Ch.3H3.6}
%We start by making an important observation, namely that the matrix $V(\theta^0)$ is block-diagonal, with blocks for $(i,j)\in\{1,2\}$ and $(i,j)\in\{3,4\}$, respectively. %ALM_V13 
We start by making an important observation, namely that the matrix $V=(V_{i,j})_{1 \leq i,j \leq 4}$ is block-diagonal, with blocks for $(i,j)\in\{1,2\}$ and $(i,j)\in\{3,4\}$, respectively. %ALM_V13 

Regarding the second block, tedious calculations (therefore carried out with \textrm{Mathematica}) for $\theta_3=\eta_{11}$ yield
\begin{align*}
%\tr \left[\Sigma^{-1}_{t}(\theta^{0})\left\{\frac{\partial \Sigma_{t}(\theta)}{\partial \theta_{3}}\right\}_{\theta=\theta^{0}}\Sigma^{-1}_{t}(\theta^{0}) \left\{
%\frac{\partial \Sigma_{t}(\theta)}{\partial\theta_{3}}\right\}_{\theta = \theta^{0}} \right]\\ %ALM_V17
\tr \left[ \left\{ \Sigma^{-1}_{t}(\theta) \frac{\partial \Sigma_{t}(\theta)}{\partial \theta_{3}} \Sigma^{-1}_{t}(\theta) 
\frac{\partial \Sigma_{t}(\theta)}{\partial\theta_{3}} \right\}_{\theta = \theta^{0}} \right]\\ %ALM_V17
&&\hspace{-5cm}=2\sin^2(ct)\frac{\left(\left(e^{\eta_{22}\sin(ct)}s_{11}-s_{12}\right)^2+2(s_{11}s_{22}-s_{12}^2)\right)}{\left(1+e^{(\eta_{11}+\eta_{22})\sin(ct)}\right)^2(s%_{11}s_{22}-s_{12}^2)}=V_{33}(t).\end{align*} %ALM_V13
_{11}s_{22}-s_{12}^2)}=V_{33}(t),\end{align*} %ALM_V13
where $V_{ij}(t)$ is the term $t$ in the sum defining $V_{ij}$, $i,j=1,2,3,4$. %ALM_V13
Since $\Sigma$ is an invertible matrix, $s_{11}s_{22}-s_{12}^2={\rm det}(\Sigma)>0$, hence this expression is clearly positive. Bounding the numerator is straightforward (since $|\sin(ct)|<1$), and $1/\left(1+e^{(\eta_{11}+\eta_{22})\sin(ct)}\right)^2$ can simply be bounded by 1. Hence, there exist a constant $\alpha$ not depending on $t$ %such that $V_{33}(t)\leq\alpha$, consequently the entry $V_{33}(\theta^0)$ is positive and finite. The same conclusion obviously holds for $V_{44}(\theta^0)$. Turning our attention towards $V_{34}(\theta^0)$, we obtain %ALM_V13
such that $V_{33}(t)\leq\alpha$, consequently the entry $V_{33}$ is positive and finite. The same conclusion obviously holds for $V_{44}$. Turning our attention towards $V_{34}$, we obtain %ALM_V13
\begin{align*}
%\tr \left[\Sigma^{-1}_{t}(\theta^{0})\left\{\frac{\partial \Sigma_{t}(\theta)}{\partial \theta_{3}}\right\}_{\theta=\theta^{0}}\Sigma^{-1}_{t}(\theta^{0}) \left\{
%\frac{\partial \Sigma_{t}(\theta)}{\partial\theta_{4}}\right\}_{\theta = \theta^{0}} \right]\\ %ALM_V17
\tr \left[ \left\{ \Sigma^{-1}_{t}(\theta) \frac{\partial \Sigma_{t}(\theta)}{\partial \theta_{3}} \Sigma^{-1}_{t}(\theta) 
\frac{\partial \Sigma_{t}(\theta)}{\partial\theta_{4}} \right\}_{\theta = \theta^{0}} \right]\\ %ALM_V17
&&\hspace{-9cm}=2\sin^2(ct)\frac{\left(s_{12}\left(e^{\eta_{22}\sin(ct)}s_{11}-s_{12}-e^{\eta_{11}\sin(ct)}s_{22}\right)-e^{(\eta_{11}+\eta_{22})\sin(ct)}(s_{11}s_{22}-2s_{12}^2)\right)}{\left(1+e^{(\eta_{11}+\eta_{22})\sin(ct)}\right)^2(s_{11}s_{22}-s_{12}^2)}=V_{34}(t).
\end{align*}
Again, it is a simple exercise to show that this term is bounded independently of $t$. %Now, we know that the diagonal elements of the block $(3,4)$ are positive; showing the positive definiteness of the %What remains to prove is that the block $(3,4)$ is positive definite; we have seen so far that the diagonal elements are $>0$, which is a minimum. 

Let us turn our attention towards the block $(1,2)$ now. 
From (\ref{32.13}) we know that
\begin{eqnarray*}
\frac{\partial e^{T}_{t}(\theta)}{\partial\theta_{i}} &=& \sum_{k=1}^{t-1}\psi_{tik}(%\theta,\theta^{0})g_{t-k}(\theta^{0})\epsilon_{t-k}. %ALM_V13
\theta,\theta^{0})g_{t-k} \epsilon_{t-k}. %ALM_V13
\end{eqnarray*} 
Following Section \ref{h36} for $i,j = 1$ with $\theta_{1} = A'_{11}$, we obtain a generalized version of \eqref{B19} as
\begin{eqnarray*}
E_{\theta^{0}}\left(\frac{\partial e^{T}_{t}(\theta)}{\partial\theta_{1}}
\Sigma^{-1}_{t}\frac{\partial e_{t}(\theta)}{\partial \theta_{1}}\right)&=& E_{\theta^{0}}\left[\left(\sum_{k=1}^{t-1}\psi_{t1k}(\theta,\theta^{0})g_{t-k}\epsilon_{t-k}\right)^{T}\Sigma^{-1}_{t}\sum_{k=1}^{t-1}\psi_{t1k}(\theta,\theta^{0})g_{t-k}\epsilon_{t-k}\right]\nonumber\\
&=& \sin^{2}(at)\sum_{k=1}^{t-1} \tr\left[ \left(\begin{array}{cc}  A_{t,1,1}^{(k-1)}  &   A_{t,1,2}^{(k-1)}  \\ 0 & 0 \end{array}\right)\Sigma_{t-k}
\left(\begin{array}{cc}  A_{t,1,1}^{(k-1)} &  A_{t,1,2}^{(k-1)}  \\ 0 & 0 \end{array}\right)^T\Sigma_{t}^{-1}\right],
\end{eqnarray*}
an expression clearly bounded (under the same conditions as those in B.1). Its positiveness follows from the fact that we consider a Mahalanobis distance in the metric $\Sigma_t^{-1}$. 

Finally, checking that the blocks $(1,2)$ and $(3,4)$ are positive definite has been done numerically in the numerical example of Section 5.2. % (TO DO).

%%%%%%%%%%%%%%%%%%%%%%%%%%%%%%%%%%%%%%%%%%%%%%%%%%%%%%%%%%%%%%%%%%%%%%%
\subsubsection{Assumption $\Hb_{3.7}$} \label{H.3.7ch2}
In the first part of this assumption we have to show that 
$$
%\frac{1}{n^{2}}\sum_{d= 1}^{n-1}\sum_{t=1}^{n-d} \sum^{t-1}_{k=1}\left\|g_{t-k}(\theta^{0})\right\|^{2}_{F}\left\|\psi^{(n)}_{tik}(\theta^{0})\right\|_{F}\left\|\psi^{(n)}_{t+d,i,k+d}(\theta^{0})\right\|_{F}=O\left(\frac{1}{n}\right),\quad \text{a.s.} %ALM_V13
\frac{1}{n^{2}}\sum_{d= 1}^{n-1}\sum_{t=1}^{n-d} \sum^{t-1}_{k=1}\left\|g_{t-k} \right\|^{2}_{F}\left\|\psi^{(n)}_{tik} \right\|_{F}\left\|\psi^{(n)}_{t+d,i,k+d} \right\|_{F}=O\left(\frac{1}{n}\right).%ALM_V13
$$
From (\ref{SC.821}) we have 
$$
%\left\|g_{t-k}(\theta^{0})\right\|^{2}_{F} =2+e^{-2\eta_{11}\sin(ct)}+e^{-2\eta_{22}\sin(ct)} $$ %ALM_V13
\left\|g_{t-k} \right\|^{2}_{F} =2+e^{-2\eta_{11}\sin(ct)}+e^{-2\eta_{22}\sin(ct)} $$ %ALM_V13
which is an easily bounded function. Consequently, we are left with exactly the same expression as for Example 4.1, which solves the question. %Much more tedious are the remaining points, but I cannot imagine there would be a problem.

\end{appendices}


\begin{thebibliography}{9}  

\bibitem{Alj_Azrak_Melard_TA} Alj, A., Azrak, R. \& M\'elard, G. (2015a). Technical appendix to ``Asymptotic properties of QML estimators for VARMA models with time-dependent coefficients, Parts 1 and 2''. Technical report.


\bibitem {AAM} Alj, A., Azrak, R., \& M\'{e}lard, G. (2015b). Asymptotic properties of QML estimators for VARMA models with time-dependent coefficients, part 2. %in 
%In {\it Time series analysis by time dependent models} (eds Azrak, R. \& M\'{e}lard, G.). In preparation. %ALM_V13
In {\it Time series analysis by time dependent models} (eds R. Azrak \& G. M\'{e}lard). In preparation. %ALM_V13

\bibitem {AJM} Alj, A., J\'{o}nasson, K. \& M\'{e}lard, G. (2015c). The exact Gaussian likelihood estimation of time-dependent VARMA models. 
%{\it Comput. Statistics and Data Analysis}, in press. %ALM_V13
{\it Comput. Statist. Data Anal.}, in press. %ALM_V13


\bibitem {AM2006} Azrak, R. \& M\'{e}lard, G. (2006). Asymptotic properties of quasi-%likelihood estimators for ARMA models with time-dependent coefficients, {\it Statist. Inference for Stochastic Processes} \textbf{9}, 279-330. %ALM_V13
likelihood estimators for ARMA models with time-dependent coefficients. {\it Stat. Inference Stoch. Process.} \textbf{9}, 279--330. %ALM_V13

\bibitem {AM2011}Azrak, R. \& M\'elard, G. (2011). Autoregressive models with time-dependent coefficients
- A comparison with Dahlhaus’ approach. ECARES working paper, Universit\'e
Libre de Bruxelles.

\bibitem {Azrak_Melard_book} Azrak, R. \& M\'{e}lard, G. (2015). {\it Time series %analysis by time dependent models}, Chapter 4, book in preparation. %ALM_V13
analysis by time dependent models}. In preparation. %ALM_V13

\bibitem {Basawa_Lund_2001}Basawa, I. V. \& Lund, R. L. S. (2001). Large sample properties of
%parameter estimates for periodic ARMA models, \textit{J. Time Ser.
%Anal.} \textbf{22}, 651-663. %ALM_V13
parameter estimates for periodic ARMA models. \textit{J. Time Series 
Anal.} \textbf{22}, 651--663. %ALM_V13

\bibitem {Basawa_Prakasa_Rao} Basawa, I. V. \& Prakasa Rao, B. L. S. (1980). \textit{ %Statistical Inference for Stochastic Processes}, Academic Press, New York. %ALM_V13
Statistical inference for stochastic processes}. Academic Press, New York. %ALM_V13

\bibitem {Bibi_Francq}Bibi, A. \& Francq, C. (2003). Consistent and asymptotically
normal estimators for cyclically time-dependent linear models.
\textit{Ann. Inst. Statist. Math.} \textbf{55}, 41--68.

\bibitem {Brockwell_Davis_91}Brockwell, P. J. \& Davis, R. A. (1991).
%\textit{Time series: Theory and Methods}, Springer, Wiley. <<<<<< %ALM_V13
\textit{Time series: theory and methods}. Springer, New York. %ALM_V13

\bibitem {Boubacar_Francq} Boubacar Mainassara, Y. \& Francq, C. (2011). Estimating structural VARMA models
%with uncorrelated but non-independent error terms. \textit{J. Multi. Anal.} \textbf{102}, 496-505. %ALM_V13
with uncorrelated but non-independent error terms. \textit{J. Multivariate Anal.} \textbf{102}, 496--505. %ALM_V13

\bibitem {BJR2008} Box, G. E. P., Jenkins, G. M. \& Reinsel G. C. (2008). 
\textit{Time series analysis, forecasting and control}, 4th edn. Wiley, New York. 

%\bibitem {Creal_Koop_Lucas_2013} Creal, D. D., Koopman, S. J., Lucas, A. (2013). %ALM_V17
\bibitem {Creal_Koop_Lucas_2013} Creal, D. D., Koopman, S. J. \& Lucas, A. (2013). %ALM_V17
 %Generalized Autoregressive Score Models with Applications, \textit{J. Appl. Econom.} \textbf{28}, 777-795.  %ALM_V17
Generalized autoregressive score models with Applications, \textit{J. Appl. Econom.} \textbf{28}, 777-795.  %ALM_V17

\bibitem {Dahlhaus96a}Dahlhaus, R. (1996a). Maximum likelihood estimation and model
selection for locally stationary processes. \textit{J.
%Nonparametr. Statist.} \textbf{6}, 171-191. %ALM_V13
Nonparametr. Stat.} \textbf{6}, 171--191. %ALM_V13

\bibitem {Dahlhaus96b}Dahlhaus, R. (1996b). On the Kullback-Leibler information
divergence of locally stationary processes. \textit{Stoch.
Process. Appl.} \textbf{62}, 139--168.

\bibitem {Dahlhaus96c}Dahlhaus, R. (1996c). Asymptotic statistical inference for
nonstationary processes with evolutionary spectra, In
%\textit{Athens Conference on Applied }\textit{Probability on Time
%Series Analysis}, (P. M. Robinson and M. Rosenblatt, eds)
%\textbf{2}. Springer, New York, pp. 145-159. %ALM_V13
\textit{Athens Conference on applied }\textit{probability and time
series analysis 2} (eds P. M. Robinson \& M. Rosenblatt), 145--159. Springer, New York. %ALM_V13

\bibitem {Dahlhaus97}Dahlhaus, R. (1997). Fitting time series models to nonstationary
processes. \textit{Ann. }\textit{Statist.} \textbf{25}, 1--37.

\bibitem {Dahlhaus2000} Dahlhaus, R. (2000). A likelihood approximation for locally stationary processes. {\it Ann. Statist.} \textbf{28}, 1762--1794.

\bibitem{Frei_Halb} Freiman, G. \& Halberstam, H. (1988). On product of sines. \textit{Acta Mathematica} \textbf{49}, 377--385. 

\bibitem{Francq_Raıssi}Francq, C. \& Ra\"issi, H. (2007). Multivariate portmanteau test for autoregressive models
with uncorrelated but nonindependent errors. 
%\textit{J. Time Ser. Anal.} \textbf{28}, 454-470. %ALM_V13
\textit{J. Time Series Anal.} \textbf{28}, 454--470. %ALM_V13

\bibitem {Grillenzoni}Grillenzoni, C. (1990). Modeling time-varying dynamical systems. 
\textit{J. Amer. }\textit{Statist. Assoc.} \textbf{85}, 499--507.

\bibitem {Golub_van_Loan}
Golub, G. \& Van Loan, C. (1996). 
%\textit{Matrix computations}, 3rd edition. Johns Hopkins %ALM_V13
\textit{Matrix computations}, 3rd edn. Johns Hopkins %ALM_V13
University Press, Baltimore. 

\bibitem{Hall_Heyde} Hall, P. \& Heyde, C. C. (1980). \textit{Martingale limit theory and its application}. Academic Press, New York.

\bibitem {Hallin86} Hallin, M. (1986). Non-stationary \textit{q}-dependent processes and time-varying moving average models: invertibility properties and the forecasting problem. \textit{Adv. }\textit{Appl. Probab.} \textbf{18}, 170--210.

\bibitem {Hallin_Ingenbleek83} Hallin, M. \& Ingenbleek, J. F. (1983). Nonstationary Yule-Walker equations. \textit{Statist. Probab. Lett.} \textbf{1}, 189--195.

\bibitem {Hamdoune}\noindent Hamdoune, S. (1995). Etude des probl\`{e}mes d'estimation de
certains mod\`{e}les ARMA \'{e}volutifs. Thesis presented at
Universit\'{e} Henri Poincar\'{e}, Nancy 1.

%\bibitem {}
\bibitem {handeis} Hannan E. J. \& Deistler, M. (1988).
%\textit{The Statistical Theory of Linear Systems}, Wiley, New York. %ALM_V13
\textit{The statistical theory of linear systems}. Wiley, New York. %ALM_V13

\bibitem{}Hindrayanto, I., Koopman, S. J. \& Ooms, M. (2010). Exact maximum likelihood estimation for non-stationary periodic time series models. \emph{Comput. Statist. Data Anal.}~{\bf 54}, 2641--2654.


\bibitem{Janous_King} Janous, W. \& King, J. (2000). More on a sine product formula. \textit{The Mathematical Gazette} \textbf{84}, 113--115.  

\bibitem {Klimko_Nelson} Klimko, L. A. \& Nelson, P. I. (1978). On conditional least squares estimation for stochastic processes. \textit{Ann. Statist. }\textbf{6}, 629--642.

\bibitem {Kohn78} Kohn, R. (1978).
Asymptotic properties of time domain Gaussian estimators.
\textit{Adv. Appl. Probab.} \textbf{10}, 339--359.

\bibitem {Kollo_von_Rosen}  Kollo, T. \& von Rosen, D. (2005). \textit{Advanced multivariate statistics with matrices}. Springer Verlag, New York.% 579, Dordrecht, 

\bibitem {Kwoun_Yajima} Kwoun, G. H. \& Yajima, Y. (1986). On an autoregressive model with time-dependent coefficients. \textit{Ann. Inst. Statist. Math. Part A} \textbf{38}, 297--309.

\bibitem{Lutke} L\"utkepohl, H. (2005). \emph{New introduction to multiple time series analysis}. Springer-Verlag, New York.


\bibitem {Melard85} M\'{e}lard, G. (1985). \textit{Analyse de donn\'ees chronologiques}. Coll. S\'eminaire de math\'ematiques sup\'erieures - S\'eminaire scientifique OTAN (NATO Advanced Study Institute) \textbf{89}, Presses de l'Universit\'e de Montr\'eal, Montr\'eal.

%\bibitem {Priestley}Priestley, M. B. (1988). \textit{Non-Linear and Non-Stationary Time
%Series Analysis}, Academic Press, New York. %ALM_V13
\bibitem {Priestley}Priestley, M. B. (1988). \textit{Non-linear and non-stationary time
series analysis}. Academic Press, New York. %ALM_V13


%\bibitem {Quenouille}Quenouille, M. H. (1957). \textit{The} \textit{Analysis}
%\textit{of} \textit{Multiple} \textit{Time} \textit{Series,}
%Griffin, London. %ALM_V13
\bibitem {Quenouille}Quenouille, M. H. (1957). \textit{The} \textit{analysis}
\textit{of} \textit{multiple} \textit{time} \textit{series}. 
Griffin, London. %ALM_V13

\bibitem {Singh_Peiris}Singh, N. \& Peiris, M. S. (1987). A note on the properties of
some nonstationary ARMA processes. \textit{Stoch. Process.
Appl.} \textbf{24}, 151--155.

%\bibitem {Stout}Stout, W. F. (1974). \textit{Almost Sure Convergence}, Academic
%Press, New York. %ALM_V13
\bibitem {Stout}Stout, W. F. (1974). \textit{Almost sure convergence}. Academic
Press, New York. %ALM_V13

\bibitem {Subba_Rao70}\noindent Subba Rao, T. (1970). The fitting of non-stationary time-series
models with time dependent parameters. \textit{J.} \textit{Roy.}
\textit{Statist.} \textit{Soc.} \textit{Ser. B} \textbf{32},
312--322.


\bibitem {Taniguchi_Kakizawa} Taniguchi, M. \& Kakizawa, Y. (2000). \textit{Asymptotic theory of statistical inference for time series}. Springer Verlag, New York.

%\bibitem[Ter\"asvirta and Yang (2014)]{Tera_Yang_2014} Ter\"asvirta, T. and Yang, Y. (2014)  %ALM_V17
\bibitem[Ter\"asvirta \& Yang (2014)]{Tera_Yang_2014} Ter\"asvirta, T. \& Yang, Y. (2014)  %ALM_V17
Linearity and misspecification tests for vector smooth transition regression models, CREATES Research Papers; No. 2014-04, Aarhus Universitet.

%\bibitem[Ter\"asvirta et al.\ (2010)]{TeraTjosGrang} Ter\"asvirta, T., Tj\o stheim, D. and Granger C. W. J. (2010) \textit{Modelling Nonlinear Economic Time Series}. Oxford: Oxford University Press. %ALM_V17
\bibitem[Ter\"asvirta et al.\ (2010)]{TeraTjosGrang} Ter\"asvirta, T., Tj\o stheim, D. \& Granger C. W. J. (2010) \textit{Modelling nonlinear economic time series}. Oxford: Oxford University Press. %ALM_V17

\bibitem{Tiao_Grupe} Tiao, G. C. \& Grupe, M. R. (1980).  {Hidden periodic autoregressive-moving average models in time series data}. \emph{J. Roy. Statist. Soc. Ser. B} {\bf 67}, 365--373.

\bibitem {Tjostheim}Tj\o stheim, D. (1986). Estimation in nonlinear time series models.
\textit{Stoch. Process. Appl. }\textbf{21}, 251--273.

\bibitem{Try} Triantafyllopoulos, K. \& Nason, G. P. (2007). A Bayesian analysis of moving average processes with time-varying parameters. \emph{Comput. Statist. Data Anal.}~{\bf 52}, 1025--1046.


\bibitem {VB_Dahlhaus_2006}Van Bellegem, S. \& Dahlhaus, R. (2006). Semiparametric estimation
by model selection for locally stationary processes. 
\textit{J. Roy. Statist. Soc. Ser. B}
\textbf{68}, 721--746.

\bibitem {VB_von_Sachs} Van Bellegem, S. \& von Sachs, R. (2004).
Forecasting economic time series with unconditional time-varying variance.
\textit{International Journal of Forecasting} \textbf{20}, 611--627. 

\bibitem {Whittle}Whittle, P. (1965). Recursive relations for predictors of
non-stationary processes. \textit{J.} \textit{Roy.}
\textit{Statist.} \textit{Soc.} \textit{Ser.} \textit{B}
\textbf{27}, 523--532.

\end{thebibliography}
\end{document}